\documentclass{article} 
\usepackage[english]{babel} 
\usepackage{amssymb,amsmath,amsthm}
\usepackage{geometry} 
\usepackage{amsmath}
\usepackage{txfonts}
\usepackage{mathdots}
\usepackage{booktabs}
\usepackage{arydshln}
\usepackage[classicReIm]{kpfonts}
\usepackage{graphicx} 
\usepackage{enumitem}
\DeclareGraphicsExtensions{.pdf,.jpeg,.png}
\emergencystretch=\maxdimen
\usepackage[utf8]{inputenc}
\usepackage[english]{babel}
\usepackage{subfigure}
\usepackage{tablefootnote}
\usepackage{comment}
\usepackage{bm}
\usepackage{cite}
\usepackage[bookmarksopen,bookmarks=true]{hyperref}
\newtheorem{Th}{Theorem}
\begin{document}
\makeatletter
\renewcommand*{\@fnsymbol}[1]{\ensuremath{\ifcase#1\or \dagger\or *\or \dagger\or \ddagger\or
    \mathsection\or \mathparagraph\or \|\or **\or \dagger\dagger
    \or \ddagger\ddagger \else\@ctrerr\fi}}
\makeatother
\author{Iman Malmir\thanks{Department of Aerospace Engineering, MUT University of Technology, Tehran, Iran. E-mail: iman.malmir@outlook.com}}
\date{ }
\renewcommand*{\thefootnote}{\fnsymbol{footnote}}

\title{A novel wavelet-based optimal linear quadratic tracker for time-varying systems with multiple delays}
\maketitle
\begin{abstract}
\noindent A new method for solving optimal tracking control of linear quadratic time-varying systems with multiple time delays in state and input variables and with combined constraints is presented in this paper. By using the relations of Chebyshev wavelets, we simulate the optimal tracking problem to a static optimization one. This alternative method is applied on different optimal tracking systems and simulation results demonstrate the effectiveness of the proposed method.
\end{abstract}

\\


\section{ Introduction}

\renewcommand*{\thefootnote}{\fnsymbol{footnote}}
As we all know, a time-delay is a common phenomenon in engineering problems, and it is usually the main source of performance degradation in various control systems. The mathematical simulation of a time-delay system leads to a system of differential equations of delayed (retarded) type \cite{Bellman.Cooke}. One class of these equations, the integro-differential equations, was first studied by Volterra who developed his theories for them and investigated time delay phenomena in various systems \cite{Malek-Zavarei.Jamshidi}. The optimal control of linear quadratic time-varying systems with delays has been considered in many research works, see \cite{Malek-Zavarei.Jamshidi}, \cite{Gorecki} and the references therein. In optimal control problems, when we try to keep the output or state near a desired output or state, we are dealing with tracking problems. We find that in both state and output time-delay or regulator systems, the desired (reference) state and input is zero and in time-delay or regulator tracking system the error is to be made zero \cite{Naidu}. Optimal time-delay tracking control as a combination of time-delay optimal control and tracking control, aims at finding the optimal control law to minimize the given performance index function to make the system output track the reference signal in an optimal way and it has been a goal pursued in many areas. The tracking system is widely used in aerospace and mechanical systems, robot control, flight and spacecraft systems, etc. For example, consider an antenna control system to track an aircraft.

Most previous studies which have been done to solve the optimal time-delay tracking control problem, used discrete-time strategies. Ref.\cite{Leondes.Shieh} presented an iterative method by using discretization to find the suboptimal control of a linear quadratic time-varying system with multiple delays. \cite{Tsai.Tsai.Guo.Chen} presented a discretization approach by using the Newton center interpolation formula and the linear interpolation techniques for systems with multiple discrete and distributed time delays. In \cite{Tang.Sun.Liu} an optimal tracking controller for discrete time-delay systems based on a sensitivity approximation approach is designed in which the problem is transformed into a series of difference equations without time-advance on delayed terms. \cite{Tang.Sun} proposed a suboptimal tracking method which obtained by finite iterations of a solution of $N$ two-point boundary value problems (TPBVPs); for continuous-time control systems it provides good methodology but obtaining the solutions of these TPBVPs is difficult. In \cite{Zhang.Tang.Han} by applying an approximation approach of differential equations, TPBVP derived from the optimal tracking control problem is transformed into a sequence of linear TPBVPs without delays. Ref.\cite{Chang.Shieh.Liu.Cofie} used a discrete-time strategy to design an optimal controller for multiple-input and multiple-output (MIMO) continuous time systems with multiple delays in states, inputs and outputs which the Chebyshev quadrature formula together with a linear interpolation method was employed to get an extended discrete-time model from the continuous-time multiple time-delays system. \cite{Huang.Tsai.Provence.Shieh} converted a continuous-time input-state delayed system into an equivalent discrete-time input-state delayed model and its extended discrete-time delay-free model. A discrete-time methodology has notable disadvantages. For example, the discretized system is described by using the extended high-order state-space equation. The dimensions of the state-space description become significantly large when the sample period is extremely small compared to the time-delay. Computational difficulties with the discrete-time strategy occur in the algorithm, and the case we have to impose some constraints to the system.

In this paper we introduce an alternative numerical method to solve constrained linear quadratic time-delay tracking optimal control systems.  By choosing a parameterization method we convert the original problem to a static optimization one. This state and control parameterization method is based on Chebyshev wavelets which consist of Chebyshev polynomials of the first kind \cite{{Daubechies}, {Mason.Handscomb}}. The main motivations and contributions of the present research are summarized as
\begin{itemize}[label={}]
  \item A continuous-time accurate model of the optimal tracking control of linear quadratic time-varying systems with multiple delays is obtained. A major advantage of continuous-time models is that they avoid dependence on a particular timescale. Moreover, many standard numerical procedures are available to solve the simulated model.
  \item In the optimal tracking control, we would like on one hand, to keep the error small, but on the other hand, we must not pay higher cost to large inputs; hence, we have to try various values of the weighting matrices. From this fact, we conclude that we need a method provides good tracking in which with no concern about the algorithm of the solution we can change these matrices.
  \item Physical considerations imply that some constraints should be imposed on the optimal tracking control systems and unconstrained systems are less involved. An efficient method for solving a time-delay optimal tracking control problem should easily be able to resolve the problem in the cases we have to impose constraints to the system.
  \item The proposed method should have a good future and a high degree of flexibility. For example, it is possible that in the state equation, there is an inverse time term like $\mathbf{x}(t_f-t)$ and the method is capable of executing in this case.
  \item A method is presented that guarantees intersample constraint satisfaction and can be easily used to solve the optimal tracking problem in situations where there are multiple delays or no delays (LQT systems), the plant matrices are time-varying and/or constants.
  \item To handle final conditions and intersample constraints, the method can incorporate them directly into the model of the problem and unlike the other method we need no separate operations of applying these constraints to the obtained solutions.
  \item An optimal tracker is presented which can be applied to a time-delay system regardless of the system stability, minimum phase properties, the dimension of the system, equal number of input and output and the types of desired states and initial functions.
\end{itemize}

The rest of the paper is organized as follows. In Section 2, we describe the basic formulation of the Chebyshev wavelets required for converting the problem. Section 3 is devoted to the application of state-control parameterization via wavelets on the linear tracking delay system. Several numerical examples are simulated in section 4.

\subsection*{Some remarks on notation}
The transpose of a matrix $\mathbf{O}$ is written $\mathbf{O}^\top$.

\noindent $\mathbf{0}$ and $\mathbf{I}$ denote the zero and identity matrices, respectively.

\noindent The operator blkdiag denotes block diagonal concatenation of matrices.

\noindent $\otimes$ denotes Kronecker product. Kronecker product of a matrix $\mathbf{O}$ and the identity matrix $\mathbf{I}_q$ is denoted by $\hat{\mathbf{O}}$, Kronecker product of $\mathbf{O}$ and $\mathbf{I}_r$ is denoted by $\check{\mathbf{O}}$, that is, ${\hat{\mathbf{O}}={\mathbf{O}}\otimes \mathbf{I}_{q},\;\check{\mathbf{O}}={\mathbf{O}}\otimes \mathbf{I}_{r} }$.

\noindent $*$ indicates optimal condition.

\noindent $\mathcal{C}[0, t_{f}]$ denotes real-valued continuous functions on the closed interval $[0, t_{f}]$.

\section{Preliminaries}
\subsection{Chebyshev wavelets}
Chebyshev polynomials of the first kind of order $m$, $T_{m} (x)$ are solution of the differential equation
\[(1-x^2) T^{\prime\prime}_{m}-x T^{\prime}_{m}+m^{2}T_{m}=0\]
and form an orthogonal set on the interval $[-1, 1]$ with respect to the weight function ${w(x)=(1-x^{2} )^{-1/2}}$. The two useful relations for these polynomials are
\begin{equation} \label{1.0}
T_{m} (x)=\cos (m\arccos x), T_{m} (x)T_{m'} (x)=\tfrac{1}{2}\left\{T_{m+m'} (x)+T_{m-m'} (x)\right\}.
\end{equation} 
Chebyshev wavelets of the first kind are defined on [0,\,1] as
\begin{equation} \label{1.1}
\psi _{nm} (t)=\left\{\begin{array}{ll} {\sqrt{\frac{2^k}{\pi }}}\wp_{m} T_{m} (2^{k} t-2n+1), & t \in \left [\frac{n-1}{2^{k-1}} , \frac{n}{2^{k-1}} \right ]\vspace{1mm}\\  \;\;\;0, & t \notin \left [\frac{n-1}{2^{k-1}}, \frac{n}{2^{k-1}} \right ], \end{array}\right. 
\end{equation} 
 where
\begin{equation} \label{1.2}
\wp_{m}=\left\{\begin{array}{l} {\;\;1 ,\hspace{.32cm} m=0} \\  {\sqrt{2},\;\; m=1,2,3,\ldots,M-1} \end{array}\right.,\,n=1,2,\ldots,2^{k-1}
\end{equation}
and form an orthogonal basis with respect to the weight function ${w_{n} (t)}$, where ${w_{n} (t)=w(2^{k} t-2n+1)}$.
We can expand a function $f(t)$ in a series of Chebyshev wavelets by
\begin{equation} \label{1.3}
f(t)=\sum _{n=1}^{\infty }\, \sum _{m=0}^{\infty }f_{nm} \psi _{nm} (t)\cong \sum _{n=1}^{2^{k-1} }\sum _{m=0}^{M-1}f_{nm} \psi _{nm} (t) =\mathbf{f} \mathbf{\Psi}(t),
\end{equation}
where $\mathbf{f}$ and $\mathbf{\Psi}(t)$ are $1 \times 2^{k-1}M$ and $2^{k-1} M\times 1$ matrices and
\begin{equation} \label{1.4}
\mathbf{f}=[f_{10}, \ldots, f_{1M-1}, f_{20}, \ldots, f_{2M-1}, \ldots, f_{2^{k-1} 0}, \ldots, f_{2^{k-1} M-1}],
\end{equation}
\begin{equation} \label{1.5}
\mathbf{\Psi}(t)=\left[\psi _{10}(t), \ldots, \psi _{1M-1} (t), \psi _{20} (t), \ldots, \psi _{2M-1} (t), \ldots, \psi _{2^{k-1} 0} (t), \ldots, \psi _{2^{k-1} M-1} (t)\right]^{\top}.
\end{equation}
The coefficients of Chebyshev scaling functions can be approximated as follows \cite{iman}
\begin{equation} \label{1.6}
f_{nm} =\frac{\wp_{m}}{\sqrt{2^{k}\pi}} \int _{0}^{\pi }f(\frac{\cos \theta +2n-1}{2^{k} }  )\, \cos m\theta \, d\theta.
\end{equation}
\begin{Th}
\textbf{(Convergence of Chebyshev wavelets expansion)} A twice differentiable function $f(t)$, defined on $[0, 1]$, with bounded second derivatives, say $|f''(t)| \le \rho$, can be expanded as an infinite sum of Chebyshev wavelets, and
this series converges uniformly to $f(t)$.
\end{Th}

\begin{proof}
For $m \ge 2$  according to eq.\eqref{1.6} we get
\[
f_{nm} =\frac{1}{\sqrt{2^{5k+1}\pi} m} \int _{0}^{\pi }f''\left (\frac{\cos \theta +2n-1}{2^{k}}\right )\left (\frac{\sin(m-1)\theta}{m-1}-\frac{\sin(m+1)\theta}{m+1}\right )\sin\theta\,\,d\theta,
\]
where integration by parts was used twice in this evaluation. Therefore for some $\varepsilon \in [0,\; 1]$ we find
\[
\begin{array}{l}{|f_{nm}|= \frac{1}{\sqrt{2^{5k+1}\pi} m} \left|\int _{0}^{\pi }f''\left (\frac{\cos \theta +2n-1}{2^{k}}\right )\left (\frac{\sin(m-1)\theta}{m-1}-\frac{\sin(m+1)\theta}{m+1}\right )\sin\theta\,\,d\theta\right|}\\{\hspace{.76cm} \le \frac{1}{\sqrt{2^{5k+1}\pi} m}\left|f^{\prime \prime}(\varepsilon)\right| \left (\frac{1}{m-1}+\frac{1}{m+1}\right )\pi  }\\{\hspace{.76cm} \le \sqrt{\frac{\pi}{2^{5k-1}}}\left (\frac{\rho}{m^2-1}\right).}\end{array}
\]
For $m=1$, from \eqref{1.6} by employing integration by parts we simply have 
\[
f_{n1}=\frac{1}{\sqrt{2^{3k-1}\pi}}\int _{0}^{\pi }f'(\frac{\cos \theta +2n-1}{2^{k}})\, \sin^{2}\theta \, d\theta.
\]
So that
\[
|f_{n1}| = \frac{1}{\sqrt{2^{3k-1}\pi}}\left| \int _{0}^{\pi }f'(\frac{\cos \theta +2n-1}{2^{k}})\, \sin^{2}\theta \, d\theta \right|.
\]
Then by choosing $t_{0} \in [0, 1]$ we can write $f^{\prime}(t)=f^{\prime}(t_{0})+\int_{t_{0}}^{t}f^{\prime\prime}(x)dx$, $|t-t_{0}| \le 1$; thus for all $t \in [0, 1]$
\[
\begin{array}{l}{\left|f^{\prime}(t)\right| \le \left|f^{\prime}(t_{0})\right|+\left|\int_{t_{0}}^{t}f^{\prime\prime}(x)dx\right|} \\ {\hspace{9.2mm}\le \left|f^{\prime}(t_{0})\right|+\rho |t-t_{0}|} \\ {\hspace{9.2mm}\le \left|f^{\prime}(t_{0})\right|+\rho.} \end{array}
\]
Assuming $|f'(t)| \le \rho_{1}$ yields immediately
\[
|f_{n1}| \le \sqrt{\tfrac{\pi }{2^{3k-1}}} \rho_{1}.
\]
For $m=0$ by the similar procedure as explained above, it is easy to verify that $f(t)$ is bounded; let us assume here that $|f(t)| \le \rho_{0}$. Hence by \eqref{1.6} we find $|f_{n0}| \le \sqrt{\tfrac{\pi}{2^{k}}}\rho_{0}$.\\
It is readily seen that $\left|\psi _{n0}(t) \right| = \sqrt{2^{k}}/\sqrt{\pi}$ and $\left|\psi _{n1}(t) \right|, \left|\psi _{nm}(t) \right| \le \sqrt{2^{k+1}}/\sqrt{\pi}$; we can then write
\[
\left|f_{n0}\right| \left|\psi_{n0}(t)\right| \le \rho_{0}, \left|f_{n1}\right| \left|\psi_{n1}(t)\right| \le \frac{\rho_{1}}{2^{k-1}}\; \text{and} \; \left|f_{nm}\right| \left|\psi_{nm}(t)\right| \le \frac{\rho}{2^{2k-1}\left (m^2-1\right)}, m \ge 2.
\]
Now from \eqref{1.1} we deduce $f(t)$ is in the form of a piecewise-defined function which we have $f(t)=f_{n}(t)$ on each subinterval $\left [\tfrac{n-1}{2^{k-1}}, \tfrac{n}{2^{k-1}} \right ]$, where $f_{n}(t)=\sum _{m=0}^{\infty }f_{nm} \psi _{nm} (t),\;t \in \left [\tfrac{n-1}{2^{k-1}}, \tfrac{n}{2^{k-1}} \right ]$. Consequently
\[
\begin{array}{l}{|f_{n}(t)|= \left |f_{n0} \psi_{n0}(t) +f_{n1} \psi_{n1}(t) +\sum _{m=2}^{\infty} f_{nm} \psi _{nm} (t) \right |} \\ {\hspace{9mm} \le \rho_{0} + \frac{\rho_{1}}{2^{k-1}} + \frac{\rho}{2^{2k-1}} \sum _{m=2}^{\infty} \frac{1}{m^2-1} = \rho_{0} + \frac{\rho_{1}}{2^{k-1}} + \frac{3}{4}} \frac{\rho}{2^{2k-1}}.\end{array}
\]
Thus, $\sum _{m=0}^{\infty }f_{nm} \psi _{nm}(t)$ is absolutely convergent on all subintervals; it means that $\sum _{n=1}^{\infty }\, \sum _{m=0}^{\infty }f_{nm} \psi _{nm} (t)$ converges to $f (t)$ uniformly and this completes the proof. Furthermore since $k \ge 2$, we conclude that
\[
|f_{n0}| \le \tfrac{\sqrt{\pi}\rho_{0}}{2}, |f_{n1}| \le \tfrac{\sqrt{2\pi}\rho_{1}}{8} \; \text{and} \; |f_{nm}| \le \tfrac{\sqrt{2\pi}\rho}{96}, m \ge 2.
\]
\end{proof}

\subsection{The operational matrix of integration for Chebyshev wavelets}
 The integration of the Chebyshev wavelet vector defined in eq. \eqref{1.5} on $[0,t ]$ can be obtained as
\begin{equation} \label{2.1} 
\int _{0}^{t}\mathbf{\Psi}(\epsilon)d\epsilon \cong \mathbf{P}\mathbf{\Psi}(t).
\end{equation} 
The matrix $\mathbf{P}$ is called $2^{k-1} M\times 2^{k-1} M$ Chebyshev wavelets operational matrix of integration. It follows from \eqref{1.5} that
$
\int _{0}^{t}\mathbf{\Psi}(\epsilon)d\epsilon=\int _{\tfrac{n-1}{2^{k-1}}}^{t}[\psi_{nm}(\epsilon)]^{\top} d\epsilon
$.
We conclude from \eqref{1.1} that when $t < n\,/2^{k-1}$, the integral is a function of the time in which the integrand defined, so it should be expanded by the wavelets of the current subinterval; when $t = n\,/2^{k-1}$, this definite integral should be expanded on all subsequent subintervals. Hence by setting $t_{n}=2^{k}t-2n+1$, we can write
\[
\text{if}\;m=0, \left\{\begin{array}{l} {\int_{\tfrac{n-1}{2^{k-1}}}^{t}\psi_{n0}(\epsilon)d\epsilon=\sqrt{\tfrac{2^k}{\pi}}\left(\tfrac{1}{2^k}T_{0}+\tfrac{1}{2^k}T_{1}(t_n)\right)=\tfrac{1}{2^k}[1,\tfrac{1}{\sqrt{2}},\underbrace{0,0,\ldots,0}_{M-2}]\bm{\varphi}_{n}^{\top}(t),\hspace{1.28cm}t < \tfrac{n}{2^{k-1}}}\\{\int_{\tfrac{n-1}{2^{k-1}}}^{t} \psi_{n0}(\epsilon)d\epsilon=\sqrt{\tfrac{2^k}{\pi}}\tfrac{2}{2^k}\sqrt{\tfrac{\pi}{2^k}}\sum_{\eta=n+1}^{2^{k-1}}\psi_{\eta 0}(t)=\tfrac{1}{2^k}\sum_{\eta=n+1}^{2^{k-1}}[2,\overbrace{0,0,\ldots,0}^{M-1}]\bm{\varphi}_{\eta}^{\top}(t),\; t = \tfrac{n}{2^{k-1}},} \end{array}\right.
\]
\[
\text{if}\;m=1, \left\{\begin{array}{l} {\int_{\tfrac{n-1}{2^{k-1}}}^{t}\psi_{n1}(\epsilon)d\epsilon=\sqrt{\tfrac{2^{k+1}}{\pi}}\tfrac{1}{2^k}\tfrac{1}{4}\left(T_{2}(t_n)-1\right)=\tfrac{1}{2^k}[-\tfrac{\sqrt{2}}{4},0,\tfrac{1}{4},\underbrace{0,0,\ldots,0}_{M-3}]\bm{\varphi}_{n}^{\top}(t),\hspace{.8cm}t < \tfrac{n}{2^{k-1}}}\vspace{-3mm}\\{\int_{\tfrac{n-1}{2^{k-1}}}^{t} \psi_{n1}(\epsilon)d\epsilon=0=\tfrac{1}{2^k}\sum_{\eta=n+1}^{2^{k-1}}[0,\overbrace{0,0,\ldots,0}^{M-1}]\bm{\varphi}_{\eta}^{\top}(t), \hspace{3.84cm}t = \tfrac{n}{2^{k-1}},}\end{array}\right.
\]
\[
\hspace{1.8mm}\text{if}\;m\ge2, \left\{\begin{array}{l} {\begin{array}{l}{\hspace{-1.9mm}\int_{\tfrac{n-1}{2^{k-1}}}^{t}\psi_{nm}(\epsilon)d\epsilon=\sqrt{\tfrac{2^{k+1}}{\pi}}\tfrac{1}{2^k}\tfrac{1}{2}\left(\tfrac{1}{m+1}T_{m+1}(t_n)-\tfrac{1}{m-1}T_{m-1}(t_n)-\tfrac{(-1)^{m+1}}{m+1}+\tfrac{(-1)^{m-1}}{m-1}\right)}\\{\hspace{2.12cm}=\tfrac{1}{2^k}[\underbrace{\tfrac{(-1)^{m-1}\sqrt{2}}{m^{2}-1}, \ldots, -\tfrac{1}{2(m-1)},0,\tfrac{1}{2(m+1)}, \ldots,0}_M]\bm{\varphi}_{n}^{\top}(t),\hspace{1.82cm}t < \tfrac{n}{2^{k-1}}}\end{array}}\\{\int_{\tfrac{n-1}{2^{k-1}}}^{t} \psi_{nm}(\epsilon)d\epsilon=\sqrt{\tfrac{2^{k+1}}{\pi}}\tfrac{1}{2^k}\left(-\tfrac{1+(-1)^{m}}{m^{2}-1}\right)\sqrt{\tfrac{\pi}{2^k}}\sum_{\eta=n+1}^{2^{k-1}}\psi_{\eta 0}(t)}\\{\hspace{2.30cm}=\tfrac{1}{2^k}\sum_{\eta=n+1}^{2^{k-1}}[-\tfrac{(1+(-1)^{m})\sqrt{2}}{m^{2}-1},\overbrace{0,0,\ldots,0}^{M-1}]\bm{\varphi}_{\eta}^{\top}(t),\hspace{2.7cm} t = \tfrac{n}{2^{k-1}},}\end{array}\right.
\]
where $\bm{\varphi}_{\kappa}(t)=[\psi _{\kappa 0}(t), \psi _{\kappa 1}(t), \ldots , \psi _{\kappa M-1} (t)]$ for $\kappa =n, n+1,n+2,\ldots,2^{k-1}$. As a result
\begin{equation} \label{2.2} 
\mathbf{P}=\frac{1}{2^{k} } \left[\begin{array}{ccc} \mathbf{p}_{1} \\ \mathbf{p}_{2} \\ \vdots \\ \mathbf{p}_{2^{k-1}} \end{array}\right],\;\mathbf{p}_{n}= \big[ \overbrace{\mathbf{0} \quad \mathbf{0} \quad \cdots \quad \mathbf{0}}^{(n-1)\;\text{times}} \quad \mathbf{L} \quad \overbrace{\mathbf{E} \quad \mathbf{E} \quad \cdots \quad  \mathbf{E}}^{(2^{k-1}-n)\;\text{times}} \big],
\end{equation} 
where $\mathbf{L}$ and $\mathbf{E}$ are $M\times M$ matrices and are in the forms ($m \ge 2$)

\begin{equation} \label{2.3} 
\mathbf{L}=\left[\begin{smallmatrix} {1} & {{\textstyle\frac{1}{\sqrt{2} }} } & {0} & {0} & {0} & {\cdots } & {0} & {0} & {0} \\ {-{\textstyle\frac{\sqrt{2} }{4}} } & {0} & {{\textstyle\frac{1}{4}} } & {0} & {0} & {\cdots } & {0} & {0} & {0} \\ {-{\textstyle\frac{\sqrt{2} }{3}} } & {-{\textstyle\frac{1}{2}} } & {0} & {{\textstyle\frac{1}{6}} } & {0} & {\cdots } & {0} & {0} & {0} \\ {\tfrac{\sqrt{2}}{8}} & {0} & {-{\textstyle\frac{1}{4}} } & {0} & {\tfrac{1}{8}} & {\cdots } & {0} & {0} & {0} \\ {\vdots} & {\vdots } & {\vdots } & {\vdots } & {\vdots} & {\ddots } & {\vdots} & {\vdots} & {\vdots} \\ {\tfrac{(-1)^{M-1}\sqrt{2}}{(M-2)^2-1}} & {0} & {0} & {0} & {0} & {\cdots } & {-\tfrac{1}{2(M-3)}} & {0} & {{\textstyle\frac{1}{2(M-1)}} } \\ {\tfrac{(-1)^{M}\sqrt{2}}{(M-1)^2-1}} & {0} & {0} & {0} & {0} & {\cdots } & {0} & {-{\textstyle\frac{1}{2(M-2)}} } & {0} \end{smallmatrix}\right]
, 
\mathbf{E}=\left[\begin{smallmatrix} {2} & {0} & {0} & {\cdots } & {0} \\ {0} & {0} & {0} & {\cdots } & {0} \\ -\tfrac{2\sqrt{2}}{3} & {0} & {0} & {\cdots } & {0} \\ {0} & {0} & {0} & {\cdots } & {0} \\ {\vdots } & {\vdots } & {\vdots } & {\ddots } & {\vdots } \\ {-\tfrac{(1+(-1)^{m})\sqrt{2}}{m^{2}-1}} & {0} & {0} & {\cdots } & {0} \\ {\vdots } & {\vdots } & {\vdots } & {\ddots } & {\vdots } \\ {-\tfrac{(1+(-1)^{M-1})\sqrt{2}}{(M-1)^2-1}} & {0} & {0} & {\cdots } & {0} \end{smallmatrix}\right].
\end{equation}

\subsection{The integration matrix of the product of Chebyshev wavelets on $\mathbf{[0,1]}$}
To transform the performance index into a quadratic form, we have to find an integration matrix of the product of two Chebyshev scaling function vectors on [0, 1], so we introduce
\begin{equation} \label{2.5} 
\mathbf{C}=\int _{0}^{1}\mathbf{\Psi}(t)\mathbf{\Psi}^{\top}(t)dt.
\end{equation} 
$\mathbf{C}$ is obtained by integrating the elements of $\mathbf{\Psi}(t)\mathbf{\Psi}^{\top}(t)$ from 0 to 1 and by using the compact support property of the wavelets: $ \forall n\ne n', \psi _{nm} (t)\psi _{n'm'} (t)=0$, where $m'=0,1,2,\ldots,M-1,\;n'=1,2,3,\ldots,2^{k-1}$. 
We see immediately that
\begin{align*}\mathbf{C}=&\int _{0}^{1}[\psi _{10}(t), \ldots, \psi _{1M-1} (t), \ldots, \psi _{2^{k-1} M-1} (t)]^{\top}[\psi _{10}(t),\ldots,\psi _{1M-1} (t), \ldots, \psi _{2^{k-1} M-1}(t)]dt\\=&\left[\begin{smallmatrix} \mathbf{C}_1 & \mathbf{0}_{M \times M} & \cdots & \mathbf{0}_{M \times M} \\  \mathbf{0}_{M \times M}  & \mathbf{C}_2 & \cdots & \mathbf{0}_{M \times M} \\ \vdots & \vdots & \ddots & \vdots \\ \mathbf{0}_{M \times M} & \mathbf{0}_{M \times M} & \cdots & \mathbf{C}_{2^{k-1}} \end{smallmatrix}\right],\end{align*}
where $\forall m,m' \ge 2$
\[\mathbf{C}_{n}=\int_{\tfrac{n-1}{2^{k-1}}}^{\tfrac{n}{2^{k-1}}}\left[ \begin{smallmatrix} \psi_{n0}(t)\psi_{n0}(t) & \psi_{n0}(t)\psi_{n1}(t) & \ldots &\psi_{n0}(t)\psi_{nm'}(t) & \ldots & \psi_{n0}(t)\psi_{nM-1}(t) \\  \psi_{n1}(t)\psi_{n0}(t)  & \psi_{n1}(t)\psi_{n1}(t) & \ldots & \psi_{n1}(t)\psi_{nm'}(t) & \cdots & \psi_{n1}(t)\psi_{nM-1}(t) \\ \vdots & \vdots & \ddots & \vdots & \ddots & \vdots \\ \psi_{nm}(t)\psi_{n0}(t) & \psi_{nm}(t)\psi_{n1}(t) & \cdots & \psi_{nm}(t)\psi_{nm'}(t) & \cdots & \psi_{nm}(t)\psi_{nM-1}(t) \\ \vdots & \vdots & \ddots & \vdots & \ddots & \vdots \\ \psi_{nM-1}(t)\psi_{n0}(t) & \psi_{nM-1}(t)\psi_{n1}(t) & \cdots & \psi_{nM-1}(t)\psi_{nm'}(t) & \cdots & \psi_{nM-1}(t)\psi_{nM-1}(t) \end{smallmatrix} \right]dt.
\]
Since $\frac{n-1}{2^{k-1} } \le t\le \frac{n}{2^{k-1} } $, it follows that $-1\le 2^{k} t-2n+1\le 1$. So we substitute $\cos \theta =2^{k} t-2n+1$. From definition of these wavelets we can identify $\mathbf{C}_{1}=\mathbf{C}_{2}=\cdots=\mathbf{C}_{2^{k-1}}$, 
hence $\mathbf{C}$ is symmetric. In general, we find
\[
\begin{array}{l}{\int_{\tfrac{n-1}{2^{k-1}}}^{\tfrac{n}{2^{k-1}}}\psi _{nm} (t)\psi _{nm'} (t)dt=\left\{\begin{array}{l} {\frac{2^k}{\pi}\int_{\tfrac{n-1}{2^{k-1}}}^{\tfrac{n}{2^{k-1}}}dt,\hspace{5.12cm}m=m'=0}\vspace{1mm} \\ {\frac{2^{k+1}}{\pi }\int_{\tfrac{n-1}{2^{k-1}}}^{\tfrac{n}{2^{k-1}}}T_{m}(2^{k}t-2n+1)T_{m'}(2^{k}t-2n+1)dt, m, m'\ne 0}\vspace{1mm} \\ {\frac{\sqrt{2}2^k}{\pi}\int_{\tfrac{n-1}{2^{k-1}}}^{\tfrac{n}{2^{k-1}}}T_{m'}(2^{k}t-2n+1)dt,\hspace{2.26cm} m=0} \end{array}\right.}\\ {\hspace{3.37cm}=\left\{\begin{array}{l} {\frac{2}{\pi},\hspace{4.1cm}m=m'=0}\\ { \frac{2}{\pi }\int_{0}^{\pi }\cos m\theta \cos m'\theta \sin \theta \, d\theta,\;\; m, m'\ne 0}\\ {0,\hspace{4.16cm} m+m' \;\;\text{is odd}}\\ {\frac{-2\sqrt{2}}{\pi(m'^{2}-1)},\hspace{3.20cm} m=0.}  \end{array}\right.} \end{array}
\]
When $m+m'$ is even, we get
\[{\int_{0}^{\pi }\cos m\theta \cos m'\theta \sin \theta \, d\theta}=\tfrac{1}{2}\left(\tfrac{1}{m+m'+1}-\tfrac{1}{m+m'-1}+\tfrac{1}{m-m'+1} -\tfrac{1}{m-m'-1}\right).\]  
Finally, the matrix $\mathbf{C}$ is
\begin{equation} \label{2.6} 
\mathbf{C}=\tfrac{2}{\pi}.\text{blkdiag}\left(\mathbf{C}_{1},\mathbf{C}_{2},\cdots,\mathbf{C}_{2^{k-1}} \right),
\end{equation} 
where for $\mathbf{C}_{n}=[\text{c}_{ij}]$, $i, j=1,2,\ldots,M$ we have
\[\text{c}_{ij}=\left\{\begin{array}{l} {l_{ij} \frac{(1-(i-1)^{2} - (j-1)^{2})}{((i+j-2)^{2}-1)((i-j)^{2}-1)} \, \, ,\, \, \, i+j = 2, 4,6,\cdots , 2M} \\ {0,\hspace{3.03cm}\, i+j = 3, 5,7,\cdots , 2M-1} \end{array}\right.\text{and}\,\; l_{ij}=\left\{\begin{array}{l} {1 ,\hspace{.21cm} i=j=1} \\ {\hspace{-2mm}\sqrt{2},\, i=1\,\,\, \text{or}\, \, j=1} \\ {2,\hspace{.21cm} i, j \ne 1.} \end{array}\right. \] 

\subsection{The product operational matrix of Chebyshev wavelets}
The useful property of the product of two Chebyshev wavelets vectors is
\begin{equation} \label{2.8} 
\mathbf{f} \mathbf{\Psi}(t) {\mathbf{\Psi}}^{\top}(t) \cong {\mathbf{\Psi}}^{\top}(t) \tilde{\mathbf{f}},
\end{equation} 
where $\tilde{\mathbf{f}}$ is called the $2^{k-1} M\times 2^{k-1} M$ product operational matrix. As we see in the previous section, for $n\ne n'$ we have $\psi _{nm} (t)\psi _{n'm'} (t)=0$, thus
\[\begin{array}{l}{\mathbf{f}\mathbf{\Psi}(t) {\mathbf{\Psi}}^{\top}(t)=\mathbf{f}[\psi _{10}(t),\,\ldots,\,\psi _{nm} (t),\,\ldots,\, \psi _{2^{k-1} M-1} (t)]^{\top}[\psi _{10}(t),\,\ldots ,\,\psi _{n'm'} (t),\,\ldots,\, \psi _{2^{k-1} M-1}(t)]}\\{\hspace{1.89cm}=[f_{10},\, \ldots,\, f_{nm},\, \ldots,\, f_{2^{k-1} M-1}]\begin{bmatrix} \bm{\Omega}_1 & \mathbf{0}_{M \times M} & \cdots & \mathbf{0}_{M \times M} \\  \mathbf{0}_{M \times M}  & \bm{\Omega}_2 & \cdots & \mathbf{0}_{M \times M} \\ \vdots & \vdots & \ddots & \vdots\\ \mathbf{0}_{M \times M} & \mathbf{0}_{M \times M} & \cdots & \bm{\Omega}_{2^{k-1}} \end{bmatrix}.}\end{array}\]
If $m \ne m'$ and $m, m' \ne0$, by \eqref{1.0} we see that
\[\begin{array}{l} {\psi _{nm} (t)\, \psi _{nm'} (t)={\textstyle\frac{2^{k+1} }{\pi }} T_{m} (2^{k} t-2n+1)T_{m'} (2^{k} t-2n+1)} \\ {\hspace{2.2cm}=\tfrac{2^{k/2} }{\sqrt{\pi }}({\textstyle\frac{2^{k/2} }{\sqrt{\pi } }} T_{m-m'} (2^{k} t-2n+1)+\tfrac{2^{k/2} }{\sqrt{\pi}} T_{m+m'} (2^{k} t-2n+1))} \\ {\hspace{2.2cm}=\sqrt{\tfrac{2^{k} }{\pi}}({\textstyle\frac{1}{\sqrt{2} }} \psi _{n\,m-m'} (t)+{\textstyle\frac{1}{\sqrt{2} }} \psi _{n\,m+m'} (t))\,,\; m+m'\le M-1}\end{array}\]
and when $m+m'\ge M$ then $\psi _{n\,m+m'} (t)$ is not defind and hence
\[\psi _{nm} (t)\, \psi _{nm'} (t)\approx \sqrt{\tfrac{2^{k} }{\pi}} {\textstyle\frac{1}{\sqrt{2} }} \psi _{n\,m-m'} (t)\, ,\, \, m+m'\ge M.\]
If $m=m'$ we obtain
\[
\psi_{nm} (t) \psi_{nm'} (t)=\left\{\begin{array}{l}{\sqrt{\tfrac{2^{k} }{\pi}}(\psi _{n0}(t)+\tfrac{1}{\sqrt{2}}\psi _{n\,m+m'} (t)),\; m+m' \le M-1} \\
{\sqrt{\tfrac{2^{k} }{\pi}}\psi _{n0}(t),\hspace{26.2mm}m+m' \ge M.}\end{array}\right.
\]
Finally, if $m\;\text{or}\;m'=0$
\[\psi_{n0} (t)\psi _{nm'} (t)=\sqrt{\tfrac{2^{k} }{\pi}} \psi _{nm'} (t),\;{\text{or}}\;\psi_{nm} (t)\psi _{n0} (t)=\sqrt{\tfrac{2^{k} }{\pi}} \psi _{nm} (t).\]
By applying these relations to the product $\mathbf{\Psi}(t)\mathbf{\Psi}^{\top}(t)$, also assuming that $2<\varsigma<M-2, \varsigma \in \mathbb{N}$, we get
\[\bm{\Omega}_{n}=\sqrt{\tfrac{2^{k} }{\pi}}\left[\begin{smallmatrix} \psi_{n0}(t) & \psi_{n1}(t) & \ldots & \psi_{n\varsigma}(t) & \psi_{n\varsigma+1}(t) & \ldots & \psi_{nM-1}(t) \\ \psi_{n1}(t) & \psi_{n0}(t)+\tfrac{1}{\sqrt{2}}\psi_{n2}(t) & \ldots & \tfrac{1}{\sqrt{2}}(\psi_{n\varsigma-1}(t)+\psi_{n\varsigma+1}(t)) & \tfrac{1}{\sqrt{2}}(\psi_{n\varsigma}(t)+\psi_{n\varsigma+2}(t)) & \cdots & \tfrac{1}{\sqrt{2}}\psi_{nM-2}(t) \\ \vdots & \vdots & \ddots & \vdots & \vdots & \ddots & \vdots \\ \psi_{n\varsigma}(t) & \tfrac{1}{\sqrt{2}}(\psi_{n\varsigma-1}(t)+\psi_{n\varsigma+1}(t)) & \cdots & \psi_{n0}(t)+\tfrac{1}{\sqrt{2}}\psi_{n2\varsigma}(t) & \tfrac{1}{\sqrt{2}}(\psi_{n1}(t)+\psi_{n2\varsigma+1}(t)) & \cdots & \tfrac{1}{\sqrt{2}}\psi_{nM-\varsigma-1}(t)\\ \psi_{n\varsigma+1}(t) & \tfrac{1}{\sqrt{2}}(\psi_{n\varsigma}(t)+\psi_{n\varsigma+2}(t)) & \cdots & \tfrac{1}{\sqrt{2}}(\psi_{n1}(t)+\psi_{n2\varsigma+1}(t)) & \psi_{n0}(t) & \cdots & \tfrac{1}{\sqrt{2}}\psi_{nM-\varsigma-2}(t) \\ \vdots & \vdots & \ddots & \vdots & \vdots & \ddots & \vdots \\ \psi_{nM-1}(t) & \tfrac{1}{\sqrt{2}}\psi_{nM-2}(t) & \cdots & \tfrac{1}{\sqrt{2}}\psi_{nM-\varsigma-1}(t) & \tfrac{1}{\sqrt{2}}\psi_{nM-\varsigma-2}(t) & \cdots & \psi_{n0}(t) \end{smallmatrix}\right].\]
When $M$ is even, $M-1$ is odd and we have $2\varsigma=M-2$. So $\varsigma=(M-2)/2, M-\varsigma-1=M/2$ and $2\varsigma+1=M-1$. When $M$ is odd, we have $2\varsigma=M-1$. Thus $\varsigma=(M-1)/2, M-\varsigma-1=(M-1)/2$ . Since $2\varsigma+1=M$, then necessarily $\psi_{n2\varsigma+1}=0$. If we take $\tilde{\mathbf{f}}=[\tilde{f}_{ab}]$, where $a, b=1,2,\ldots,2^{k-1}M$, then by using the elements of the matrices $\bm{\Omega}_{n}$ and equating coefficients of same Chebyshev wavelets, the matrix $\tilde{\mathbf{f}}$ takes the form
\begin{equation} \label{2.9} 
\tilde{\mathbf{f}}=\sqrt{\tfrac{2^{k} }{\pi}}.\text{blkdiag}\left(\tilde{\mathbf{f}}_{1} ,\tilde{\mathbf{f}}_{2},\, \cdots,\, \tilde{\mathbf{f}}_{2^{k-1}}\right),
\end{equation}
where by setting
\[
\zeta=M-\varsigma-1\
\]
and
\[
\xi=2\varsigma+1
\]
we have
\[
\tilde{\mathbf{f}}_{n} =\tiny  \left[\begin{smallmatrix} {f_{n0} } & {f_{n1} } & {f_{n2} } & {\cdots } & {f_{n\varsigma}} & {f_{n\varsigma+1}} & {\cdots } & {f_{nM-3}} & {f_{nM-2}} & {f_{nM-1} } \\ {f_{n1} } & {f_{n0} +\tfrac{1}{\sqrt{2}} f_{n2} } & {\tfrac{1}{\sqrt{2}} (f_{n1} + f_{n3} )} & {\cdots } & {\tfrac{1}{\sqrt{2}} (f_{n\varsigma-1}+f_{n\varsigma+1})} & {\tfrac{1}{\sqrt{2}} (f_{n\varsigma} +f_{n\varsigma+2})} & {\cdots } & {\tfrac{1}{\sqrt{2}} (f_{nM-4}+f_{nM-2})} & {\tfrac{1}{\sqrt{2}} (f_{nM-3} +f_{nM-1})} & {\tfrac{1}{\sqrt{2}} f_{nM-2} } \\ {f_{n2} } & {\tfrac{1}{\sqrt{2}} (f_{n1} +f_{n3} )} & {f_{n0} +\tfrac{1}{\sqrt{2}} f_{n4} } & {\cdots } & {\tfrac{1}{\sqrt{2}} (f_{n\varsigma-2}+f_{\varsigma+2})} & {\tfrac{1}{\sqrt{2}} (f_{n\varsigma-1}+f_{\varsigma+3})} & {\cdots } & {\tfrac{1}{\sqrt{2}} (f_{nM-5} +f_{nM-1})} & {\tfrac{1}{\sqrt{2}} f_{nM-4}} & {\tfrac{1}{\sqrt{2}} f_{nM-3} } \\ {\vdots } & {\vdots } & {\vdots } & {\ddots } & {\vdots } & {\vdots} & {{\mathinner{\mkern2mu\raise1pt\hbox{.}\mkern2mu\raise4pt\hbox{.}\mkern2mu\raise7pt\hbox{.}\mkern1mu}} } & {\vdots } & {\vdots } & {\vdots } \\ {f_{n\varsigma}} & {\tfrac{1}{\sqrt{2}} (f_{n\varsigma-1}+f_{n\varsigma+1})} & {\tfrac{1}{\sqrt{2}} (f_{n\varsigma-2}+f_{\varsigma+2})} & {\cdots } & {f_{n0} +\tfrac{1}{\sqrt{2}} f_{n 2\varsigma} } & {\tfrac{1}{\sqrt{2}} (f_{n1} + f_{n\xi } )} & {\cdots } & {\tfrac{1}{\sqrt{2}} f_{n\zeta-2}} & {\tfrac{1}{\sqrt{2}} f_{n\zeta-1}} & {\tfrac{1}{\sqrt{2}} f_{n\zeta}} \\ {f_{n\varsigma+1}} & {\tfrac{1}{\sqrt{2}} (f_{n\varsigma} +f_{n\varsigma+2})} & {\tfrac{1}{\sqrt{2}} (f_{n\varsigma-1}+f_{n\varsigma+3})} & {\cdots } & {\tfrac{1}{\sqrt{2}} (f_{n1} +f_{n\xi } )} & {f_{n0} } & {\cdots } & {\tfrac{1}{\sqrt{2}} f_{n\zeta-3}} & {\tfrac{1}{\sqrt{2}} f_{n\zeta-2}} & {\tfrac{1}{\sqrt{2}} f_{n\zeta-1}} \\ {f_{n\varsigma+2}} & {\tfrac{1}{\sqrt{2}} (f_{n\varsigma+1}+f_{n\varsigma+3})} & {\tfrac{1}{\sqrt{2}} (f_{n\varsigma}+f_{n\varsigma+4})} & {\cdots } & {\tfrac{1}{\sqrt{2}} f_{n2} } & {\tfrac{1}{\sqrt{2}} f_{n1} } & {\cdots } & {\tfrac{1}{\sqrt{2}} f_{n\zeta-4}} & {\tfrac{1}{\sqrt{2}} f_{n\zeta-3}} & {\tfrac{1}{\sqrt{2}} f_{n\zeta-2}} \\ {\vdots } & {\vdots } & {\vdots } & {{\mathinner{\mkern2mu\raise1pt\hbox{.}\mkern2mu\raise4pt\hbox{.}\mkern2mu\raise7pt\hbox{.}\mkern1mu}} } & {\vdots } & {\vdots } & {\ddots } & {\vdots } & {\vdots } & {\vdots } \\ {f_{nM-2} } & {\tfrac{1}{\sqrt{2}} (f_{nM-3} +f_{nM-1} )} & {\tfrac{1}{\sqrt{2}} f_{nM-4} } & {\cdots } & {\tfrac{1}{\sqrt{2}} f_{n\zeta-1}} & {\tfrac{1}{\sqrt{2}} f_{n\zeta-2}} & {\cdots } & {\tfrac{1}{\sqrt{2}} f_{n1}} & {f_{n0} } & {\tfrac{1}{\sqrt{2}} f_{n1} } \\ {f_{nM-1} } & {\tfrac{1}{\sqrt{2}} f_{nM-2} } & {\tfrac{1}{\sqrt{2}} f_{nM-3} } & {\cdots } & {\tfrac{1}{\sqrt{2}} f_{n\zeta}} & {\tfrac{1}{\sqrt{2}} f_{n\zeta-1}} & {\cdots } & {\tfrac{1}{\sqrt{2}} f_{n2}} & {\tfrac{1}{\sqrt{2}} f_{n1}} & {f_{n0}} \end{smallmatrix}\right],
\]
and for different values of $M$, we take
\[\varsigma =\left\{\begin{array}{l} {(M-2)/2,\;M=4,6,8,\ldots } \\ {(M-1)/2,\;M=3,5,7,\ldots } \end{array}\right.,\; \zeta=\left\{\begin{array}{l} {{M/2 ,\hspace{.82cm} M=4,6,8,\ldots } } \\ \hspace{-1.1mm}{(M-1)/2,\hspace{1.7mm} M=3,5,7,\ldots } \end{array}\right. \,\text{and}\; f_{n\xi } =\left\{\begin{array}{l} {f_{nM-1},\;\, M=4,6,8,\ldots } \\ {0,\hspace{.82cm} M=3,5,7,\ldots} \end{array}\right.. \]

We use some approximations in above lines; it is clear that a necessary condition for the error of these approximations to become very small is that we use large $M$. Also when $f(t)$ is a polynomial function of degree $\mathsf{n}$, we must take $M \ge \mathsf{n}+2$; when we have other types such that an exponential function, a trigonometric function, ... we have to use large $M$.

\subsection{The delay operational matrix of Chebyshev wavelets}
\noindent The delay Chebyshev scaling function $\mathbf{\Psi}(t-h_{v})$ is the shifted function of $\mathbf{\Psi}(t)$ and it is given by
\begin{equation} \label{2.5.1} 
\mathbf{\Psi}(t-h_{v})=\left\{\begin{array}{ll} \mathbf{0}, & 0 \le t < h_v \\ \mathbf{D}_{v}\mathbf{\Psi}(t), &  h_v \le t \le 1. \end{array}\right. 
\end{equation} 
It follows from the discussion in \cite{iman} that by assuming $n_{v} = 2^{k-1}  h_v$, where $n_{v} \in {\mathbb{N}}$, we can write
\[
\left\{\begin{array}{ll} \bm{\varphi}_{n}^{\top}(t-h_v)=\mathbf{I}_{M}\bm{\varphi}_{n+n_{v}}^{\top}(t), & n \le 2^{k-1}-n_{v}  \\
 \bm{\varphi}_{n}^{\top}(t-h_v)=\mathbf{0}, & n > 2^{k-1}-n_{v}. \end{array}\right.
\]
Thus we conclude that the $2^{k-1} M\times 2^{k-1} M$ delay matrix $\mathbf{D}_{v}$ is
\begin{equation} \label{2.5.3} 
\mathbf{D}_{v}=\left[\begin{array}{c;{4pt/2pt}r}
   \mbox{$\mathbf{0}_{2^{k-1}M\times n_v M}$}\;\;\;
       & \hspace{-4.2mm}\begin{matrix}\mathbf{I}_{(2^{k-1}-n_v)M}\\
         \hspace{7mm}\mathbf{0}_{n_dM \times (2^{k-1}-n_v)M}
       \end{matrix}\\
        \end{array}
     \right].
\end{equation} 

\section{Simulation of the optimal tracking problem}
Consider a linear time-varying system with multiple time delays described by
\begin{equation} \label{3.1} 
\dot{\mathbf{x}}(t)=\mathbf{A}(t)\mathbf{x}(t)+\sum_{\mu=1}^{V}\mathbf{A}_{\mu}(t)\mathbf{x}(t-h_{\mu} )+\mathbf{B}(t)u(t)+\sum_{\nu=1}^{W}\mathbf{B}_{\nu}(t)\mathbf{u}(t-h_{\nu} )\, \, ,\, \, \, \, 0\le t\le t_{f}
\end{equation} 
\begin{equation} \label{3.2}
\left\{\begin{array}{ll} \mathbf{x}(t)=\mathbf{f}(t), & -h_{x} \le t\le 0  \\ \mathbf{u}(t)=\mathbf{g}(t), & -h_{u}\le t\le 0, \end{array}\right.  
\end{equation}
\begin{equation} \label{3.3} 
\mathbf{x}(0)=\mathbf{x}_{0}  
\end{equation} 
and a quadratic performance index as
\begin{equation} \label{3.j} 
J=\tfrac{1}{2} [\mathbf{x}(t_{f})-\mathbf{r}(t_{f})]^{\top} \mathbf{T}[\mathbf{x}(t_{f})-\mathbf{r}(t_{f})]+\tfrac{1}{2} \int _{0}^{t_{f}}\left\{  [\mathbf{x}(t)-\mathbf{r}(t)]^{\top}\mathbf{Q}[\mathbf{x}(t)-\mathbf{r}(t)]+\mathbf{u}^{\top}(t)\mathbf{R}\mathbf{u}(t)\right\} dt,
\end{equation}
where $\mathbf{x}(t)$ and $\mathbf{u}(t)$ are $q$- and $r$-dimensional state and control vectors, respectively, $\mathbf{A}(t)$, $\mathbf{A}_{\mu}(t)$, $\mu=1, 2, \ldots, V$, $\mathbf{B}(t)$ and $\mathbf{B}_{\nu}(t)$, $\nu=1, 2, \ldots, W$ are piecewise-continuous matrices of compatible dimensions, $h_{\mu}$ and $h_{\nu}$ denote time delays, $h_x$ and $h_u$ are the supremes of $h_{\mu}$ and $h_{\nu}$, respectively, $\mathbf{f}(t)$ is a $q$-dimensional initial state vector function, $\mathbf{g}(t)$ is an $r$-dimensional initial control vector function, $\mathbf{x}_0$ is an initial condition vector, $\mathbf{Q}$ is a positive semi-definite matrix, $\mathbf{R}$ is a positive definite matrix and $\mathbf{r}(t)$ is a $q$-dimensional desired or reference state vector. The main purpose of the matrix $\mathbf{T}$ is to ensure that the error at the terminal time is as small as possible. So, this matrix should be positive semi-definite. Our objective is to control this system in such a way that the state $\mathbf{x}(t)$ tracks the desired state $\mathbf{r}(t)$ as close as possible during the time interval $[0, t_{f} ]$. The optimal tracking problem is to find $\mathbf{u}^{*}(t)$, $\mathbf{x}^{*}(t)$ and $J^*$ for the time-delay system \eqref{3.1}--\eqref{3.3} such that the performance index in \eqref{3.j} is minimized.

First we must change the range of the independent variable $t$ such that $0\le t\le 1$; assume $\tau =t/t_{f}$. This implies that for the new independent variable $\tau \in [0, 1]$ the state equation given by eq. \eqref{3.1} has to be changed by an additional factor $t_f$ that arises from the chain rule as for $t = \tau.t_f$, the terminal time $t_f$ represents the deriviative of the interior, that is,
\begin{equation}\label{newState}
\dot{\mathbf{x}}(\tau)=t_f\Big(\mathbf{A}(\tau)\mathbf{x}(\tau)+\sum_{\mu=1}^{V}\mathbf{A}_{\mu}(\tau)\mathbf{x}(\tau-\tau_{\mu} )+\mathbf{B}(\tau)u(\tau)+\sum_{\nu=1}^{W}\mathbf{B}_{\nu}(\tau)\mathbf{u}(\tau-\tau_{\nu})\Big)\, \, ,\, \, \, \, 0\le \tau \le 1.
\end{equation}
Also we take $\tau_{\mu}=h_{\mu}/t_{f}$, $\tau_{\nu}=h_{\nu}/t_{f}$, $n_{\mu} = 2^{k-1}  \tau_{\mu}$ and $n_{\nu} = 2^{k-1}  \tau_{\nu}$.
\noindent Let us define a new state vector as
\begin{equation} \label{3.4} 
\bar{\mathbf{x}}(\tau)=\mathbf{x}(\tau)-\mathbf{r}(\tau).
\end{equation}
We parameterize this new state and the control vectors as follows
\begin{equation} \label{3_5} 
\bar{\mathbf{x}}(\tau)\cong\hat{\mathbf{\Psi} }^{\top} (\tau)\bar{\mathbf{X}},\; \mathbf{u}(\tau)\cong\check{\mathbf{\Psi}}^{\top} (\tau)\mathbf{U},
\end{equation}
where $\bar{\mathbf{X}}$ and $\mathbf{U}$ are $2^{k-1} qM\times 1$ and $2^{k-1} rM\times 1$ column vectors of unknown parameters and
\begin{equation} \label{3_6} 
\bar{\mathbf{X}}=\big{[}\bar{X}_{10}^{1} ,\,\ldots,\, \bar{X}_{10}^{q} ,\,\ldots,\, \bar{X}_{1M-1}^{1} ,\ldots ,\, \bar{X}_{1M-1}^{q} ,\,\ldots,\,\bar{X}_{2^{k-1} M-1}^{1} ,\,\ldots,\,\bar{X}_{2^{k-1} M-1}^{q}\big{]}^{\top},
\end{equation} 
\begin{equation} \label{3_7} 
\mathbf{U}=\big{[}U_{10}^{1} ,\,\ldots,\,U_{10}^{r} ,\, \ldots,\,U_{1M-1}^{1} ,\ldots ,\,U_{1M-1}^{r} ,\,\ldots,\,U_{2^{k-1} M-1}^{1} ,\,\ldots,\,U_{2^{k-1} M-1}^{r}\big{]}^{\top}.
\end{equation}
We expand the initial and desired state by using \eqref{1.3} as
\begin{equation} \label{3.8.0}
\mathbf{x}_{0} =\hat{\mathbf{\Psi}}^{\top} (\tau)\mathbf{X}_0,
\end{equation}
\begin{equation} \label{3.8.1}
\mathbf{r}(\tau)=\hat{\mathbf{\Psi}}^{\top}(\tau)\mathbf{\Gamma},
\end{equation}
where $\mathbf{X}_{0}$ and $\mathbf{\Gamma}$ are known $2^{k-1} qM\times 1$ column vectors given by 
\begin{equation} \label{3_9} 
\mathbf{X}_0=\sqrt{\tfrac{\pi}{2^{k}}}\left[\mathbf{X}_{10}^{0}, \mathbf{X}_{20}^{0}, \ldots, \mathbf{X}_{2^{k-1} 0}^{0}\right]^{\top}  
, 
\mathbf{X}_{n0}^{0} =[\mathbf{x}_{0}^{\top},\overbrace{0,0,\ldots ,0}^{q(M-1)}],
\end{equation}
\begin{equation} \label{3.10} 
\mathbf{\Gamma}=\big{[}\Gamma_{10}^{1} ,\,\ldots,\, \Gamma_{10}^{q} ,\,\ldots,\,\Gamma_{1M-1}^{1} ,\ldots ,\, \Gamma_{1M-1}^{q} ,\,\ldots,\,\Gamma_{2^{k-1} M-1}^{1} ,\,\ldots,\,\Gamma_{2^{k-1} M-1}^{q}\big{]}^{\top},
\end{equation}
in which $\Gamma_{nm}^{\alpha}$, $\alpha=1,2,\ldots,q$ can be obtained by eq.\eqref{1.6}. Using \eqref{2.5.1}, we can write
\begin{equation} \label{3.11}
\mathbf{r}(\tau-\tau_{\mu})=\hat{\mathbf{\Psi }}^{\top} (\tau)\hat{\mathbf{D}}_{\mu}^{\top} \mathbf{\Gamma}.
\end{equation}
If $\left\{\begin{array}{l}{0 \le \tau \le \tau_{\mu}}\\{0 \le \tau \le \tau_{\nu}}\end{array}\right.$ then $\left\{\begin{array}{l} {-\tau_{\mu} \le \tau-\tau_{\mu} \le 0}\\ {-\tau_{\nu} \le \tau-\tau_{\nu} \le 0} \end{array}\right.$, so according to \eqref{3.2} we have $\left\{\begin{array}{l}{\hspace{.05cm}\mathbf{x}(\tau-\tau_{\mu})=\mathbf{f}(\tau-\tau_{\mu})}\\{\mathbf{u}(\tau-\tau_{\nu})=\mathbf{g}(\tau-\tau_{\nu})}\end{array}\right.$; by \eqref{1.3}, we can write
\begin{equation}\label{funfg}
\mathbf{f}(\tau-\tau_{\mu})=\hat{\mathbf{\Psi }}^{\top}(\tau)\mathbf{F}_{\mu}, \;\mathbf{g}(\tau-\tau_{\nu})=\check{\mathbf{\Psi }}^{\top}(\tau)\mathbf{G}_{\nu},
\end{equation}
where $\mathbf{F}_{\mu}$ and $\mathbf{G}_{\nu}$ are, respectively, $2^{k-1} qM\times 1$ and $2^{k-1} rM\times 1$ column vectors of constants defined by
\begin{equation} \label{3.12} 
\mathbf{F}_{\mu}=\big{[}F_{10}^{1\mu}, \ldots, F_{10}^{q\mu}, \ldots, F_{1M-1}^{1\mu}, \ldots, F_{1M-1}^{q\mu}, \ldots, F_{n_{\mu}M-1 }^{1\mu}, \ldots,  F_{n_{\mu} M-1}^{q\mu}, \overbrace{0, 0, 0, \ldots, 0}^{(2^{k-1} -n_{\mu})qM}\big{]}^{\top},
\end{equation} 
\begin{equation} \label{3.13} 
\mathbf{G}_{\nu}=\big{[}G_{10}^{1\nu}, \ldots, G_{10}^{r\nu}, \ldots, G_{1M-1}^{1\nu}, \ldots, G_{1M-1}^{r\nu}, \ldots, G_{n_{\nu} M-1\, }^{1\nu}, \ldots, G_{n_{\nu} M-1}^{r\nu}, \overbrace{0, 0, 0, \ldots, 0}^{(2^{k-1} -n_{\nu})rM}\big{]}^{\top}.
\end{equation}
$F_{nm}^{\alpha\mu}$ and for $\beta=1,2,\ldots,r$, $G_{nm}^{\beta\nu}$ can be calculated using formula \eqref{1.6}. Thus from \eqref{2.5.1}, we find
\begin{align} \label{3_14} 
\nonumber\mathbf{x}(\tau-\tau_{\mu})&=\left\{\begin{array}{l}{\mathbf{f}(\tau-\tau_{\mu}),\hspace{2.5cm}0 \le \tau \le \tau_{\mu}}\\{\hat{\mathbf{\Psi }}^{\top}(\tau)\hat{\mathbf{D}}_{\mu}^{\top}\bar{\mathbf{X}}+\mathbf{r}(\tau-\tau_{\mu}),\quad \tau_{\mu} \le \tau \le 1}\end{array} \right.\\&=\hat{\mathbf{\Psi }}^{\top}(\tau)\mathbf{F}_{\mu}+\hat{\mathbf{\Psi }}^{\top}(\tau)\hat{\mathbf{D}}_{\mu}^{\top}\bar{\mathbf{X}}+\hat{\mathbf{\Psi }}^{\top} (\tau)\hat{\mathbf{D}}_{\mu}^{\top} \mathbf{\Gamma},
\end{align}
\begin{align} \label{3.15} 
\nonumber\mathbf{u}(\tau-\tau_{\nu})&=\left\{\begin{array}{l}{\mathbf{g}(\tau-\tau_{\nu}),\hspace{.8cm}0 \le \tau \le \tau_{\nu}}\\{\check{\mathbf{\Psi }}^{\top} (\tau)\check{\mathbf{D}}_{\nu}^{\top}\mathbf{U},\quad \tau_{\nu} \le \tau \le 1}\end{array}\right.\\&=\check{\mathbf{\Psi }}^{\top} (\tau)\mathbf{G}_{\nu}+\check{\mathbf{\Psi }}^{\top} (\tau)\check{\mathbf{D}}_{\nu}^{\top}\mathbf{U}.
\end{align}
Now we express the time-varying matrices in \eqref{newState} in terms of Chebyshev scaling function. So
\begin{align}
\nonumber\mathbf{A}(\tau)&=[\mathbf{A}_{10}, \ldots, \mathbf{A}_{1M-1}, \mathbf{A}_{20}, \ldots, \mathbf{A}_{2M-1}, \ldots, \mathbf{A}_{2^{k-1}\, 0}, \ldots, \mathbf{A}_{2^{k-1}\, M-1}]\hat{\mathbf{\Psi}}(\tau) \\ &=\mathbf{A}\hat{\mathbf{\Psi}}(\tau)
\end{align}
and
\begin{align}
\nonumber\mathbf{B}(\tau)&=[\mathbf{B}_{10}, \ldots, \mathbf{B}_{1M-1}, \mathbf{B}_{20}, \ldots, \mathbf{B}_{2M-1}, \ldots, \mathbf{B}_{2^{k-1}\, 0}, \ldots, \mathbf{B}_{2^{k-1}\, M-1} ]\check{\mathbf{\Psi}} (\tau) \\ &=\mathbf{B}\check{\mathbf{\Psi}}(\tau).
\end{align}
For $\mu=1, 2, \ldots, V$ we have
\begin{align}
\nonumber\mathbf{A}_{\mu}(\tau)&=[{\mathbf{A}_{\mu}}_{10}, \ldots, {\mathbf{A}_{\mu}}_{1M-1}, {\mathbf{A}_{\mu}}_{20}, \ldots, {\mathbf{A}_{\mu}}_{2M-1}, \ldots, {\mathbf{A}_{\mu}}_{2^{k-1}\, 0}, \ldots, {\mathbf{A}_{\mu}}_{2^{k-1}\, M-1} ]\hat{\mathbf{\Psi}}(\tau) \\ &=\mathbf{A}_{\mu}\hat{\mathbf{\Psi}}(\tau);
\end{align}
and for $\nu=1, 2, \ldots, W$
\begin{align}
\nonumber\mathbf{B}_{\nu}(\tau)&=[{\mathbf{B}_{\nu}}_{10}, \ldots, {\mathbf{B}_{\nu}}_{1M-1}, {\mathbf{B}_{\nu}}_{20}, \ldots, {\mathbf{B}_{\nu}}_{2M-1}, \ldots, {\mathbf{B}_{\nu}}_{2^{k-1}\, 0}, \ldots, {\mathbf{B}_{\nu}}_{2^{k-1}\, M-1} ]\check{\mathbf{\Psi}}(\tau) \\ &=\mathbf{B}_{\nu}\check{\mathbf{\Psi}}(\tau).
\end{align}
Then, as was mentioned in cite{iman}, we integrate equation \eqref{newState} from 0 to $\tau$, substitute these definitions and use \eqref{2.8} and\eqref{2.1}, finding
\begin{equation} \label{3.15} 
\begin{array}{l} {\hat{\mathbf{\Psi}}^{\top} (\tau)\bar{\mathbf{X}}-\hat{\mathbf{\Psi}}^{\top} (\tau)\mathbf{X}_{0}+\hat{\mathbf{\Psi}}^{\top} (\tau)\mathbf\Gamma=t_{f}\big{\{}\hat{\mathbf{\Psi}}^{\top} (\tau)\hat{\mathbf{P}}^{\top} \tilde{\mathbf{A}}\bar{\mathbf{X}}+\hat{\mathbf{\Psi}}^{\top} (\tau)\hat{\mathbf{P}}^{\top} \tilde{\mathbf{A}}\mathbf\Gamma+\hat{\mathbf{\Psi}}^{\top} (\tau)\sum_{\mu=1}^{V}\big(\hat{\mathbf{P}}^{\top} \tilde{\mathbf{A}}_{\mu}\mathbf{F}_{\mu}+} \\ {\hspace{2cm}\hat{\mathbf{P}}^{\top} \tilde{\mathbf{A}}_{\mu}\hat{\mathbf{D}}_{\mu}^{\top} \bar{\mathbf{X}}+\hat{\mathbf{P}}^{\top} \tilde{\mathbf{A}}_{\mu}\hat{\mathbf{D}}_{\mu}^{\top} \mathbf\Gamma\big)+\hat{\mathbf{\Psi}}^{\top} (\tau)\hat{\mathbf{P}}^{\top} \tilde{\mathbf{B}}\mathbf{U}+\hat{\mathbf{\Psi}}^{\top} (\tau)\sum_{\nu=1}^{W}\big(\hat{\mathbf{P}}^{\top} \tilde{\mathbf{B}}_{\nu}\mathbf{G}_{\nu}+\hat{\mathbf{P}}^{\top} \tilde{\mathbf{B}}_{\nu}\check{\mathbf{D}}_{\nu}^{\top} \mathbf{U}\big)\big{\}}.} \end{array}
\end{equation} 
Thus
\begin{equation} \label{QP1}
\begin{array}{l} {\big[t_{f}\big(\hat{\mathbf{P}}^{\top} \tilde{\mathbf{A}}+\sum_{\mu=1}^{V} \hat{\mathbf{P}}^{\top} \tilde{\mathbf{A}}_{\mu}\hat{\mathbf{D}}_{\mu}^{\top}\big)-\mathbf{I}_{qs}\big]{\hspace{.15mm}} \bar{\mathbf{X}}+t_f\big[\hat{\mathbf{P}}^{\top} \tilde{\mathbf{B}}+\sum_{\nu=1}^{W}\hat{\mathbf{P}}^{\top} \tilde{\mathbf{B}}_{\nu}\check{\mathbf{D}}_{\nu}^{\top}\big]{\hspace{.05mm}}\mathbf{U}=\mathbf\Gamma-\mathbf{X}_{0}-} \\ {\hspace{4cm}t_{f}\big{\{}\hat{\mathbf{P}}^{\top} \tilde{\mathbf{A}}\mathbf\Gamma+\sum_{\mu=1}^{V}\big(\hat{\mathbf{P}}^{\top} \tilde{\mathbf{A}}_{\mu}\hat{\mathbf{D}}_{\mu}^{\top} \mathbf\Gamma+\hat{\mathbf{P}}^{\top} \tilde{\mathbf{A}}_{\mu}\mathbf{F}_{\mu}\big)+\sum_{\nu=1}^{W}\hat{\mathbf{P}}^{\top} \tilde{\mathbf{B}}_{\nu}\mathbf{G}_{\nu}\big{\}}.} \end{array}
\end{equation}
where we let $s=2^{k-1} M$.

From the definition of this wavelet, we see the fact that the time interval $[0,1]$ is divided into $2^{k-1} $ subintervals. In order to ensure continuity in the obtained states across these subintervals, the following compatibility constraint is added at the interface points ($\tau_{\iota}$) of each subinterval:

for
\[\tau_{\iota} =\frac{\iota}{2^{k-1} }, \; \iota=1,\, 2,\, ...,\, 2^{k-1} -1\]

we must have
\[\mathbf{x}(\tau_\iota^{-})=\mathbf{x}(\tau_\iota^{+}).\]
We assume that $\mathbf{r}(t)$ is in $\mathcal{C}[0,\,t_{f}]$, so for all $\alpha$, it is necessary that
\[\bar{x}_{\alpha}(\tau_\iota^{-})=\bar{x}_{\alpha}(\tau_\iota^{+})\]
and therefore,
\[\begin{array}{l}{\begin{bmatrix} \psi_{\iota\,0} & \psi_{\iota\,1}(\tau_{\iota}) & \cdots & \psi_{\iota\,M-1}(\tau_{\iota}) \end{bmatrix} \begin{bmatrix} \bar{X}_{\iota\,0}^{\alpha} & \bar{X}_{\iota\,1}^{\alpha} & \cdots & \bar{X}_{\iota\,M-1}^{\alpha} \end{bmatrix}^{\top}-}\\{\hspace{3cm}\begin{bmatrix} \psi_{\iota+1\,0} & \psi_{\iota+1\,1}(\tau_{\iota}) & \cdots & \psi_{\iota+1\,M-1}(\tau_{\iota}) \end{bmatrix}\begin{bmatrix} \bar{X}_{\iota+1\,0}^{\alpha} & \bar{X}_{\iota+1\,1}^{\alpha} & \cdots & \bar{X}_{\iota+1\,M-1}^{\alpha} \end{bmatrix}^{\top}=0.}\end{array}\]
Hence the compatibility constraint for the defined state is expressed as
\begin{equation} \label{QP2} 
\hat{\mathbf{\Psi}}_{c} \bar{\mathbf{X}}=\mathbf{0}_{(2^{k-1}-1)q\times 1},
\end{equation} 
where
\begin{equation} \label{3.17} 
\mathbf{\Psi}_{c}=\left[\begin{array}{ccccccc} {\bm{\varphi}_{1} (\tau_{1} )} & \hspace{-2mm}{-\bm{\varphi}_{2} (\tau_{1} )} & \mathbf{0}_{1\times M} & \mathbf{0}_{1\times M} & {\cdots } & \mathbf{0}_{1\times M} & \mathbf{0}_{1\times M} \\ \mathbf{0}_{1\times M} & {\bm{\varphi}_{2} (\tau_{2} )} & \hspace{-2mm}{-\bm{\varphi}_{3} (\tau_{2} )} & \mathbf{0}_{1\times M} & {\cdots } & \mathbf{0}_{1\times M} & \mathbf{0}_{1\times M} \\ \mathbf{0}_{1\times M} & \mathbf{0}_{1\times M} & {\bm{\varphi}_{3} (\tau_{3} )} & \hspace{-2mm}-{\bm{\varphi}_{3} (\tau_{3} )} & {\cdots } & \mathbf{0}_{1\times M} & \mathbf{0}_{1\times M} \\ {\vdots } & {\vdots } & {\vdots } & {\vdots } & {\ddots } & \vdots  & \vdots  \\ \mathbf{0}_{1\times M} & \mathbf{0}_{1\times M} & \mathbf{0}_{1\times M} & \mathbf{0}_{1\times M} & {\cdots } & \bm{\varphi}_{2^{k-1}} (\tau_{2^{k-1} -1} ) & \hspace{-2mm}{-\bm{\varphi}_{2^{k-1}} (\tau_{2^{k-1} -1} )} \end{array}\right].
\end{equation}

Setting \eqref{3.4} in the performance index \eqref{3.j}, we find
\begin{equation}
\begin{array}{l}{J=\tfrac{1}{2}\bar{\mathbf{X}}^{\top}\hat{\mathbf{\Psi}}(1)\mathbf{T}\hat{\mathbf{\Psi}}^{\top}(1)\bar{\mathbf{X}}+\tfrac{1}{2}t_{f}\int_{0}^{1}[\bar{\mathbf{X}}^{\top}\hat{\mathbf{\Psi}}(\tau)\mathbf{Q}\hat{\mathbf{\Psi }}^{\top}(\tau)\bar{\mathbf{X}}+ \mathbf{U}^{\top}\check{\mathbf{\Psi}}(\tau)\mathbf{R}\check{\mathbf{\Psi }}^{\top}(\tau)\mathbf{U}\,d\tau]} \\ {\hspace{.24cm}= \tfrac{1}{2}\{\bar{\mathbf{X}}^{\top}(\mathbf{\Psi}(1)\mathbf{\Psi}^{\top}(1)\otimes \mathbf{T})\bar{\mathbf{X}}+t_{f}\int_{0}^{1}[\bar{\mathbf{X}}^{\top}(\mathbf{\Psi}(\tau) \mathbf{\Psi }^{\top}(\tau) \otimes \mathbf{Q})\bar{\mathbf{X}}+ \mathbf{U}^{\top}(\mathbf{\Psi}(\tau) \mathbf{\Psi }^{\top}(\tau) \otimes \mathbf{R})\mathbf{U}\,d\tau]\}} \\ {\hspace{.24cm}=\tfrac{1}{2}\{\bar{\mathbf{X}}^{\top}(\mathbf{\Psi}(1)\mathbf{\Psi}^{\top}(1)\otimes \mathbf{T}+t_{f}\mathbf{C} \otimes \mathbf{Q})\bar{\mathbf{X}}+ \mathbf{U}^{\top}(t_{f}\mathbf{C} \otimes \mathbf{R})\mathbf{U}\}.}\end{array}
\end{equation}
As a result
\begin{equation}\label{QP3}
J=\tfrac{1}{2}[\begin{array}{cc}{\bar{\mathbf{X}}} & {\mathbf{U}} \end{array}] \left[\begin{array}{cc} {t_{f}\mathbf{C}\otimes \mathbf{Q}+(\mathbf{\Psi}(1)\mathbf{\Psi}^{\top}(1)\otimes \mathbf{T})} & {\mathbf{0}} \\ \mathbf{0} & {t_{f}\mathbf{C} \otimes \mathbf{R}} \end{array}\right] [\begin{array}{cc}{\bar{\mathbf{X}}} & {\mathbf{U}} \end{array}]^{\top}.
\end{equation}
Taking \eqref{QP3}, \eqref{QP2} and \eqref{QP1} together, the optimal tracking control problem is transformed into a quadratic programming (QP) problem:

\[
\min_{\bm{\chi}}\;\tfrac{1}{2}\bm{\chi}^{\top} \bm{\aleph} \bm{\chi}
\] 
\[\hspace{-9mm}{\text{subject\,\,to}}\;\, \bm{\Lambda} \bm{\chi}=\mathbf{b},\] 
where we set
\begin{equation} \label{3.18} 
\bm{\chi}^{\top}=[\begin{array}{cc} {\bar{\mathbf{X}}} & {\mathbf{U}} \end{array}],
\end{equation} 
\begin{equation} \label{3.19} 
\bm{\aleph}=\left[\begin{array}{cc} {t_{f}\mathbf{C}\otimes \mathbf{Q}+(\mathbf{\Psi}(1)\mathbf{\Psi}^{\top}(1)\otimes \mathbf{T})} & {\mathbf{0}_{qs\times rs} } \\ \mathbf{0}_{rs\times qs} & {t_{f}\mathbf{C} \otimes \mathbf{R}} \end{array}\right],
\end{equation} 
\begin{equation} \label{3.20} 
\bm{\Lambda}=\left[\begin{array}{cc} {t_{f}(\hat{\mathbf{P}}^{\top} \tilde{\mathbf{A}}+\sum_{\mu=1}^{V}\hat{\mathbf{\mathbf{P}}}^{\top} \tilde{\mathbf{A}}_{\mu}\hat{\mathbf{D}}_{\mu}^{\top})-\mathbf{I}_{qs} } & {t_{f}(\hat{\mathbf{P}}^{\top} \tilde{\mathbf{B}}+\sum_{\nu=1}^{W}\hat{\mathbf{P}}^{\top} \tilde{\mathbf{B}}_{\nu}\check{\mathbf{D}}_{\nu}^{\top}) } \\ {\hat{\mathbf{\Psi}}_{c} } & {\mathbf{0}_{(2^{k-1} -1)q\times rs} } \end{array}\right],
\end{equation} 
\begin{equation} \label{3.21} 
\mathbf{b}=\left[\begin{array}{c} {\mathbf\Gamma-\mathbf{X}_{0}-t_{f}\left\{\hat{\mathbf{P}}^{\top} \tilde{\mathbf{A}}\mathbf\Gamma+\sum_{\mu=1}^{V}\big(\hat{\mathbf{P}}^{\top} \tilde{\mathbf{A}}_{\mu}\hat{\mathbf{D}}_{\mu}^{\top} \mathbf\Gamma+\hat{\mathbf{P}}^{\top} \tilde{\mathbf{A}}_{\mu}\mathbf{F}_{\mu}\big)+\sum_{\nu=1}^{W}\hat{\mathbf{P}}^{\top} \tilde{\mathbf{B}}_{\nu}\mathbf{G}_{\nu}\right\}}\\ {\mathbf{0}_{(2^{k-1} -1)q\times 1} } \end{array}\right].
\end{equation}

Our goal is to find $\bm{\chi}$ from solving the latter optimization problem which is static in nature. Standard numerical methods are available to solve this QP problem. Hence we do not need a special program. We can use the quadprog function provided by the optimization toolbox in MATLAB. This toolbox presents widely used algorithms to solve constrained and unconstrained optimization problems. We need to set what MATLAB solver to use with the algorithm field in the optimization options. In this work we use the interior-point-convex algorithm of the quadprog function in MATLAB R2012b; the proposed algorithm provides an accurate solution, and is fast and stable. For more detailed information about handling various cases of plant matrices, constraints and ..., see \cite{iman}.

\section{Numerical examples}
\subsection{Example 1}
We are interested in finding the optimal state and control which cause the time-delay system
\begin{equation} \label{ex0.1} 
\dot{x}(t)=t^2x(t)-3tx(t-{\tfrac{1}{2}} )+2u(t)+u(t-{\tfrac{1}{2}} )\, ,\;\; 0\le t\le 1 
\end{equation} 
\begin{equation} \label{ex0.2} 
x(t)=t^2+1\,,\;\; -{\textstyle\frac{1}{2}} \le \, t\le 0 
\end{equation} 
\begin{equation} \label{ex0.3} 
u(t)=t+1\,,\hspace{2.8mm} -{\textstyle\frac{1}{2}} \le \, t\le 0 
\end{equation} 
to follow the desired state
\begin{equation}\label{ex0.4}
r(t)=\left\{\begin{array}{l}{9t^2-6t+1, \hspace{0.5cm}\,0 \le t < 0.5}\\{0.25, \hspace{1.5cm}\,\,0.5 \le t <1,} \end{array}\right.
\end{equation}
while minimize the performance index
\begin{equation} \label{ex0.5} 
J=\tfrac{1}{4}[x(1)-r(1)]^2+ \int _{0}^{1}\left\{[x(t)-r(t)]^{2}+R_{d}u^{2} (t)\right\}dt,
\end{equation}
where

a. $R_{d} = 0.005,$ \hspace{5.20cm} b. $R_{d}=\dfrac{0.005}{5t+1}.$
\\

In case a, we have $R = 0.01$. We first set $k=2$ and $M=5$; next, let us identify the required matrices in the present tracking system. Obviously $ T = \tfrac{1}{2}, Q = 2, B =2, B_{1}=1$. From the findings given in sections 2 and 3 we have

\noindent $\mathbf{P}=\tfrac{1}{4}\left[\begin{smallmatrix} \mathbf{L} & \mathbf{E} \\ \mathbf{0} & \mathbf{L} \end{smallmatrix}\right], \mathbf{C}=\tfrac{2}{\pi}.{\text{blkdiag}}(\mathbf{C}_1,\mathbf{C}_2)$, where
\[
\mathbf{L}= {\footnotesize\begin{bmatrix} 1 & \tfrac{1}{\sqrt{2}} & 0 & 0 & 0 \\ -\tfrac{\sqrt{2}}{4} & 0 & \frac{1}{4} & 0 & 0 \\ -\tfrac{\sqrt{2}}{3} & -\tfrac{1}{2} & 0 & \tfrac{1}{6} & 0 \\ \tfrac{\sqrt{2}}{8} & 0 & -\frac{1}{4} & 0 & \frac{1}{8} \\ -\tfrac{\sqrt{2}}{15} & 0 & 0 & -\frac{1}{6} & 0 \end{bmatrix}},\mathbf{E}= {\footnotesize\begin{bmatrix} 2 & 0 & 0 & 0 & 0 \\ 0 & 0 & 0 & 0 & 0 \\ -\tfrac{2\sqrt{2}}{3} & 0 & 0 & 0 & 0 \\ 0 & 0 & 0 & 0 & 0 \\ -\tfrac{2\sqrt{2}}{15} & 0 & 0 & 0 & 0 \end{bmatrix}}, \mathbf{C}_{1}=\mathbf{C}_{2}= {\footnotesize\left [\begin{smallmatrix} 1 & {0} & \frac{-\sqrt{2}}{3} & {0} & \frac{-\sqrt{2}}{15}  \\ {0} & \frac{2}{3} & {0} & \frac{-2}{5} & {0} \\ \frac{-\sqrt{2} }{3} & {0} & \frac{14}{15} & {0} & \frac{-38}{105} \\ {0} & \frac{-2}{5} & {0} & \frac{34}{35} & {0} \\ \frac{-\sqrt{2}}{15} & {0} & \frac{-38}{105} & {0} & \frac{62}{63}  \end{smallmatrix} \right],}
\]
\[
\mathbf{X}_0=\tfrac{\sqrt{\pi}}{2}\left[1,0,0,0,0,1,0,0,0,0\right]^{\top},\,\mathbf{\Gamma}=\left[\tfrac{11\sqrt{\pi}}{64},-\tfrac{3\sqrt{2\pi}}{32},\tfrac{9\sqrt{2\pi}}{128},0,0,\tfrac{\sqrt{\pi}}{8},0,0,0,0\right]^{\top},
\]
$
\mu=\nu=1$, $f(t-\tfrac{1}{2})=t^2-t+\tfrac{5}{4}$, $g(t-\tfrac{1}{2})=t+\tfrac{1}{2}$, $n_{\mu}=n_{\nu}=1$ and $t_{\iota} =\tfrac{1}{2}$, thus
\[
\mathbf{F}_{1}=\left[\tfrac{35\sqrt{\pi}}{64},-\tfrac{\sqrt{2\pi}}{32},\tfrac{\sqrt{2\pi}}{128},0,0,0,0,0,0,0\right]^{\top}
,
\mathbf{G}_{1}=\left[\tfrac{3\sqrt{\pi}}{8},\tfrac{\sqrt{2\pi}}{16},0,0,0,0,0,0,0,0\right]^{\top},
\]
\[
\mathbf{D}_{1}=\left[\begin{array}{c;{4pt/2pt}r}
   \mbox{$\mathbf{0}_{10\times 5}$} 
       & \begin{smallmatrix}\mathbf{I}_{5}\\
        \mathbf{0}_{5 \times 5}
       \end{smallmatrix}\\
        \end{array}
     \right]
, 
\mathbf{\Psi}_{c}=\left[ \tfrac{2}{\sqrt{\pi }}, \tfrac{2\sqrt{2}}{\sqrt{\pi}}, \tfrac{2\sqrt{2}}{\sqrt{\pi}}, \tfrac{2\sqrt{2}}{\sqrt{\pi}}, \tfrac{2\sqrt{2}}{\sqrt{\pi}}, \tfrac{-2}{\sqrt{\pi }}, \tfrac{2\sqrt{2}}{\sqrt{\pi}}, \tfrac{-2\sqrt{2}}{\sqrt{\pi}}, \tfrac{2\sqrt{2}}{\sqrt{\pi}}, \tfrac{-2\sqrt{2}}{\sqrt{\pi}} \right],
\]
$\mathbf{\Psi}(1)\mathbf{\Psi}^{\top}(1)=\left[\begin{smallmatrix} {\mathbf{0}_{5\times 5} } & {\mathbf{0}_{5\times 5} } \\ {\mathbf{0}_{5\times 5} } & {\mathbf{Y}} \end{smallmatrix}\right]$; $A(t)=t^2$ and $A_1(t)=-3t$, so $\tilde{\mathbf{A}}={\footnotesize\tfrac{2}{\sqrt{\pi}}}.{\text{blkdiag}}(\tilde{\mathbf{A}}_1,\tilde{\mathbf{A}}_2), \tilde{\mathbf{A}}_{1}={\footnotesize\tfrac{2}{\sqrt{\pi}}}.{\text{blkdiag}}(\tilde{\mathbf{A}}_{1_1},\tilde{\mathbf{A}}_{1_2})$, where
\[
\mathbf{Y}= {\footnotesize\left [\begin{smallmatrix} \tfrac{4}{\pi} & \tfrac{4\sqrt{2}}{\pi} & \tfrac{4\sqrt{2}}{\pi} & \tfrac{4\sqrt{2}}{\pi} & \tfrac{4\sqrt{2}}{\pi} \\ \tfrac{4\sqrt{2}}{\pi} & \frac{8}{\pi} & \frac{8}{\pi} & \frac{8}{\pi} & \frac{8}{\pi} \\ \tfrac{4\sqrt{2}}{\pi} & \frac{8}{\pi} & \frac{8}{\pi} & \frac{8}{\pi} & \frac{8}{\pi} \\ \tfrac{4\sqrt{2}}{\pi} & \frac{8}{\pi} & \frac{8}{\pi} & \frac{8}{\pi} & \frac{8}{\pi} \\ \tfrac{4\sqrt{2}}{\pi} & \frac{8}{\pi} & \frac{8}{\pi} & \frac{8}{\pi} & \frac{8}{\pi} \end{smallmatrix} \right], \tilde{\mathbf{A}}_{n}={\footnotesize\left [\begin{smallmatrix} A_{n0} & A_{n1} & A_{n2} & A_{n3} & A_{n4} \\ A_{n1} & A_{n0}+\tfrac{A_{n2}}{\sqrt{2}} & \tfrac{A_{n1}+A_{n3}}{\sqrt{2}} & \tfrac{A_{n2}+A_{n4}}{\sqrt{2}} & \tfrac{A_{n3}}{\sqrt{2}} \\ A_{n2} & \tfrac{A_{n1}+A_{n3}}{\sqrt{2}} & A_{n0}+\tfrac{A_{n4}}{\sqrt{2}} & \tfrac{A_{n1}}{\sqrt{2}} & \tfrac{A_{n2}}{\sqrt{2}} \\ A_{n3} & \tfrac{A_{n2}+A_{n4}}{\sqrt{2}} & \tfrac{A_{n1}}{\sqrt{2}} & A_{n0} & \tfrac{A_{n1}}{\sqrt{2}} \\ A_{n4} & \tfrac{A_{n3}}{\sqrt{2}} & \tfrac{A_{n2}}{\sqrt{2}} & \tfrac{A_{n1}}{\sqrt{2}} & A_{n0}  \end{smallmatrix}\right ]},}
\]
\[
A_{10}=\tfrac{3\sqrt{\pi}}{64}, A_{11}=\tfrac{\sqrt{2\pi}}{32}, A_{12}=\tfrac{\sqrt{2\pi}}{128}, A_{13}=A_{14}=0, A_{20}=\tfrac{19\sqrt{\pi}}{64}, A_{21}=\tfrac{3\sqrt{2\pi}}{32}, A_{22}=\tfrac{\sqrt{2\pi}}{128}, A_{23}=A_{24}=0,
\]
\[
A_{1_{10}}=\tfrac{-3\sqrt{\pi}}{8}, A_{1_{11}}=\tfrac{-3\sqrt{2\pi}}{16}, A_{1_{12}}=A_{1_{13}}=A_{1_{14}}=0, A_{1_{20}}=\tfrac{-9\sqrt{\pi}}{8}, A_{1_{21}}=\tfrac{-3\sqrt{2\pi}}{16}, A_{1_{22}}=A_{1_{23}}=A_{1_{24}}=0.
\]
Substituting the above-mentioned results into \eqref{3.19}--\eqref{3.21} and calling the quadprog algorithm in MATLAB yields $J^*=0.008801$, and
\[x^{*}(t)=\left\{\begin{array}{ll} - 6.73743t^4 + 7.14247t^3 + 6.23845t^2 - 5.53175t + 0.99999
, & t \in [0, 0.5] \\ 7.18440t^4 - 22.77482t^3 + 26.94463t^2 - 14.05418t + 2.95420
, &  t \in [0.5, 1], \end{array}\right.\]
\[u^{*}(t)=\left\{\begin{array}{ll} - 4.66593t^4 - 8.72582t^3 + 8.58942t^2 + 7.62985t - 3.01649
, & t \in [0, 0.5] \\ - 5.49792t^4 + 50.80960t^3 - 89.94716t^2 + 51.69140t - 7.44087
, &  t \in [0.5, 1]. \end{array}\right.\]

\begin{figure}[!ht]
\centering
\subfigure[$R=0.01$]{
    \includegraphics[scale=.5]{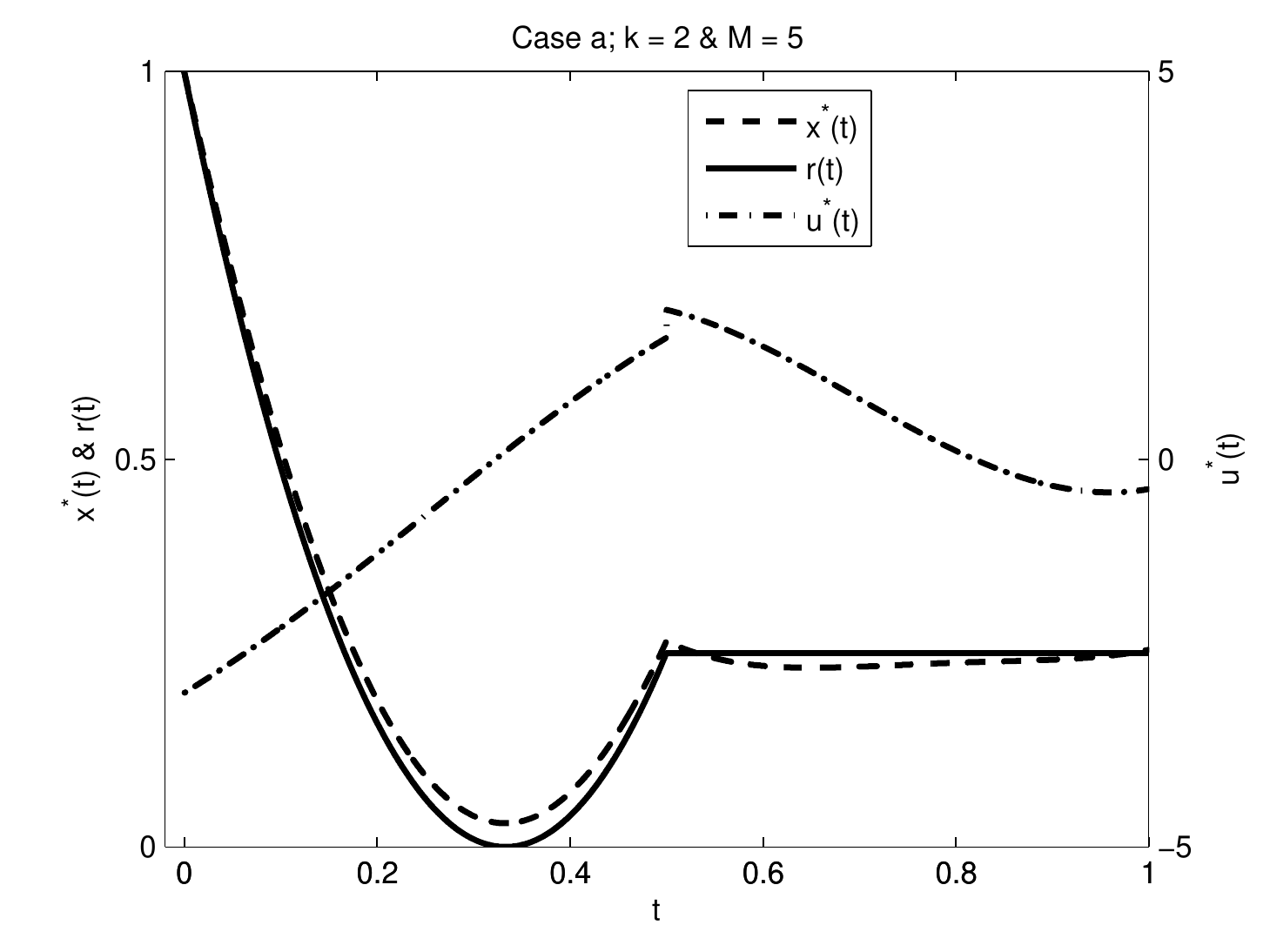}
    \label{Fig:0.1}
}
\subfigure[$R=\tfrac{0.01}{5t+1}$]{
    \includegraphics[scale=.5]{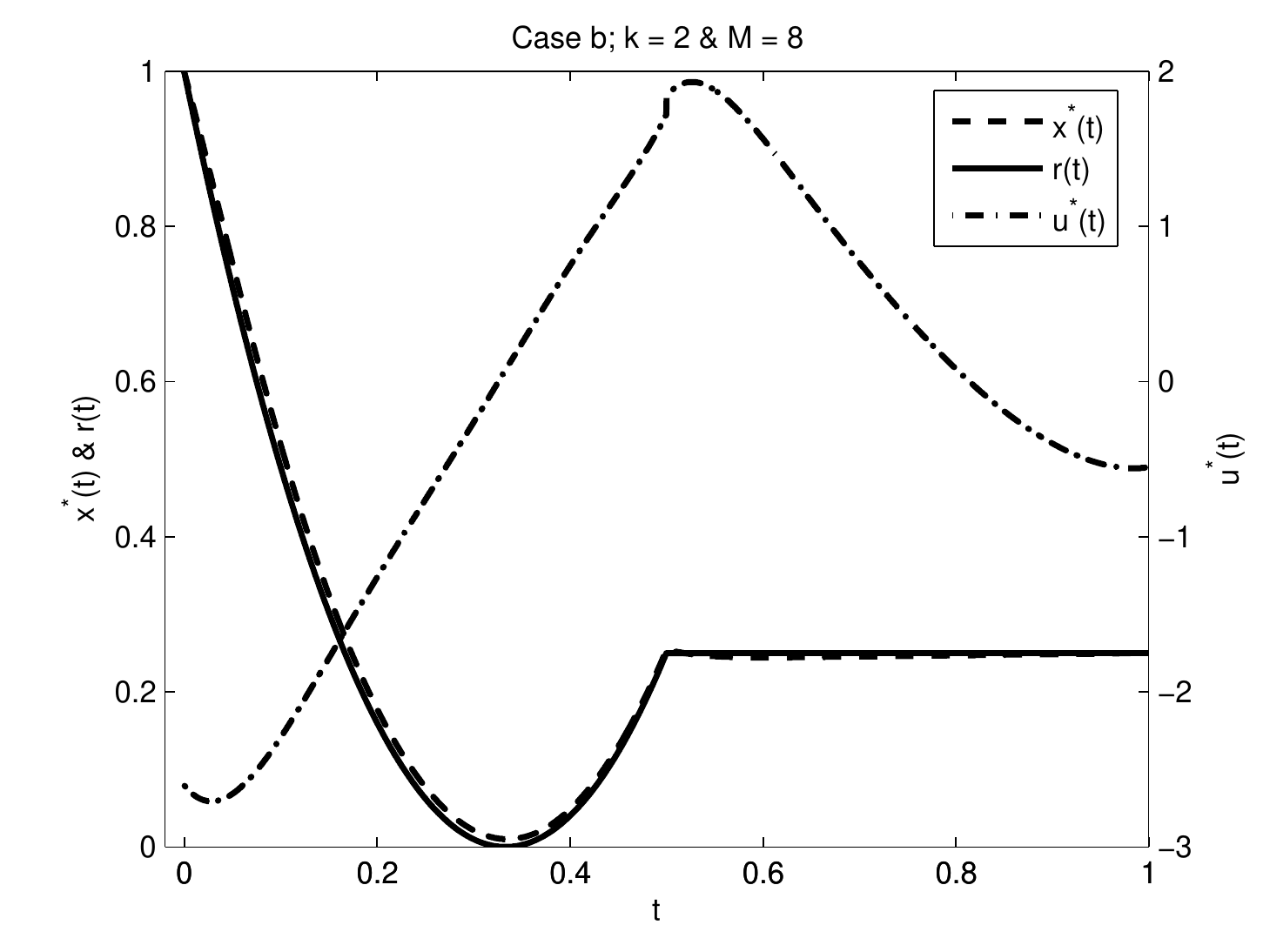}
    \label{Fig:0.2}
}
\caption[Optional caption for list of figures]{Reference state, optimal state and control for Example 1.}
\label{Fig:0}
\end{figure}

In case b, we have $R(t)=\frac{0.01}{5t+1}$. At first glance solving this problem by the method may not look easy, but it is; in order to do this, by letting $R(t)=\mathbf{R}'\mathbf{\Psi}(t)$, then using \eqref{2.8} and \eqref{2.5}, finding a approximation of $\bm{\aleph}$ is not difficult. Thus, taking $k=2$ and $M=8$ gives us $J^*=0.004968$. The optimal state and control, and the reference state are shown in Fig.\ref{Fig:0}. If we look at the obtained curves, it is obvious that we get a better tracking of the desired state with lower cost in case b.

\subsection{Example 2}
\noindent Consider the system with small delay terms contained in the state and control vectors \cite{Leondes.Shieh}. The problem involves the minimization of
\begin{equation} \label{5_8_} 
J=\int _{0}^{15}\left\{ [x_{1} (t)-0.2t]^{2} +0.025u^{2} (t)\right\} dt  
\end{equation} 
subjected to the system of delayed differential equations and initial conditions such as
\begin{equation} \label{5_9_} 
\dot{x}_{1} (t)=0.05x_{1} (t-1)+x_{2} (t)+0.01u(t-0.5),
\end{equation} 
\begin{equation} \label{5_10_} 
\dot{x}_{2} (t)=2x_{1} (t)+0.01x_{2} (t-1)-x_{2} (t)+u(t)-0.05u(t-0.5),
\end{equation} 
\begin{equation} \label{5_11_} 
\left. \begin{array}{l} {x_{1} (t)=-4 ,\;-1\le t\le 0} \\ {x_{2} (t)=0,\hspace{2.95mm}-1\le t\le 0} \\ {u(t)=0,\hspace{1.4mm} -0.5\le t\le 0.} \end{array}\right\} 
\end{equation} 
for the case which $x(t_{f}=15)$ is free and admissible optimal control and states are unbounded.

\begin{figure}[!ht]
\centering
\subfigure[Reference state and optimal states]{
    \includegraphics[scale=.5]{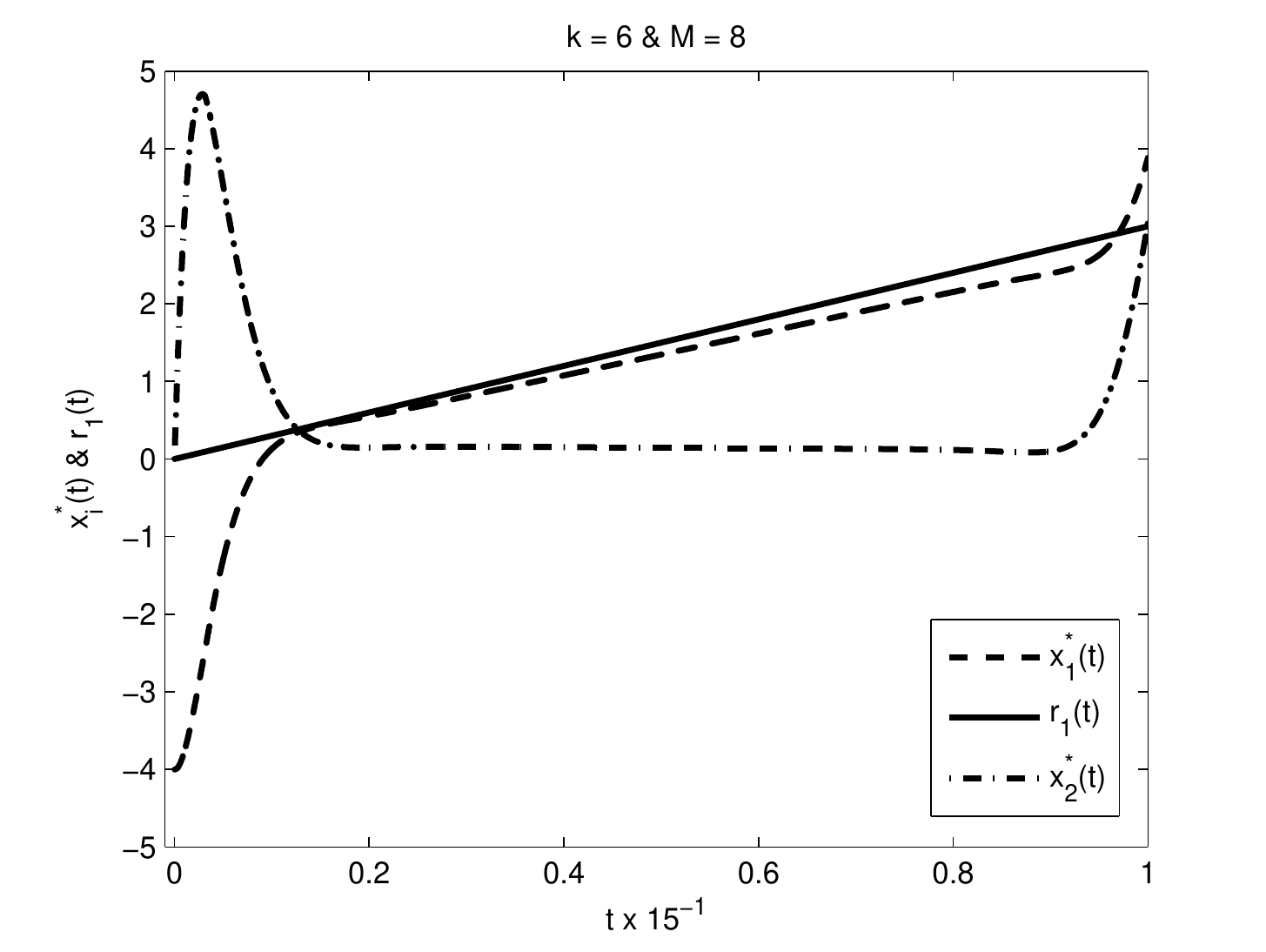}
    \label{Fig:1.1}
}
\subfigure[Optimal control]{
    \includegraphics[scale=.5]{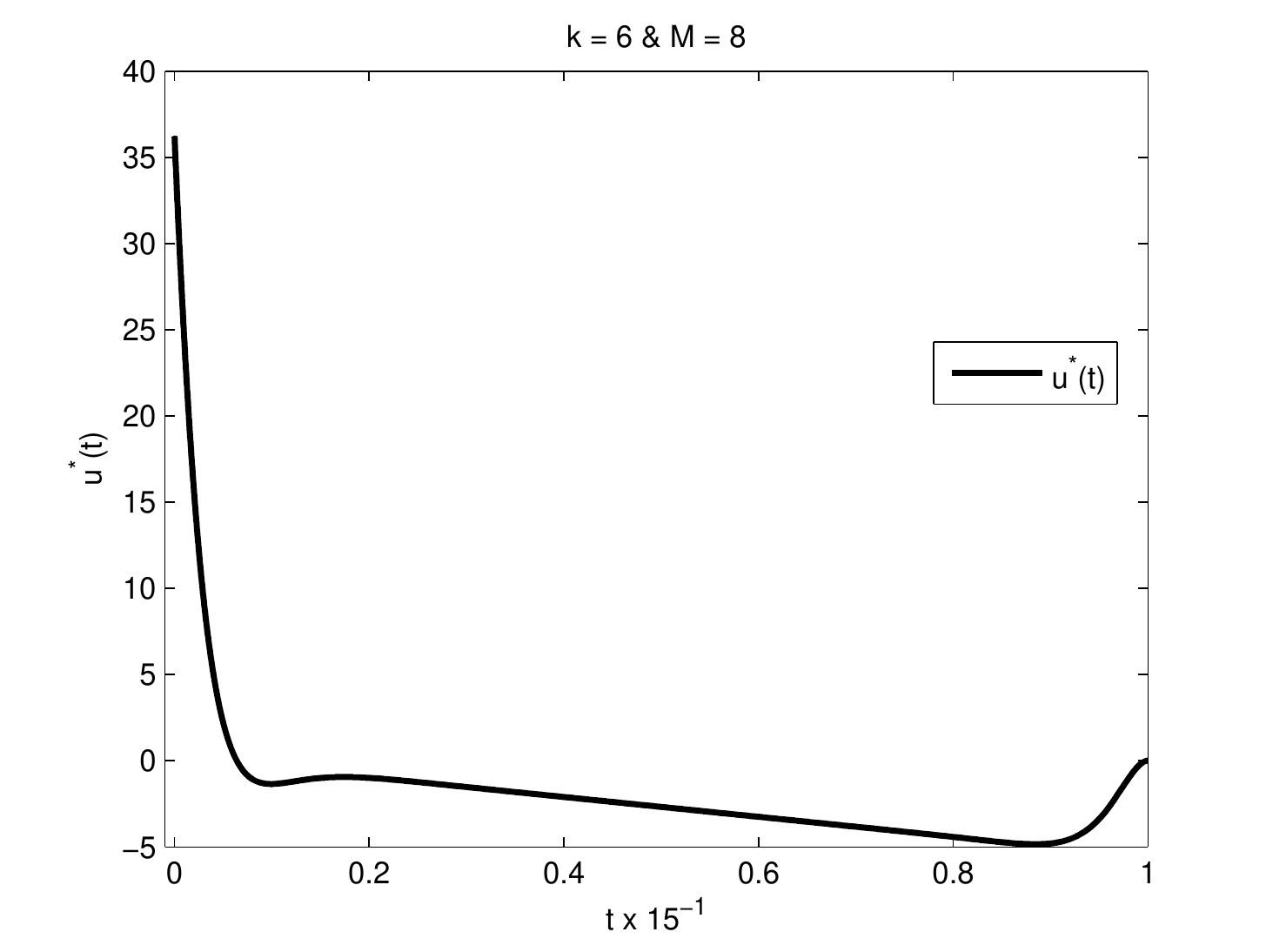}
    \label{Fig:1.2}
}
\caption[Optional caption for list of figures]{Optimal states and control for Example 2.}
\label{Fig:1}
\end{figure}

The performance index indicates that the state $x_1(t)$ is to be kept close to the reference state $r_1(t)=0.2t$ and since there is
no condition on the state $x_2(t)$, we let $r_2(t)=0$. After rescaling the time interval by setting $\tau =t/15$, we set $\bar{x}_{1}(\tau )=x_{1}(\tau )-3\tau $. Thus
\[\mathbf{r}(\tau) =\left[\begin{array}{c} {3\tau} \\ {0 } \end{array}\right]\;\text{and}\;\bar{x}_{1}(\tau -{\textstyle\frac{1}{15}} )=\left\{\begin{array}{l} {-4,\, \, \, \, \, \, \, 0\le \tau \le {\textstyle\frac{1}{15}} } \\ {\,\,0,\, \, \, \, \, \, \, \, \, {\textstyle\frac{1}{15}} \le \tau \le 1.} \end{array}\right. \] 
By choosing $k=6$ and $M=8$, formulas \eqref{3.10} and \eqref{3.12} give
\[
\mathbf\Gamma =3\Big[\tfrac{\sqrt{\pi }}{512},0,{\textstyle\frac{\sqrt{2\pi } }{1024}} ,0,\overbrace{0,\ldots,0}^{12},{\textstyle\frac{3\sqrt{\pi } }{512}} ,0,{\textstyle\frac{\sqrt{2\pi } }{1024}},0,\overbrace{0,\ldots,0}^{12} ,\ldots ,{\textstyle\frac{63\sqrt{\pi } }{512}} ,0,{\textstyle\frac{\sqrt{2\pi } }{1024}},0, \overbrace{0,\ldots,0}^{12}\Big]^{\top},
\]
\[
\mathbf{F}_{1}=\Big[\tfrac{-\sqrt{\pi}}{2},0,\overbrace{0,\ldots,0}^{14},\tfrac{-\sqrt{\pi } }{2},0,\overbrace{0,\ldots,0}^{14} , \overbrace{0,\ldots,0}^{480}\Big]^{\top}.
\]
Also we have
\[\mathbf{A}=\begin {bmatrix} 0& \hspace{.2cm}1 \\ 2&-1 \end {bmatrix}, \mathbf{A}_{1}=\begin{bmatrix} 0.05&0\\0&0.01\end{bmatrix}, \mathbf{B}=\begin{bmatrix}0\\1 \end{bmatrix}, \mathbf{B}_{1}=\begin{bmatrix} \hspace{.2cm}0.01\\-0.05\end{bmatrix}, \mathbf{T}=\mathbf{0}_{2 \times 2}, \mathbf{Q}=\begin{bmatrix}2&0\\0&0\end{bmatrix} \, \text{and} \, R=0.05.\]
Now we can solve the problem. The simulation curves are presented in Fig.\ref{Fig:1.1} and Fig.\ref{Fig:1.2}.

\begin{figure}[!ht]
\centering
\subfigure[Reference state and optimal states]{
    \includegraphics[scale=.5]{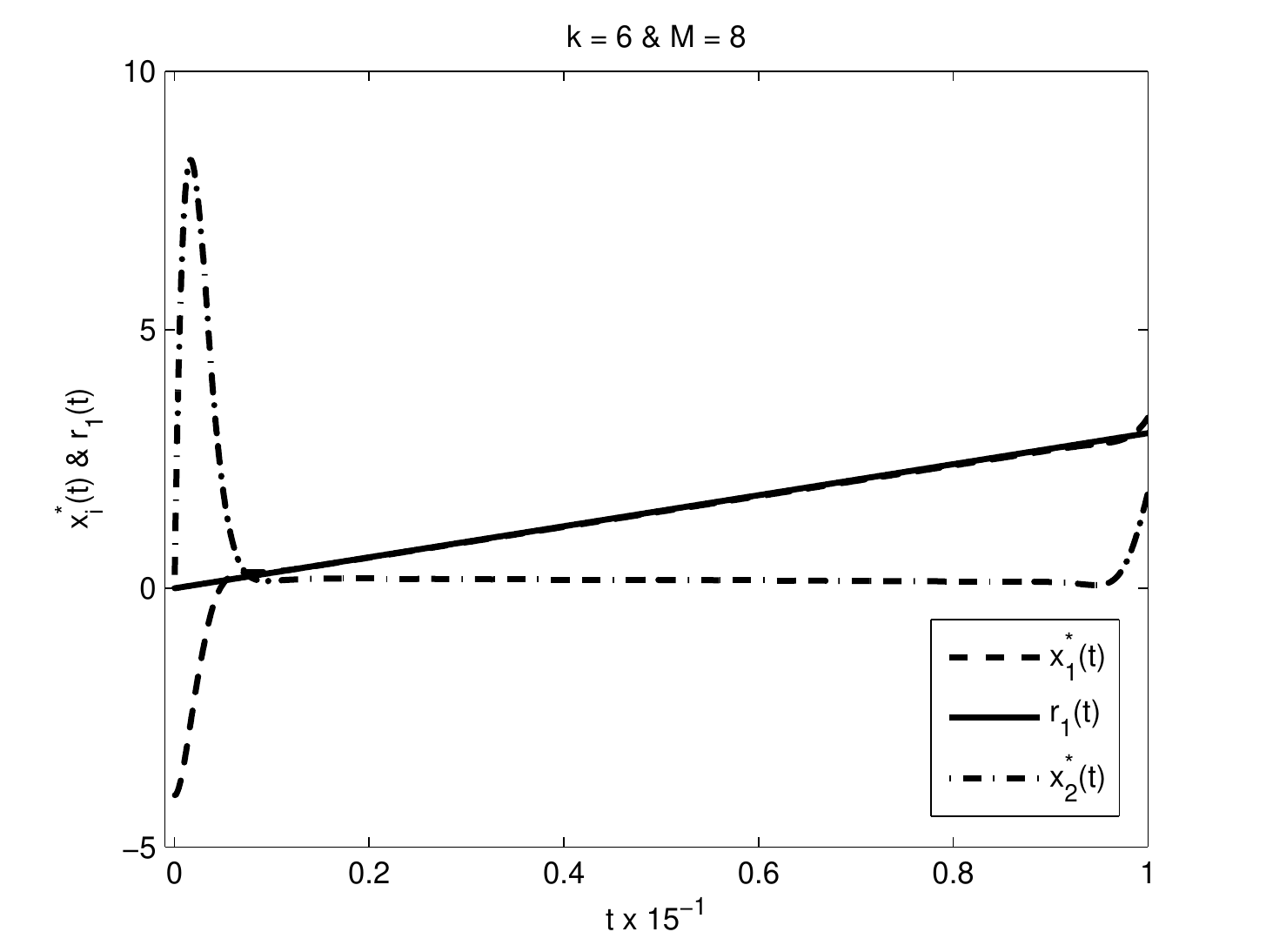}
    \label{Fig:1.1b}
}
\subfigure[Optimal control]{
    \includegraphics[scale=.5]{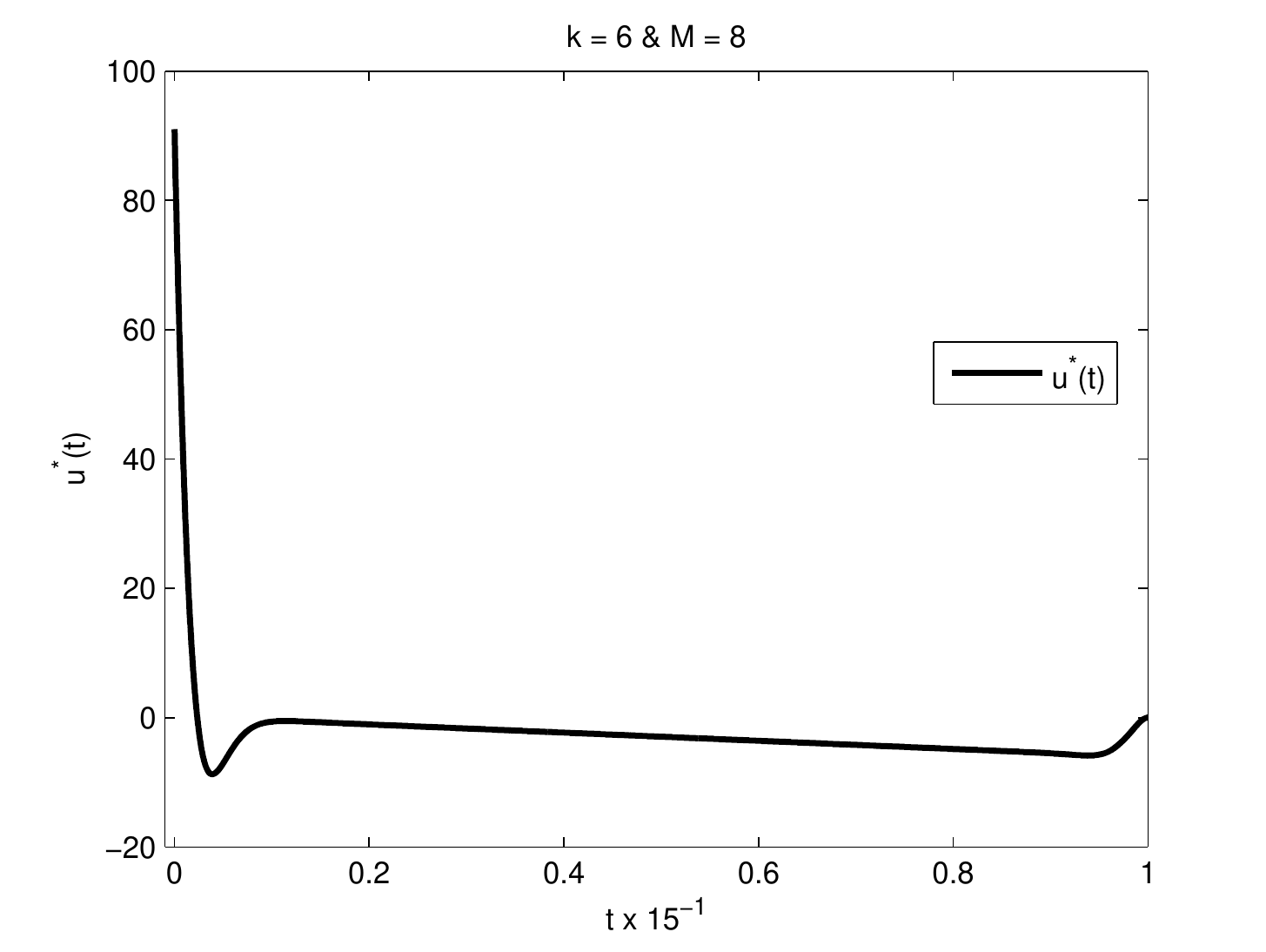}
    \label{Fig:1.2b}
}
\caption[Optional caption for list of figures]{Optimal states and control for Example 2 with new $\mathbf{Q}$.}
\label{Fig:1b}
\end{figure}

The optimal value of the cost functional \eqref{5_8_} is found to be $J^{*}=16.636902$ and a comparison of $J^*$ is given in Table \ref{tab:ex1.1}. It is clear that the obtained result is in good agreement. For performing a better tracking by the system we have to increase the value of weighting matrix $\mathbf{Q}$ such that:
\[
\mathbf{Q}_{\text{new}}=10\mathbf{Q}=\left[\begin{smallmatrix}20&0\\0&0\end{smallmatrix} \right].
\]
$\mathbf{x}^*(t)$ and $u^*(t)$ for this weighting matrix are plotted in Figs.\ref{Fig:1.1b}--\ref{Fig:1.2b}; also we get $J^*=60.853249$. These results mean that when we increase the values of the weighting matrix $\mathbf{Q}$, the state of the new system is able to track the reference state better with lower error, but we have to pay higher cost for larger control effort. Suppose, to achieve a better tracking, that instead of increasing the value of the error weighted matrix $\mathbf{Q}$, we decrease the value of the control weighted matrix $R$, such as:
\[ 
R_{\text{new}}=0.1R=0.005.
\]
This gives $J^*=6.0853249$, so we get a lower cost while the graphs of $\mathbf{x}^*(t)$ and $u^*(t)$ are exactly the same as those obtained with new $\mathbf{Q}$ (Fig.\ref{Fig:1b}). We want on one hand, to keep the new state small and on the other hand, we must not pay higher cost to large controls; this leads us immediately to the conclusion that we have to try various values of the weighting matrix $\mathbf{R}$.
\begin{table}[h]
\centering
\caption{Comparison of $J^{*}$ for Example 2}\label{tab:ex1.1}
\begin{tabular}{ll}
\toprule \noalign{\vskip 1mm} 
Source & optimal performance criteria ($J^{*}$) \\\noalign{\vskip 1mm} 
\midrule
C.T. Leondes, E. Shieh {\small\cite{Leondes.Shieh}} & 17.2993 \\
this research & 16.636902 \\ 
\bottomrule
\end{tabular}
\end{table}

\subsection{Example 3}
This example studied in \cite{Tang.Li.Zhao}. Consider the following optimal tracking problem for a time-delay system:
\begin{equation} \label{Ex3.0} 
\begin{array}{l} \dot{\mathbf{x}}(t) = \begin{bmatrix} {0} & {1} \\ {-1 } & {1} \end{bmatrix} \mathbf{x}(t) + \begin{bmatrix} {-1} & {0} \\ {0.6} & {-1.5} \end{bmatrix} \mathbf{x}(t-h_{x} ) + \begin{bmatrix} {0} \\ {2} \end{bmatrix} u(t), \end{array}
\end{equation} 
\begin{equation} \label{Ex3.1} 
\mathbf{x}(t)=\left[\begin{array}{cc} {3} & {0} \end{array}\right]^{\top} \, ,\, \, -h_{x} \le t\le 0 
\end{equation}
\begin{equation} \label{Ex3.2} 
y(t)=\left[\begin{array}{cc} {2} & {0} \end{array}\right]\mathbf{x}(t),\;\bar{y}(t)=\left[\begin{array}{cc} {1} & {0} \end{array}\right]\mathbf{z}(t),
\end{equation}
\begin{equation} \label{Ex3.3} 
\dot{\mathbf{z}}(t) = \begin{bmatrix} \hspace{-.8mm}{0} & \hspace{2mm}{0.1}\\{-0.2} & {-0.3} \end{bmatrix} \mathbf{z}(t), \mathbf{z}(0)=\left[\begin{array}{cc} {0.3} & {0} \end{array}\right]^{\top}.
\end{equation}
The problem is to find the optimal states and control for the given time-delay system, which minimizes the quadratic performance index
\begin{equation} \label{Ex3.4} 
J=\tfrac{1}{2} \int _{0}^{t_{f} }\left\{Q[\bar{y}(t)-y(t)]^2+Ru^2(t)\right\} dt,
\end{equation}
where $y(t)\in{\mathbb{R}}$ is the output of system, $\bar{y}(t)\in{\mathbb{R}}$  is the reference input in which tracked by $y(t)$ and is given by \eqref{Ex3.2}, $Q=2$ and $R=1$. The time-delay and terminal time take different values for the following cases:

Case 1. $h_{x} =1,3,5,\,\,\hspace{1.027cm} t_{f} =20$,

Case 2. $h_{x} \hspace{.01cm} =3,5,9,15,30,\,\, t_{f} =60$.
\\

\begin{figure}[!ht]
\centering
\subfigure[Output and reference input]{
    \includegraphics[scale=.5]{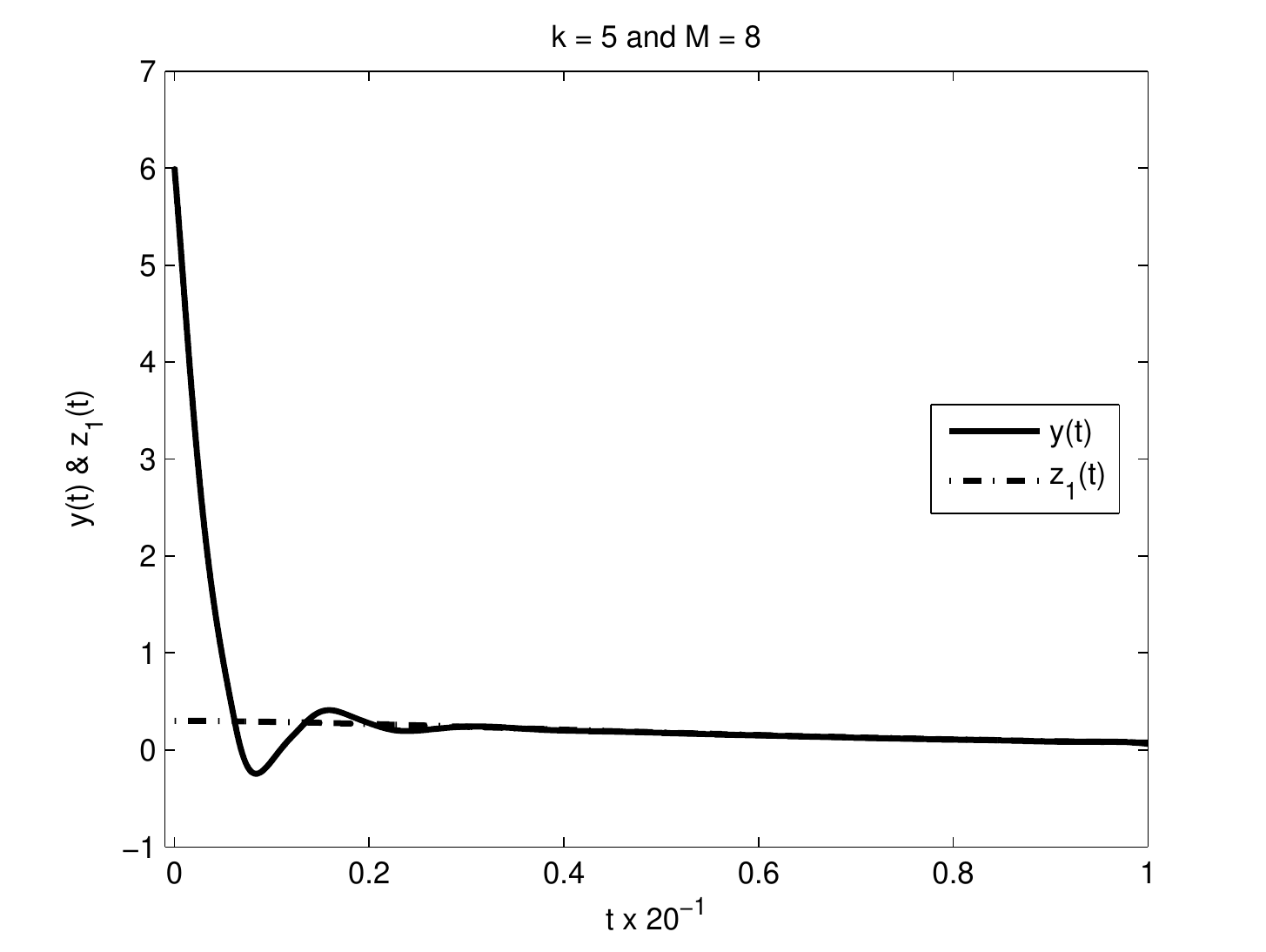}
    \label{Fig:2.2.1}
}
\subfigure[Output error]{
    \includegraphics[scale=.5]{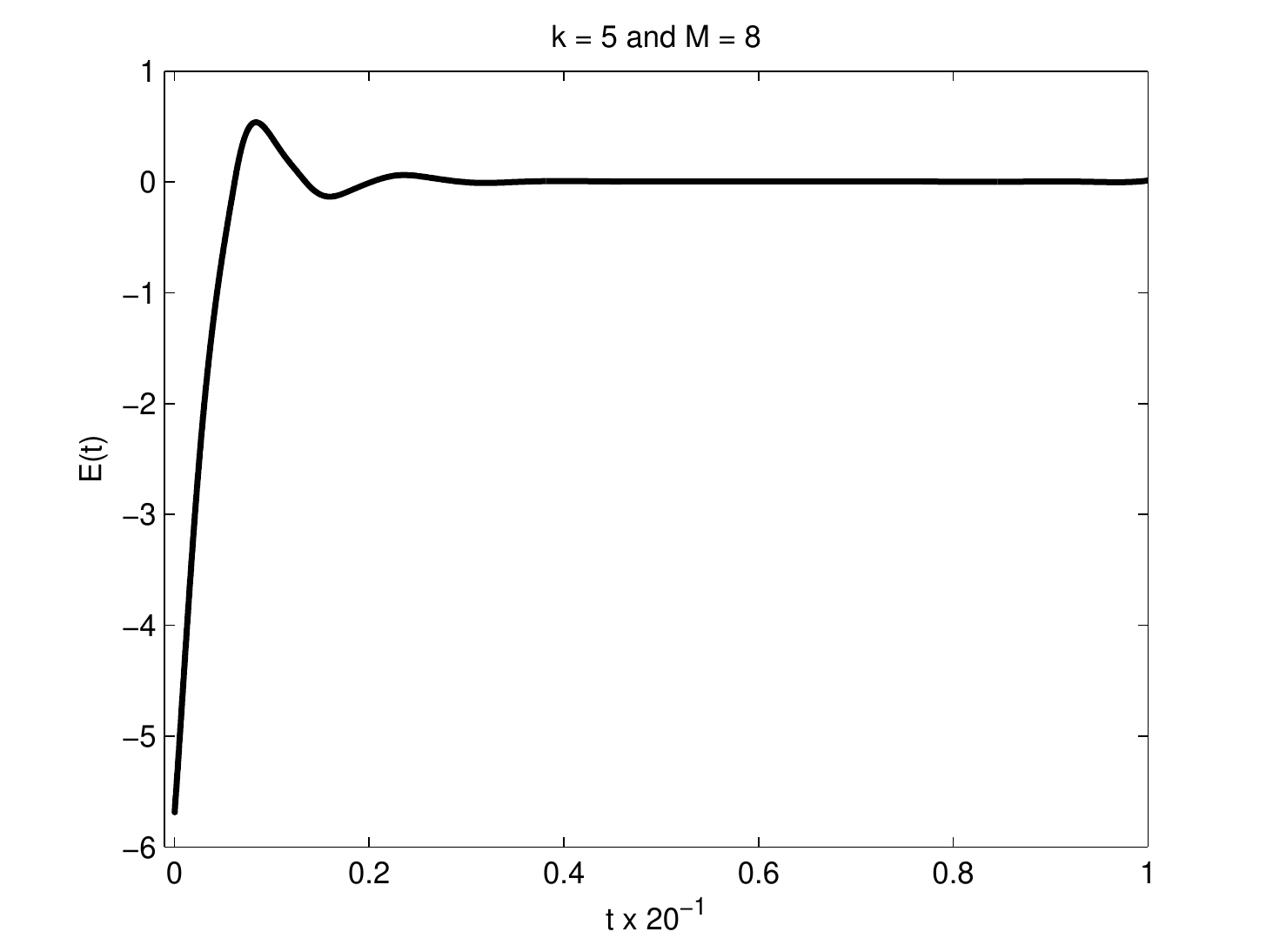}
    \label{Fig:2.2.2}
}
\caption[Optional caption for list of figures]{Output, reference and output error for Example 3; $h_{x}=1$, and $t_{f}=20$.}
\label{Fig:2.2}
\end{figure}

\begin{figure}[!ht]
\centering
\subfigure[Optimal states]{
    \includegraphics[scale=.5]{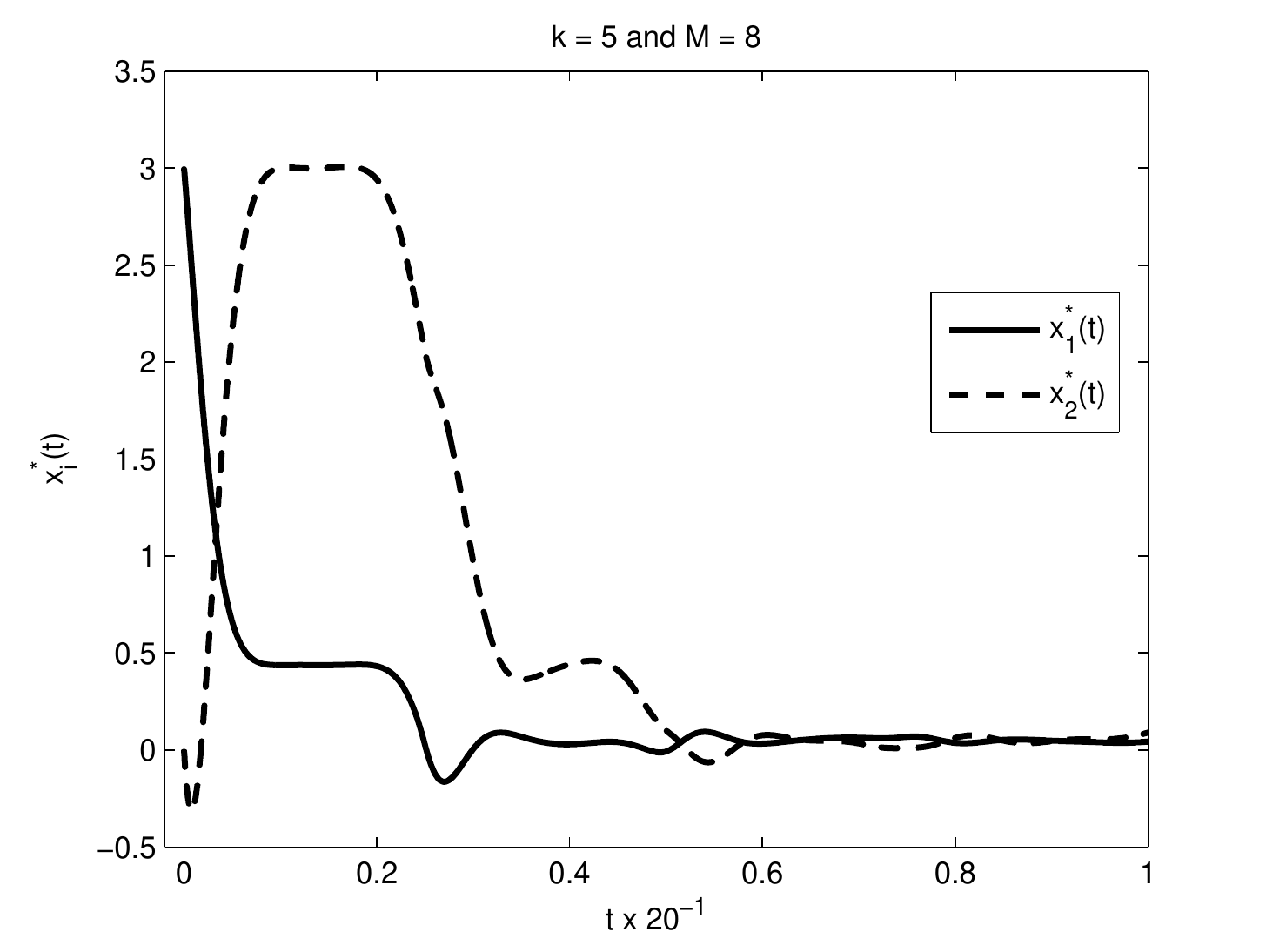}
    \label{Fig:3.1}
}
\subfigure[Optimal control]{
    \includegraphics[scale=.5]{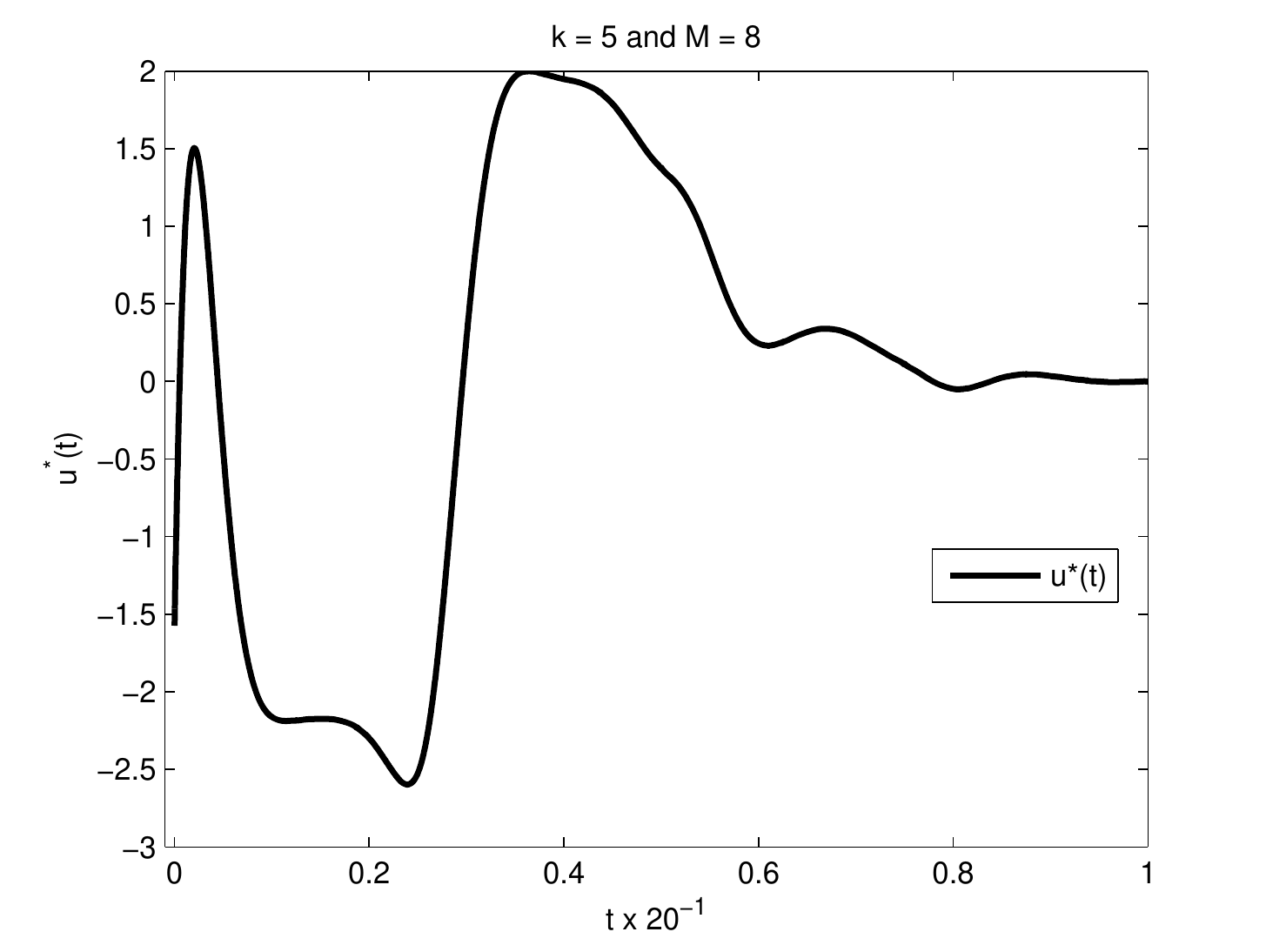}
    \label{Fig:3.2}
}
\caption[Optional caption for list of figures]{Optimal states and control for Example 3; $h_{x}=5$, and $t_{f}=20$.}
\label{Fig:3}
\end{figure}

\begin{figure}[!ht]
\centering
\subfigure[Optimal states]{
    \includegraphics[scale=.5]{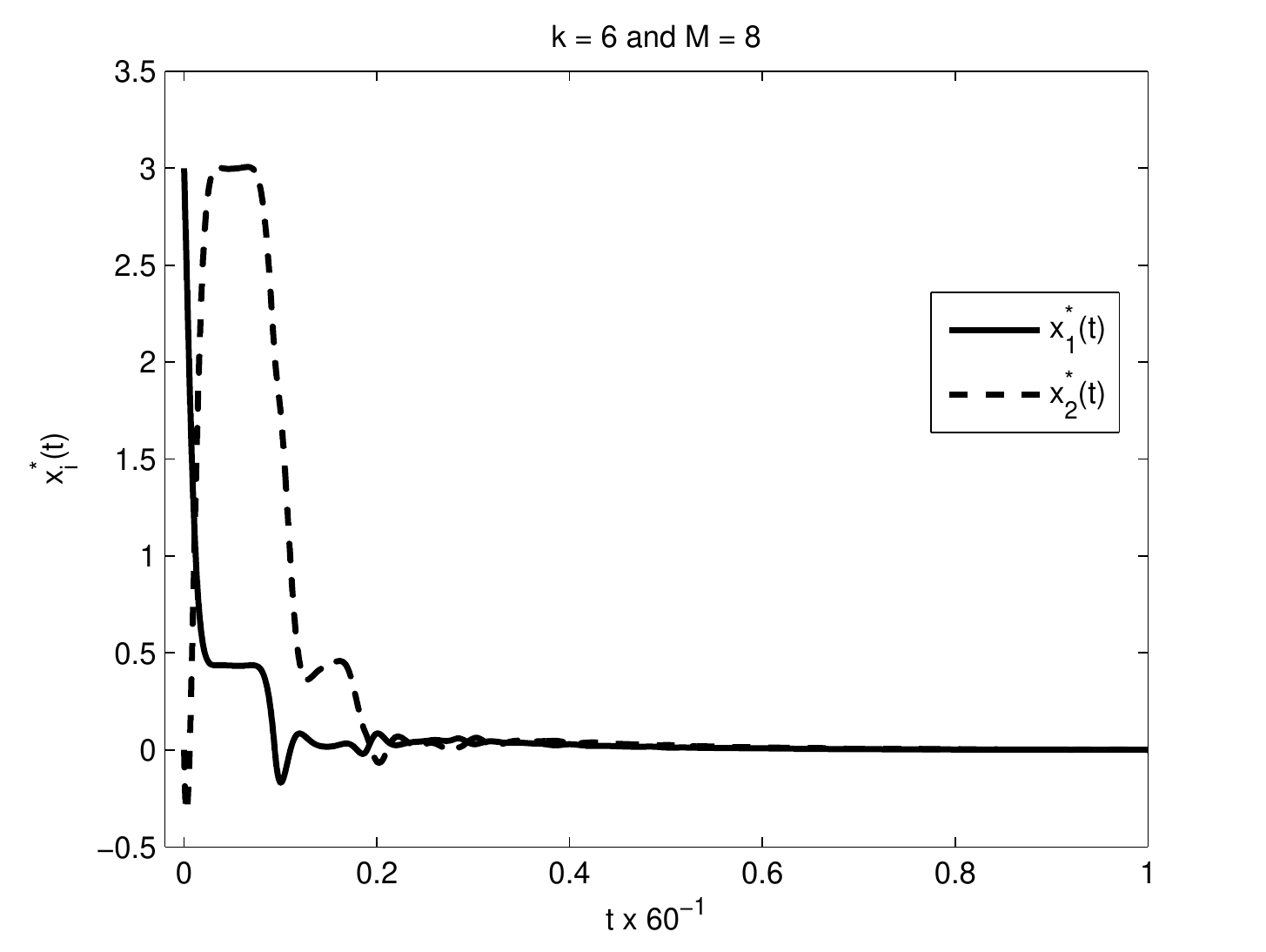}
    \label{Fig:4.1}
}
\subfigure[Optimal control]{
    \includegraphics[scale=.5]{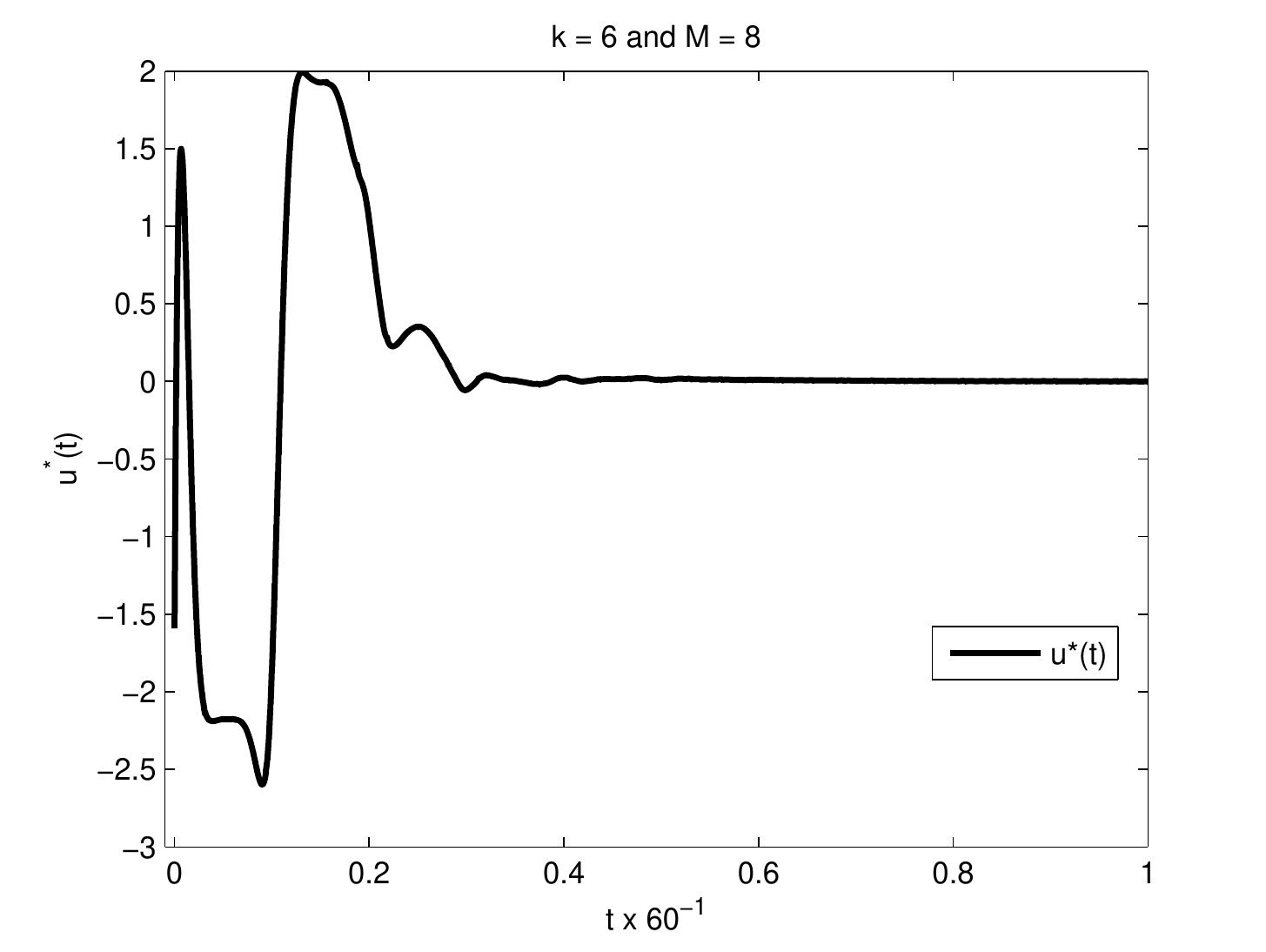}
    \label{Fig:4.2}
}
\caption[Optional caption for list of figures]{Optimal states and control for Example 3; $h_{x}=5$, and $t_{f}=60$.}
\label{Fig:4}
\end{figure}

If we set $\tau=t/t_{f}$, then the problem is converted to minimizing
\[J=\tfrac{1}{2}t_{f}  \int _{0}^{1}\bigg{\{}[\mathbf{x}(\tau)-\mathbf{r}(\tau)]^{\top}\begin{bmatrix} 8 & 0 \\ 0 & 0 \\ \end{bmatrix}[\mathbf{x}(\tau)-\mathbf{r}(\tau)]+u^{2}(\tau)\bigg{\}} d\tau \]
subject to the rescaled equations of \eqref{Ex3.0}--\eqref{Ex3.1}, where $\mathbf{r}(\tau)=\begin{bmatrix} 0.3{\text{e}}^{-0.1t_{f}\tau}-0.15{\text{e}}^{-0.2t_{f}\tau}& 0 \end{bmatrix}^{\top}$.
Using the proposed method, we solve the problem and obtain the results presented in Table ~\ref{tab:1}. The simulation curves of the obtained optimal states and control, the system output, the reference input, and the output error $E(t)$, where $E(t)=\bar{y}(t)-y(t)$, in the case $h_{x}=1$, $t_{f}=20$ are presented in Figs.\ref{Fig:2.2.1} and \ref{Fig:2.2.2}, respectively. Moreover, the graphs in Fig.\ref{Fig:3} and Fig.\ref{Fig:4} show $\mathbf{x}^*(t)$ and $ u^*(t)$ for $h_{x}=5$, $t_{f}=20$ and $h_{x}=5$, $t_{f}=60$, respectively. Table \ref{tab:1} illustrates the fact that the optimal performance index $J^*$ increases when the delay has increased and also value of the terminal time affects the value of the optimal index.

For $h_{x}=5$, and $t_{f}=20$ we get
\[
u^{*}(t)=\left\{\begin{array}{l} {-1.599503365t^{7}+6.939514532t^{6}-11.60025575t^{5}+4.164045166t^{4}+}\\{\,\,\hspace{1.5cm}18.68190416t^{3}-33.54410922t^{2}+18.08744678t-1.573550114,\hspace{.65cm}0 \le t<1.25}\\ {-1.579836463t^7+20.63373286t^6-114.3377350t^5+348.5948834t^4-}\\{\,\,\hspace{1.5cm}632.7083298t^3+687.1038586t^2-418.0542793t+110.1474518,\hspace{.35cm}1.25 \le t<2.5}\\ {- 1.584990728t^7+34.60298102t^6-322.7748087t^5+1667.517134t^4 -}\\{\,\,\hspace{1.5cm}5152.724483t^3 + 9523.371342t^2 - 9747.654972t + 4260.201064 ,\hspace{.35cm}2.5 \le t<3.75} \\ {- 1.580527130t^7 + 48.32648723t^6 - 632.1998092t^5 + 4587.257276t^4 -}\\{\,\,\hspace{1.5cm}19940.35153t^3 + 51928.97846t^2 - 75018.47153t + 46375.23574,\hspace{.55cm}3.75 \le t<5} \\ {- 1.535049521t^7+59.84697089t^6-998.0386619t^5+9228.228715t^4 -}\\{\,\,\hspace{1.5cm}51093.93715t^3+169398.1459t^2-311413.6496t+244898.9146
,\hspace{.55cm} 5 \le t<6.25}\\ {-1.611490389t^7+77.63429549t^6-1602.064754t^5+18357.38172t^4-}\\{\,\,\hspace{1.5cm}126143.9807t^3+519806.1911t^2-1189327.564t+1165540.148
,\hspace{.35cm}6.25 \le t<7.5} \\ {-1.597515699t^7+90.80070095t^6-2210.793157t^5+29889.94394t^4-}\\{\,\,\hspace{1.5cm}242350.0352t^3+1178423.616t^2-3181822.612t+3680123.555
,\hspace{.35cm}7.5 \le t<8.75} \\ {-1.617103941t^7+106.1787842t^6-2986.937676t^5+46666.73726t^4-}\\{\,\,\hspace{1.5cm}437325.8510t^3+2458209.746t^2-7674048.193t+10264009.91
,\hspace{.35cm}8.75 \le t<10} \\ {-1.642777755t^7+122.6187074t^6-3921.410142t^5+69653.27007t^4-}\\{\,\,\hspace{1.5cm}742130.0185t^3+4743069.046t^2-16836729.18t+25607790.77
,\hspace{.25cm}10 \le t<11.25} \\ {-1.586669620t^7+131.7788079t^6-4689.513324t^5+92690.80862t^4-}\\{\,\,\hspace{1.5cm}1098994.082t^3+7816336.188t^2-30877068.27t+52262383.11
,\, 11.25 \le t<12.5} \\ {-1.600070099t^7+146.9909776t^6-5786.158598t^5+126515.1074t^4-}\\{\,\,\hspace{1.5cm}1659475.150t^3+13057967.99t^2-57073349.02t+106890704.7
,\, 12.5 \le t<13.75} \\ {-1.570346408t^7+157.8118056t^6-6795.858796t^5+162560.2060t^4-}\\{\,\,\hspace{1.5cm}2332774.967t^3+20082606.41t^2-96035563.42t+196790673.1
,\hspace{.25cm}13.75 \le t<15} \\ {-1.593060270t^7+174.7227767t^6-8211.764765t^5+214386.4141t^4-}\\{\,\,\hspace{1.5cm}3357804.386t^3+31550825.23t^2-164679416.4t+368329684.0
,\hspace{.25cm}15 \le t<16.25} \\ {-1.597408464t^7+188.6101815t^6-9543.136587t^5+268224.4423t^4-}\\{\,\,\hspace{1.5cm}4522833.504t^3+45753885.01t^2-257114647.9t+619160087.7
,\hspace{.05cm}16.25 \le t<17.5} \\ {-1.596265471t^7+202.4942374t^6-11007.87894t^5+332416.2990t^4-}\\{\,\,\hspace{1.5cm}6022438.89t^3+65459683.76t^2-395241956.0t+1022669998.0
,\,17.5 \le t<18.75} \\ {-1.618309199t^7+219.4848422t^6-12756.63312t^5+411870.1854t^4-}\\{\,\,\hspace{1.5cm}7978130.478t^3+92716789.1t^2-598560876.0t+1655947901.0
,\hspace{.25cm}18.75 \le t\le 20.} \end{array}\right.
\]
\begin{figure}[!ht]
\centering
\subfigure[Optimal states]{
    \includegraphics[scale=.5]{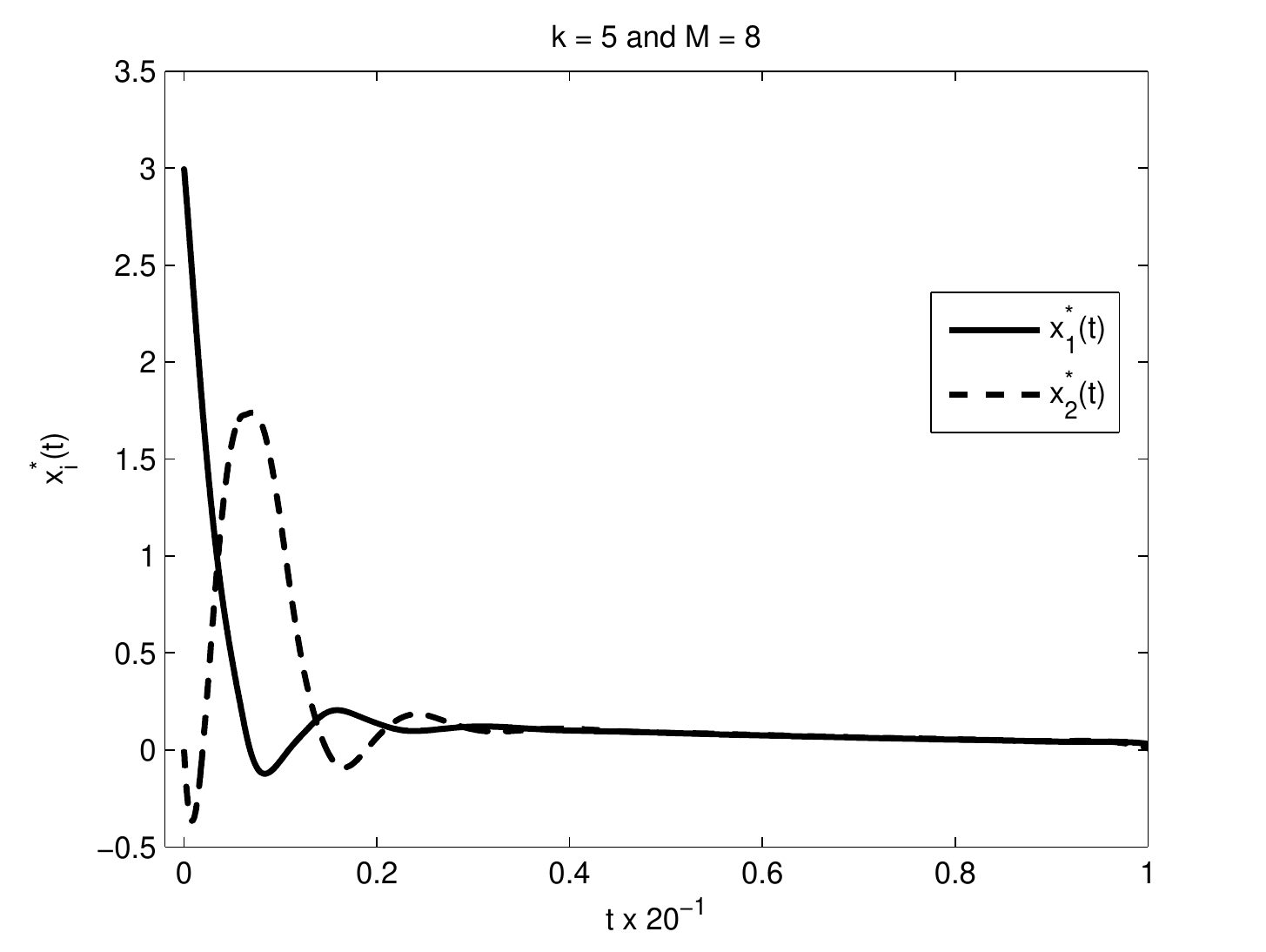}
    \label{Fig:2.1.1}
}
\subfigure[Optimal control]{
    \includegraphics[scale=.5]{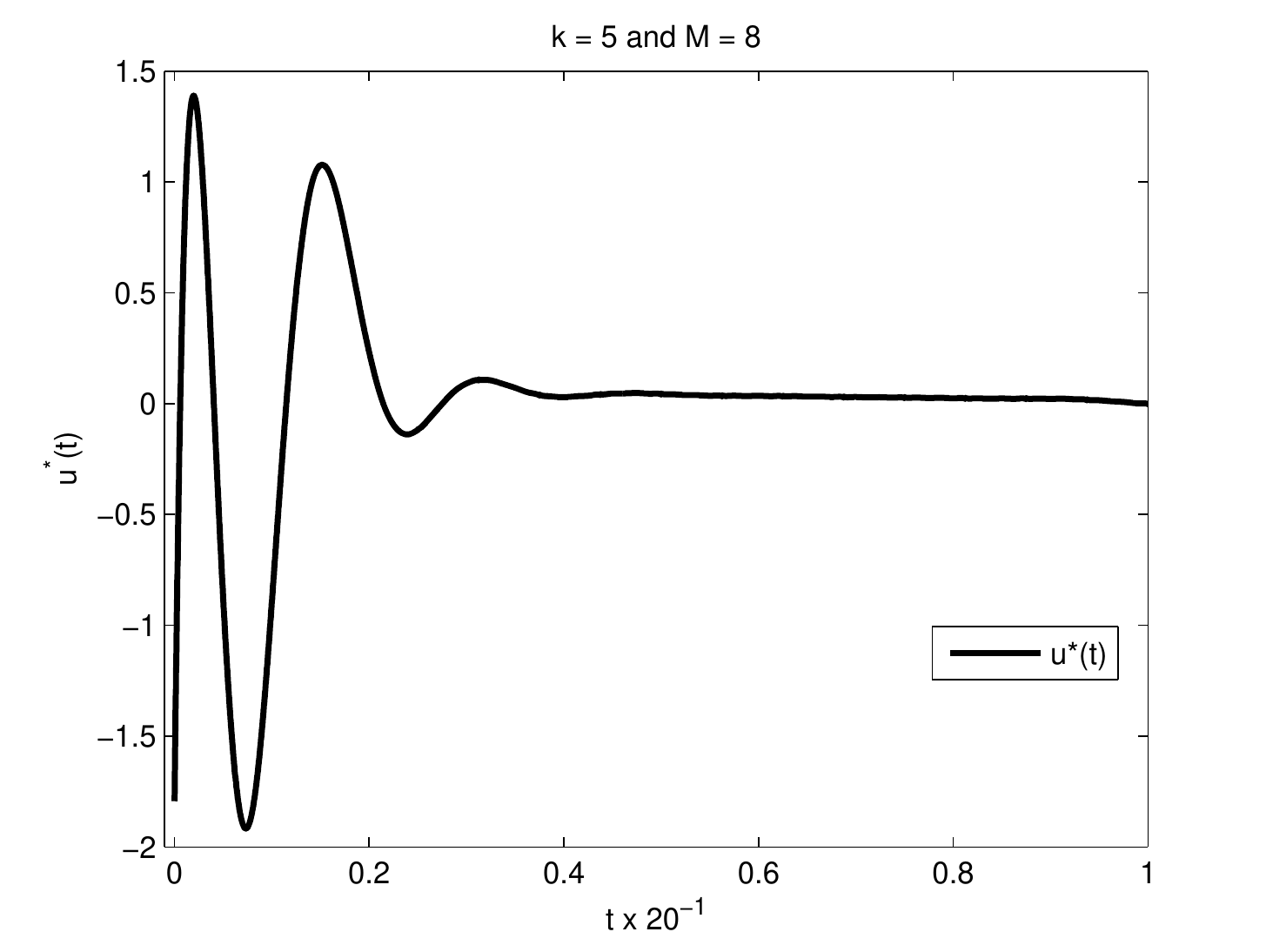}
    \label{Fig:2.1.2}
}
\caption[Optional caption for list of figures]{Optimal states and control for Example 3; $h_{x}=1$, and $t_{f}=20$.}
\label{Fig:2.1}
\end{figure}

\begin{table}[!ht]
\centering
\caption{optimal performance index for Example 3} \label{tab:1}
\begin{tabular}{cccccccccc}
\toprule
  & & $t_{f}=20$ & & & & & $t_{f}=60$ & \\\noalign{\vskip 1mm}
\cline{2-4}
\cline{6-10}
\noalign{\vskip 2mm} 
$h_{x}$ & 1 & 3 & 5 & & 3 & 5 & 15 & 30 & 50 \\
\noalign{\vskip 2mm}
$J^{*}$ &  12.5079 & 21.5463 & 30.2279 & & 24.4495 & 33.1185 & 76.8341 & 145.2474 & 161.9482\\
\bottomrule
\end{tabular}
\end{table}

Comparing the obtained curves of the output error and control with those were obtained in \cite{Tang.Li.Zhao}, we can see that the convergence rate of the presented method is higher than the convergence rate of the method in \cite{Tang.Li.Zhao}; however like the mentioned method, the proposed approach for systems with different delays and long time-delay is effective. The implementation of the method is very easy and convenient. The new optimal tracker can be successfully applied to the tracking systems regardless of the number of delays, and the number and the types of reference inputs and initial functions.

\subsection{Example 4}
Consider a non-square multi-input multi-output controllable and observable system (see \cite{Huang.Tsai.Provence.Shieh})
\begin{align}\label{Ex4.1}
\nonumber \dot{\mathbf{x}}(t)&=\begin{bmatrix} \hspace{.22cm}0.809 & -2.060 & \hspace{.22cm}0.325 & \hspace{.22cm}0.465 & \hspace{.22cm}0.895\\ \hspace{.22cm}6.667 & \hspace{.22cm}0.200 & \hspace{.22cm}1.333 & 0 & \hspace{.22cm}0.667\\ -1.291 & \hspace{.22cm}0.458 & -1.072 & -2.326 & -0.199\\-0.324 & \hspace{.22cm}0.824 & \hspace{.22cm}1.670 & -1.186 & -0.358\\ -3.509 & -4.316 & -0.702 & 0 & -8.351 \end{bmatrix}\mathbf{x}(t)\\ & +\begin{bmatrix}0 & 0 & 0 & 0 & 0 \\ -0.164 & 0 & 0 & 0 & 0 \\\hspace{.22cm} 0.729 & 0 & 0 & \hspace{.22cm}0.533 & -0.045 \\ 0 & 0 & 0 & -0.266 & \hspace{.22cm}0.167 \\ \hspace{.22cm}1.407 & 0 & 0 & 0 & -1.120 \end{bmatrix}\mathbf{x}(t-h_{x})+ \begin{bmatrix} \hspace{.22cm}0.955 & -0.379 \\ -1.667 & -1.667 \\ -0.212 & \hspace{.22cm}1.195 \\ \hspace{.22cm}0.618 & \hspace{.22cm}0.052 \\ \hspace{.22cm}0.877 & \hspace{.22cm}1.403 \end{bmatrix}\mathbf{u}(t-h_{u})
\end{align}
with the performance index
\begin{equation} \label{Ex4.2}
J=\int_{0}^{t_{f}} \left \{[\mathbf{y}(t)-\mathbf{r}(t)]^{\top}\mathbf{Q}[\mathbf{y}(t)-\mathbf{r}(t)]+\mathbf{u}^{\top}(t)\mathbf{I}_{2}\mathbf{u}(t) \right \} dt,
\end{equation} 
where
\[
\mathbf{y}(t)=\begin{bmatrix} 2 & 0 & 1 & 0 & 0 \\ 0 & 1.5 & 0 & 1.2 & 1 \end{bmatrix}\mathbf{x}(t),\,\mathbf{r}(t)=\begin{bmatrix} \sin(t) & \cos(t) \end{bmatrix}^{\top},\,\mathbf{x}(0)=\begin{bmatrix} 0.05 & 0.05 & 0.05 & 0.05 & 0.05 \end{bmatrix}^{\top}.
\]
$\mathbf{x}(t)\in{\mathbb{R}}^{5}$ is the state vector, $\mathbf{y}(t)\in{\mathbb{R}}^{2}$ is the output vector, and $\mathbf{u}(t)\in{\mathbb{R}}^{2}$ is the control vector. The problem is to find the optimal states and controls for the time-delay system \eqref{Ex4.1}, which minimizes \eqref{Ex4.2}.

We define $x_{6}(t)=2x_{1}(t)+x_{3}(t)$ and $x_{7}(t)=1.5x_{2}(t)+1.2x_{4}(t)+x_{5}(t)$. Also we set: $\mathbf{Q}=10^3\mathbf{I}_{2}$ and

(a). $h_x=0.15,\; h_u=0.05,\;\text{and}\; t_f=6$,\hspace{1.8cm} (b). $h_x=0.375,\; h_u=0.125,\;\text{and}\; t_f=8.$

\begin{figure}[!ht]
\centering
\subfigure[Optimal states]{
    \includegraphics[scale=.5]{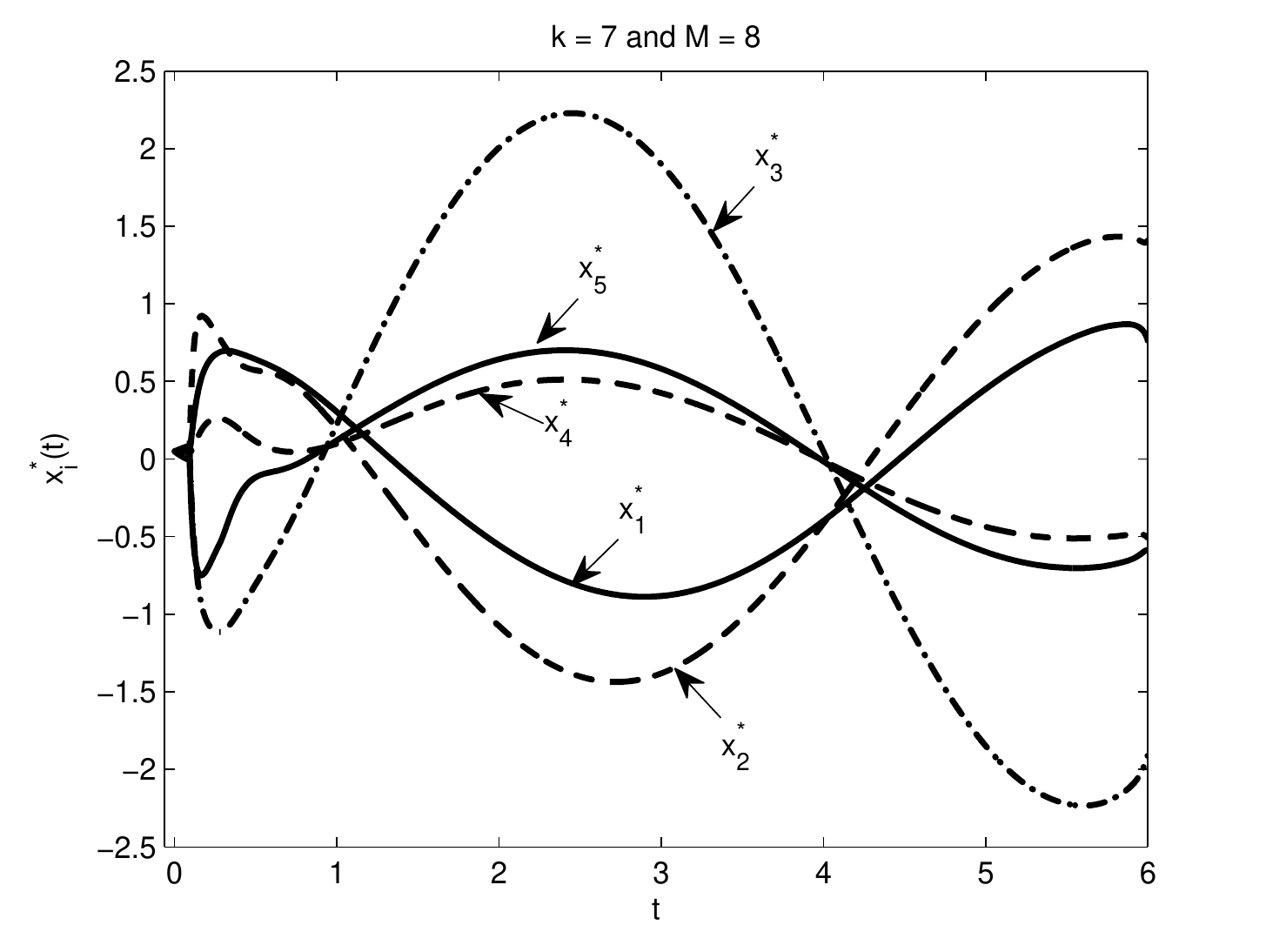}
    \label{Fig:5.1}
}
\subfigure[Transient response of Fig.\ref{Fig:5.1}]{
    \includegraphics[scale=.5]{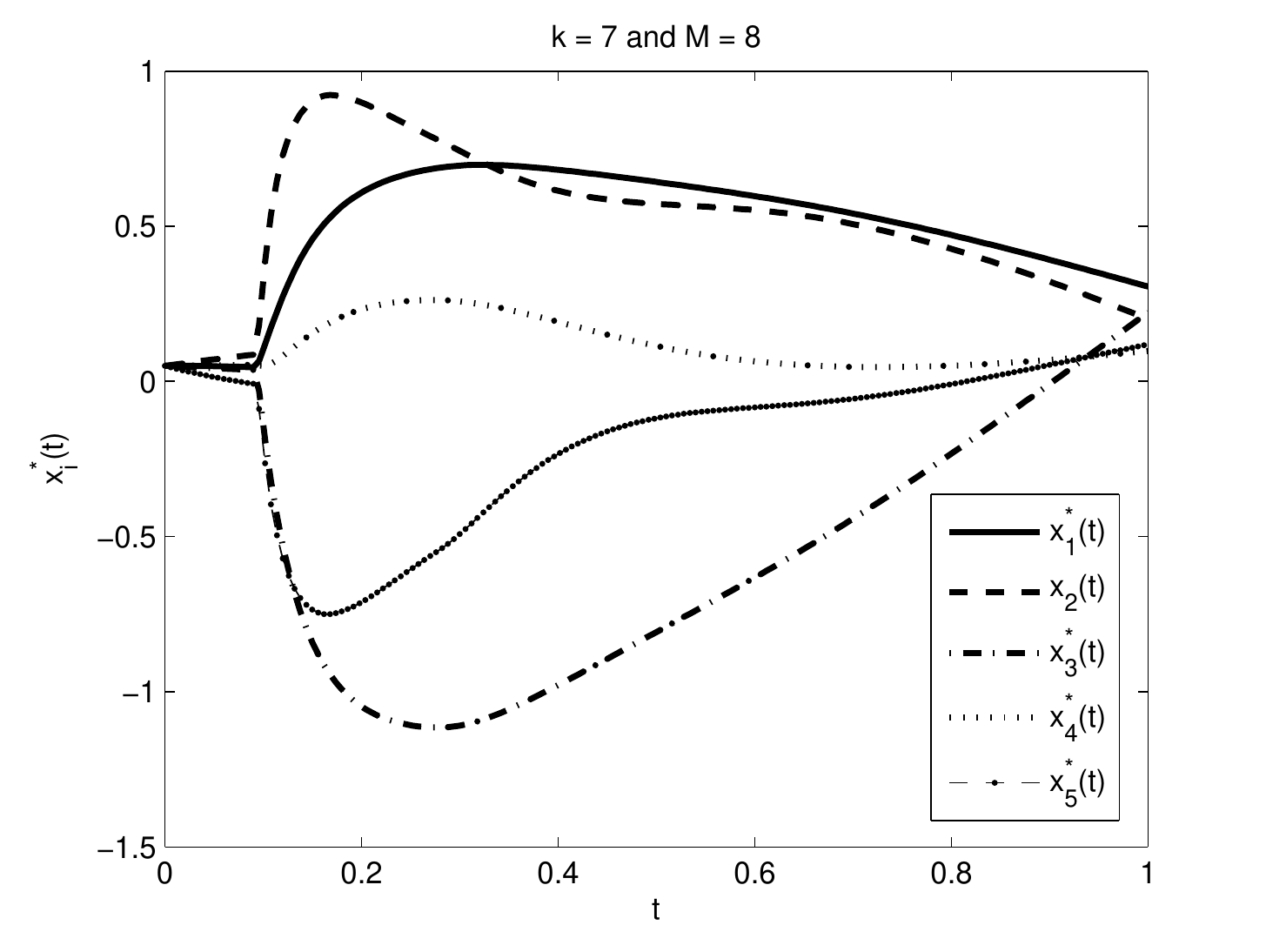}
    \label{Fig:5.2}
}
\caption[Optional caption for list of figures]{Optimal trajectories for Example 4(a).}
\label{Fig:5}
\end{figure}

\begin{figure}[!ht]
\centering
\subfigure[Optimal controls]{
    \includegraphics[scale=.5]{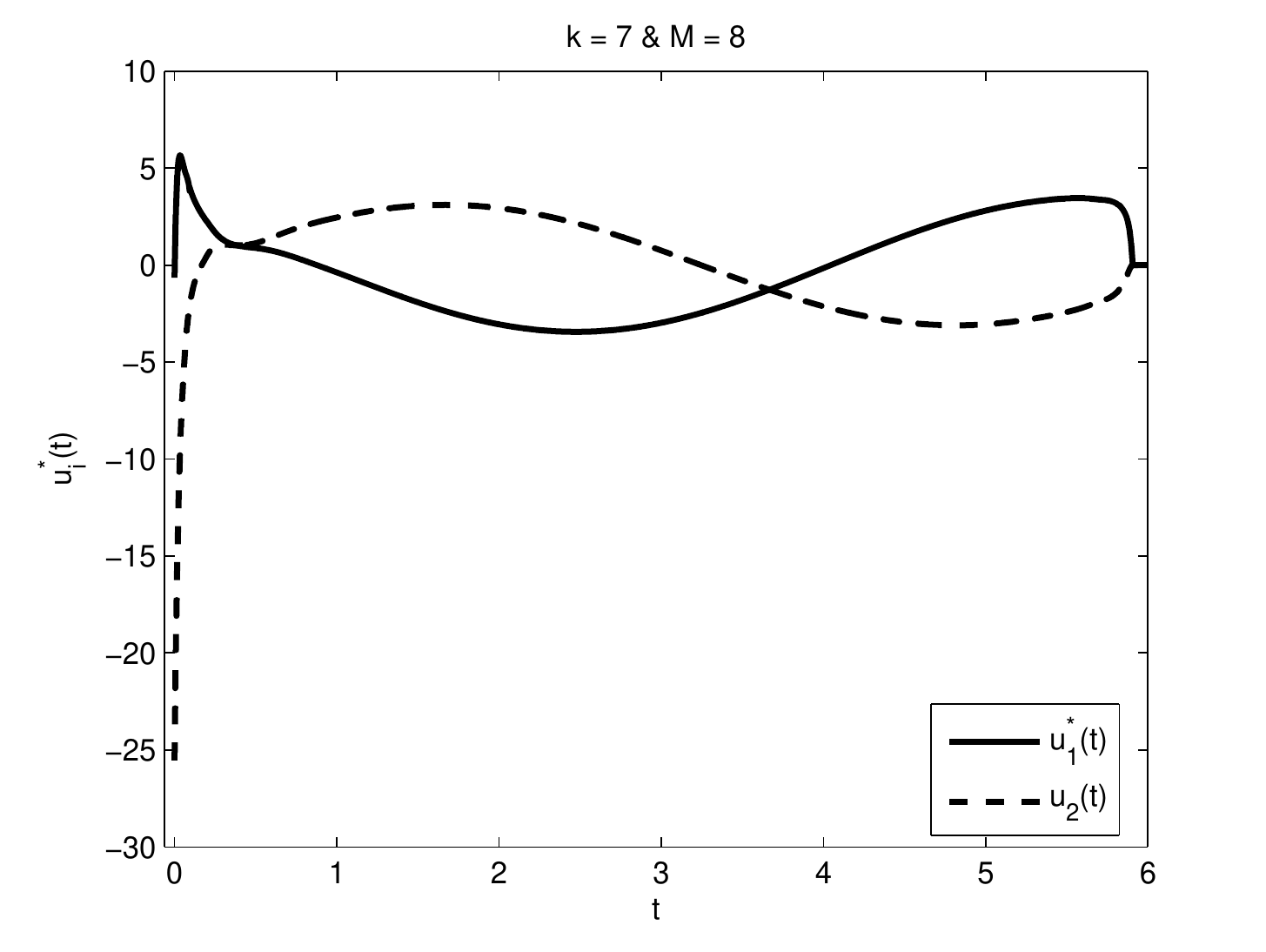}
    \label{Fig:6.1}
}
\subfigure[Outputs and reference inputs]{
    \includegraphics[scale=.5]{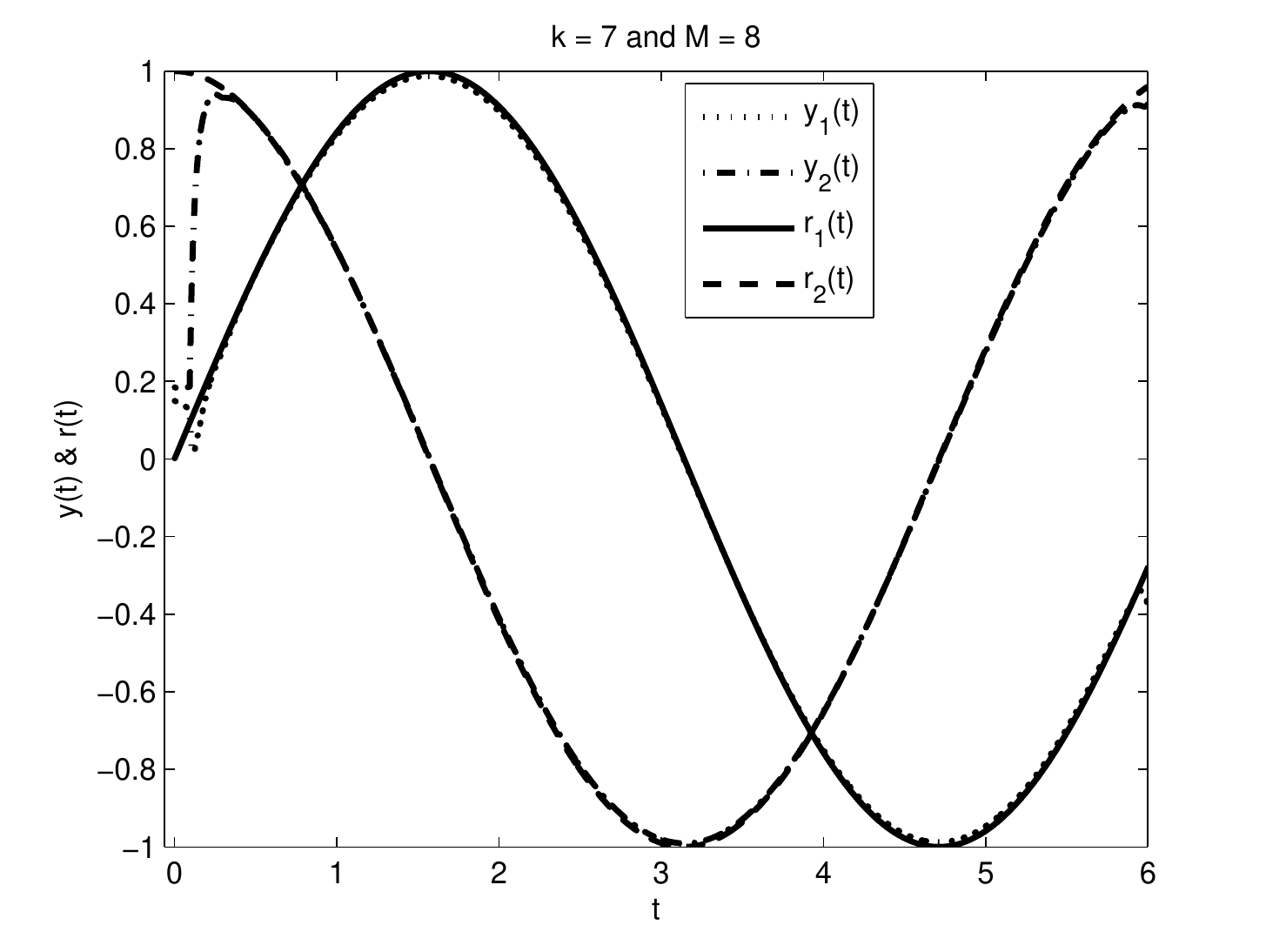}
    \label{Fig:6.2}
}
\caption[Optional caption for list of figures]{$\mathbf{u}^* (t)$, $\mathbf{y}(t)$ and $\mathbf{r}(t)$ for Example 4(a).}
\label{Fig:6}
\end{figure}

\begin{figure}[!ht]
\centering
\subfigure[Optimal states]{
    \includegraphics[scale=.5]{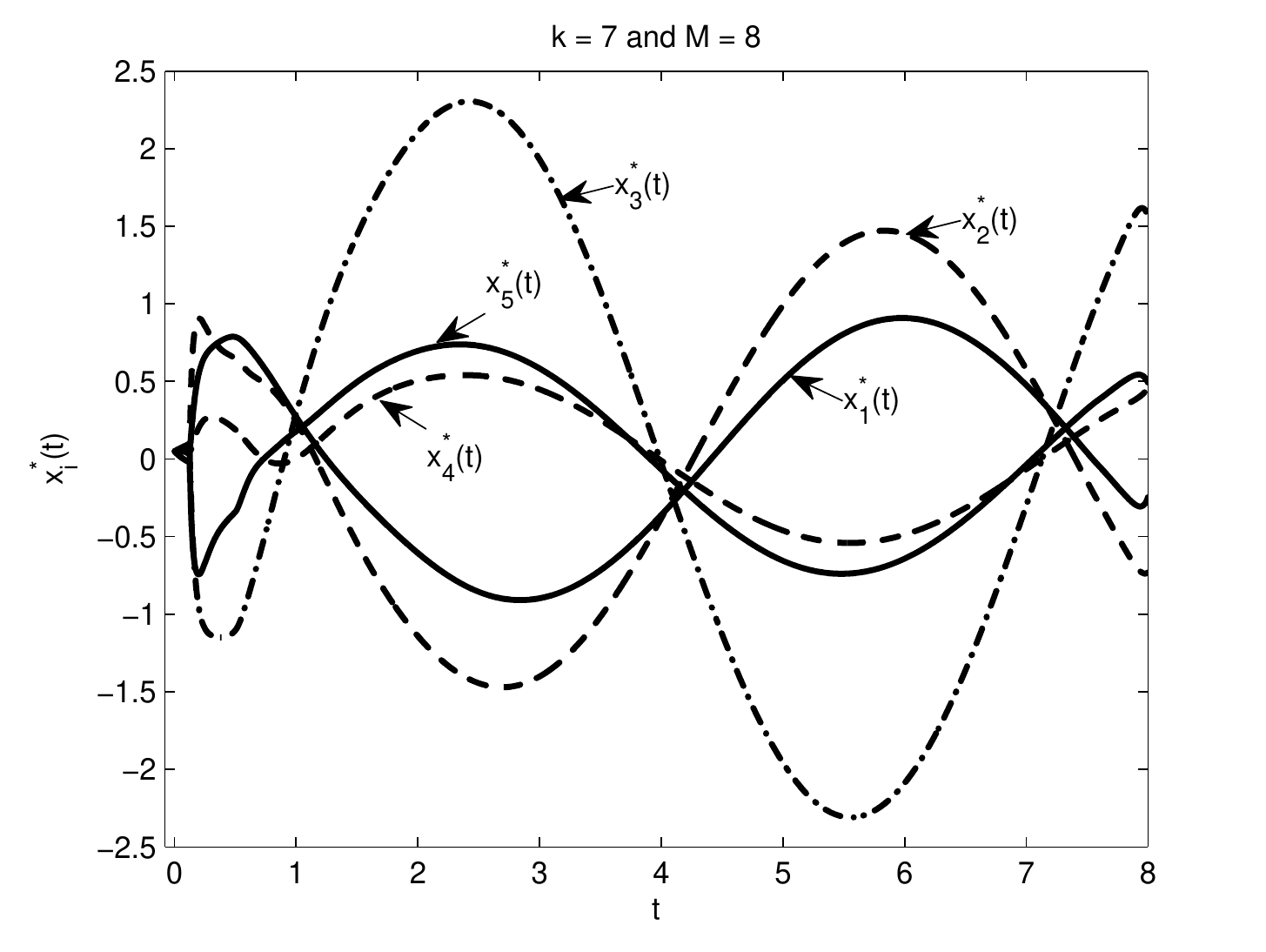}
    \label{Fig:7.1}
}
\subfigure[Optimal controls]{
    \includegraphics[scale=.5]{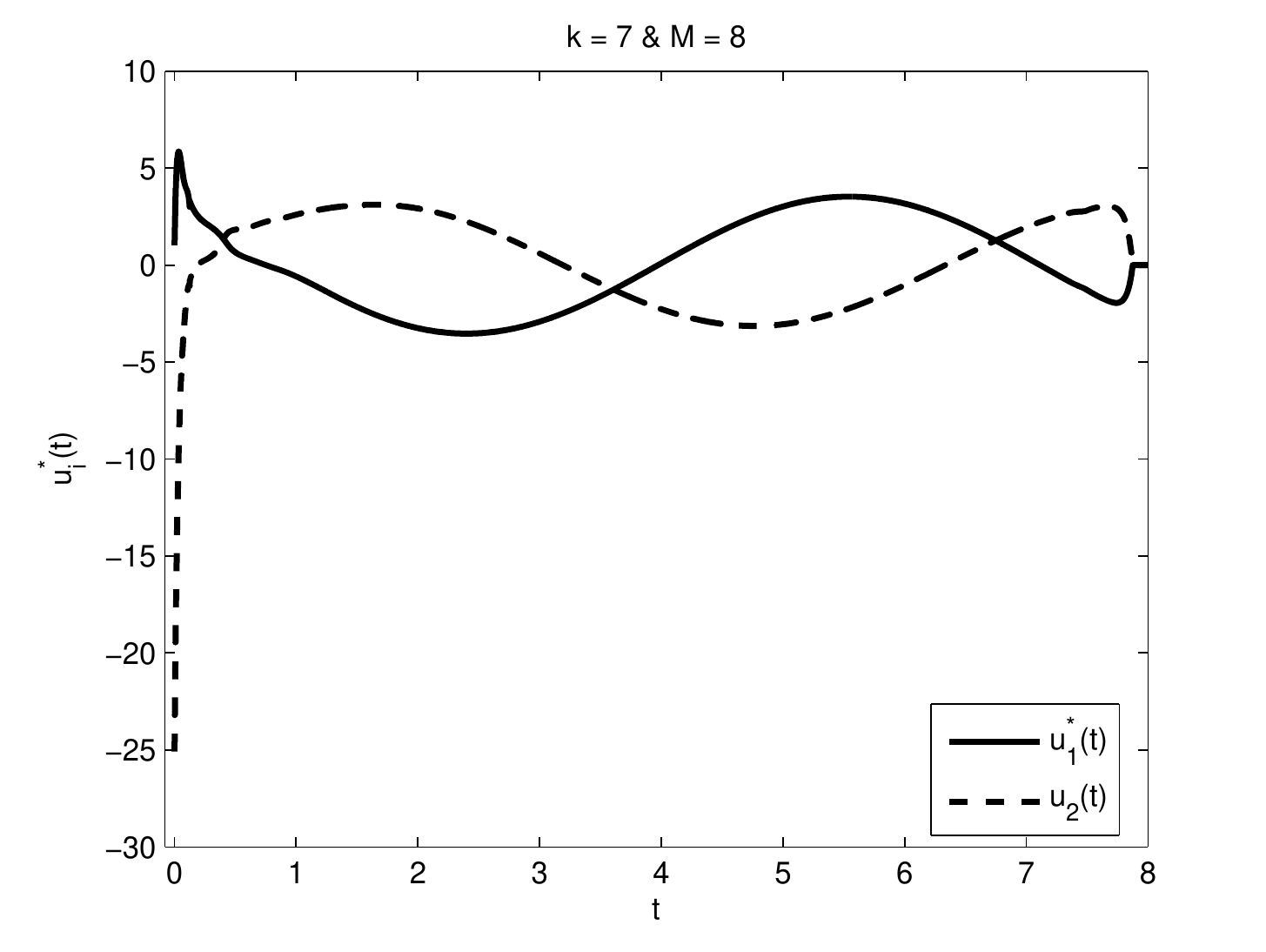}
    \label{Fig:7.2}
}
\caption[Optional caption for list of figures]{Optimal states and controls for Example 4(b).}
\label{Fig:7}
\end{figure}

\noindent Using these assumption, we get the optimal curves. The obtained results are shown in Figs.\ref{Fig:5}--\ref{Fig:7}.

\subsection{Example 5}
This example is adopted from \cite{Tsai.Wu.Lee.Guo.Su}. Consider a MIMO system
\begin{equation} \label{Ex5.1}
\dot{\mathbf{x}}(t)=\begin{bmatrix} -9.0 & \hspace{.22cm}4.0 & \hspace{.22cm}4.5 & -2.0 \\ -3.0 & \hspace{.22cm}0.4 & \hspace{.22cm}0.7 & -6.0\\ \hspace{.22cm}5.0 & \hspace{.22cm}0.3 & \hspace{.22cm}5.0 & \hspace{.22cm}3.0 \\ \hspace{.22cm}4.0 & -2.5 & \hspace{.22cm}2.0 & \hspace{.22cm}3.0  \end{bmatrix}\mathbf{x}(t)+\begin{bmatrix} \hspace{.22cm}1.0 & \hspace{.22cm}1.5 & \hspace{.22cm}0.0 \\ \hspace{.22cm}0.3 & \hspace{.22cm}2.0 & \hspace{.22cm}0.4 \\\hspace{.22cm} 0.3 & -0.3 & \hspace{.22cm}0.0 \\ -0.3 & -1.0 & \hspace{.22cm}0.5 \end{bmatrix}\mathbf{u}(t),
\end{equation}
\begin{equation} \label{Ex5.2}
\mathbf{y}(t)=\begin{bmatrix} \hspace{.2cm}1 & \hspace{.2cm}0 & \hspace{.2cm}2 & \hspace{.1cm}1 \\ -1 & \hspace{.2cm}1 & \hspace{.2cm}0 & -1\ \end{bmatrix}\mathbf{x}(t)+\begin{bmatrix} 0 & \hspace{.2cm}0 &  0 \\  0 & -1 & 0 \end{bmatrix}\mathbf{u}(t)
\end{equation}
with the performance index
\begin{equation} \label{Ex5.3}
J=\tfrac{1}{2}\int_{0}^{3}\left\{[\mathbf{y}(t)-\mathbf{r}(t)]^{\top}\mathbf{Q}[\mathbf{y}(t)-\mathbf{r}(t)]+\mathbf{u}^{\top}(t)\mathbf{R}\mathbf{u}(t) \right\} dt,
\end{equation}
where $\mathbf{x}(t)\in{\mathbb{R}}^{4}$, $\mathbf{y}(t)\in{\mathbb{R}}^{2}$, and $\mathbf{u}(t)\in{\mathbb{R}}^{3}$ are the state, output, and control vector, respectively. The initial conditions is $\mathbf{x}(0)=\begin{bmatrix} -0.25 & -0.5 & 0.25 & -0.3 \end{bmatrix}^{\top}$.
The desired input $\mathbf{r}(t)=\begin{bmatrix} r_1(t) & r_2(t) \end{bmatrix} ^{\top}$ is given by
\begin{equation}\label{Ex5.4}
r_{1}(t)=\left\{\begin{array}{ll} \cos(2\pi t), & 0 \le t <1 \\ 0.5t^{2}(1-t), & 1 \le t <2 \\ 0.5\,\cos(4\pi t)+1 & 2 \le t \le 3  \end{array}\right.\,\,\text{and}\,\,\,\,\,\,\,\,r_{2}(t)=\left\{\begin{array}{ll} 1.2t^{2}(1-t), & 0 \le t <1 \\ \cos(2\pi t),& 1 \le t <2 \\ 0.2\,\sin(4\pi t)-0.5, & 2 \le t \le 3. \end{array}\right.
\end{equation}
The problem is to find $\mathbf{x}^*(t)$ and $\mathbf{u}^*(t)$ which minimizes \eqref{Ex5.3} subject to the system \eqref{Ex5.1} and \eqref{Ex5.2} with the initial conditions and the reference input \eqref{Ex5.4}. We choose $\mathbf{Q}=10^4\mathbf{I}_{2}$ and $\mathbf{R}=\mathbf{I}_3$.

We can use the proposed algorithm in this problem with two following assumptions:
\begin{enumerate}
\item We define $x_{5}(t)=x_{1}(t)+2x_{3}(t)+x_{4}(t)$ and $x_{6}(t)=-x_{1}(t)+x_{2}(t)-x_{4}(t)-u_{2}(t)$, so the problem is reformulated to $\dot{\mathbf{x}}(t)=\mathbf{A}'\mathbf{x}(t)+\mathbf{B}'\mathbf{u}(t)+\mathbf{B}_{u}\dot{\mathbf{u}}(t)$,

\item First we set $2x_{3}(t)\rightarrow x_{3}(t)$ then ${x}_{2}(t)-u_{2}(t) \rightarrow {x}_{2}(t)$, therefore the problem is reformulated to $\dot{\mathbf{x}}(t)=\mathbf{A}''\mathbf{x}(t)+(\mathbf{B}''+\mathbf{A}''\mathbf{A}_{u})\mathbf{u}(t)-\mathbf{A}_{u}\dot{\mathbf{u}}(t)$.
\end{enumerate}
These two assumptions are equivalent, just we have to add $(\mathbf{I}_{2^{k-1}M} \otimes \mathbf{B}_{u}\mathbf{I}_r)$ or $(-\mathbf{I}_{2^{k-1}M} \otimes \mathbf{A}_{u}\mathbf{I}_r)$ to $\bm{\Lambda}_{12}$ and solve the problem. Using the first assumptions we can write
\[\mathbf{A}'=\begin{bmatrix} \,\mathbf{A} & \mathbf{0}_{4 \times 2} \\ \hspace{.2cm}\mathbf{A}_1 & \mathbf{0}_{1 \times 2} \\ \hspace{.2cm}\mathbf{A}_2 & \mathbf{0}_{1 \times 2} \end{bmatrix}, \mathbf{B}'=\begin{bmatrix} \mathbf{B} \\ \hspace{.2cm}\mathbf{B}_1 \\ \hspace{.2cm}\mathbf{B}_2 \end{bmatrix},\,\mathbf{Q}=\begin{bmatrix}\mathbf{0}_{4 \times 4} & \mathbf{0}_{4 \times 2}\\ \mathbf{0}_{2 \times 4} & 10^4\mathbf{I}_2 \end{bmatrix}\,\text{and} \;\mathbf{B}_u=\begin{bmatrix} \hspace{.2cm} \mathbf{0}_{5 \times 3} \\ \mathbf{B}_3 \end{bmatrix},\]
where $\mathbf{A}$ and $\mathbf{B}$ are defined in \eqref{Ex5.1} and
\[\mathbf{A}_1=\begin{bmatrix} 5.0 &2.1 & 16.5 & 7.0 \end{bmatrix},\,\mathbf{A}_2=\begin{bmatrix} 2.0 & -1.1 & -5.8 & -7.0 \end{bmatrix},\]
\[\mathbf{B}_1=\begin{bmatrix}1.3 & -0.1 & 0.5 \end{bmatrix},\,\mathbf{B}_2=\begin{bmatrix}-0.4 & 1.5 & -0.1\end{bmatrix},\,\text{and}\; \mathbf{B}_3=\begin{bmatrix}0 & -1 & 0\end{bmatrix}.\]

\begin{figure}[!ht]
\centering
\subfigure[Optimal states]{
    \includegraphics[scale=.5]{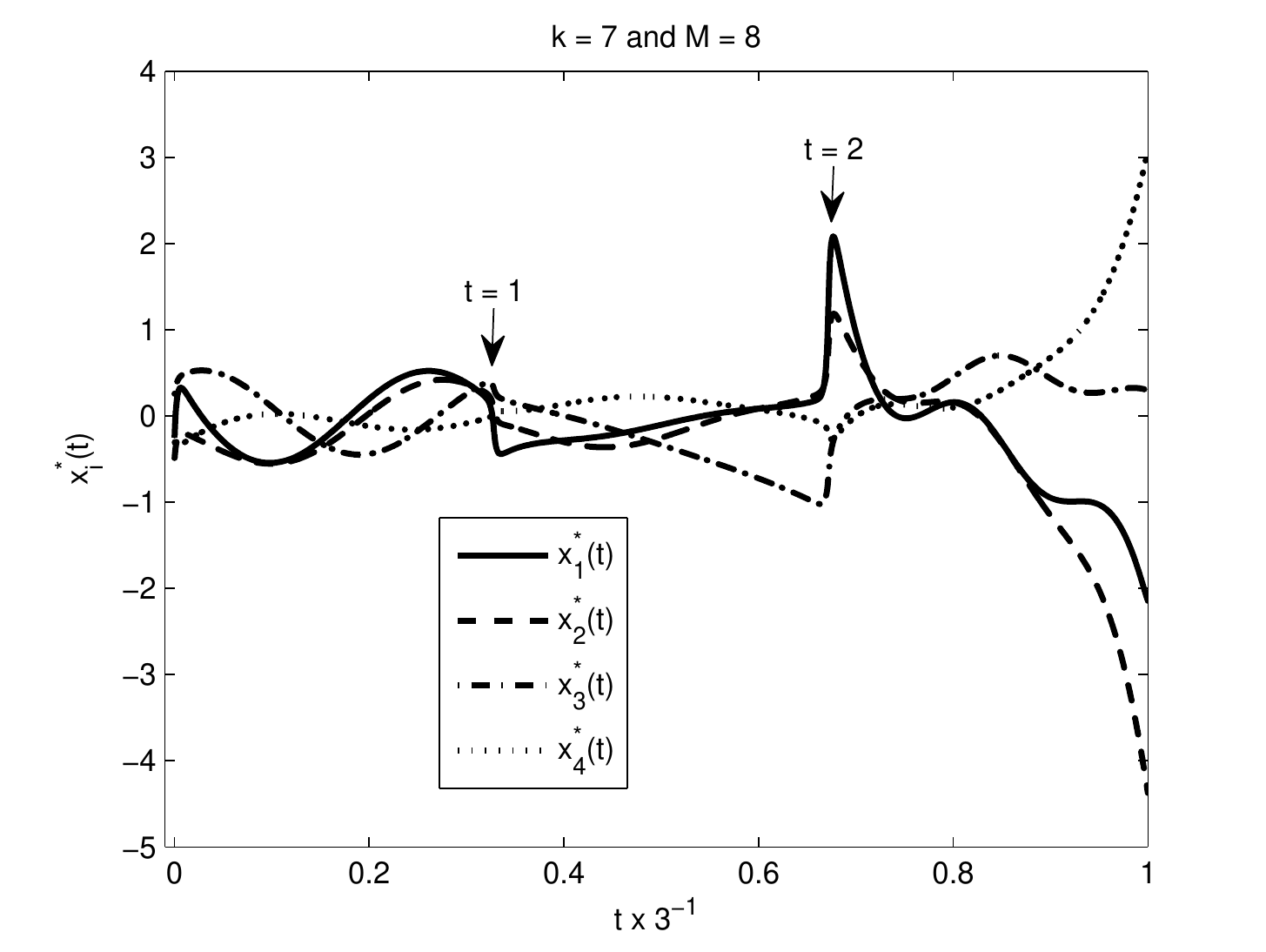}
    \label{Fig:8.1}
}
\subfigure[Optimal control $u^*_1(t)$]{
    \includegraphics[scale=.5]{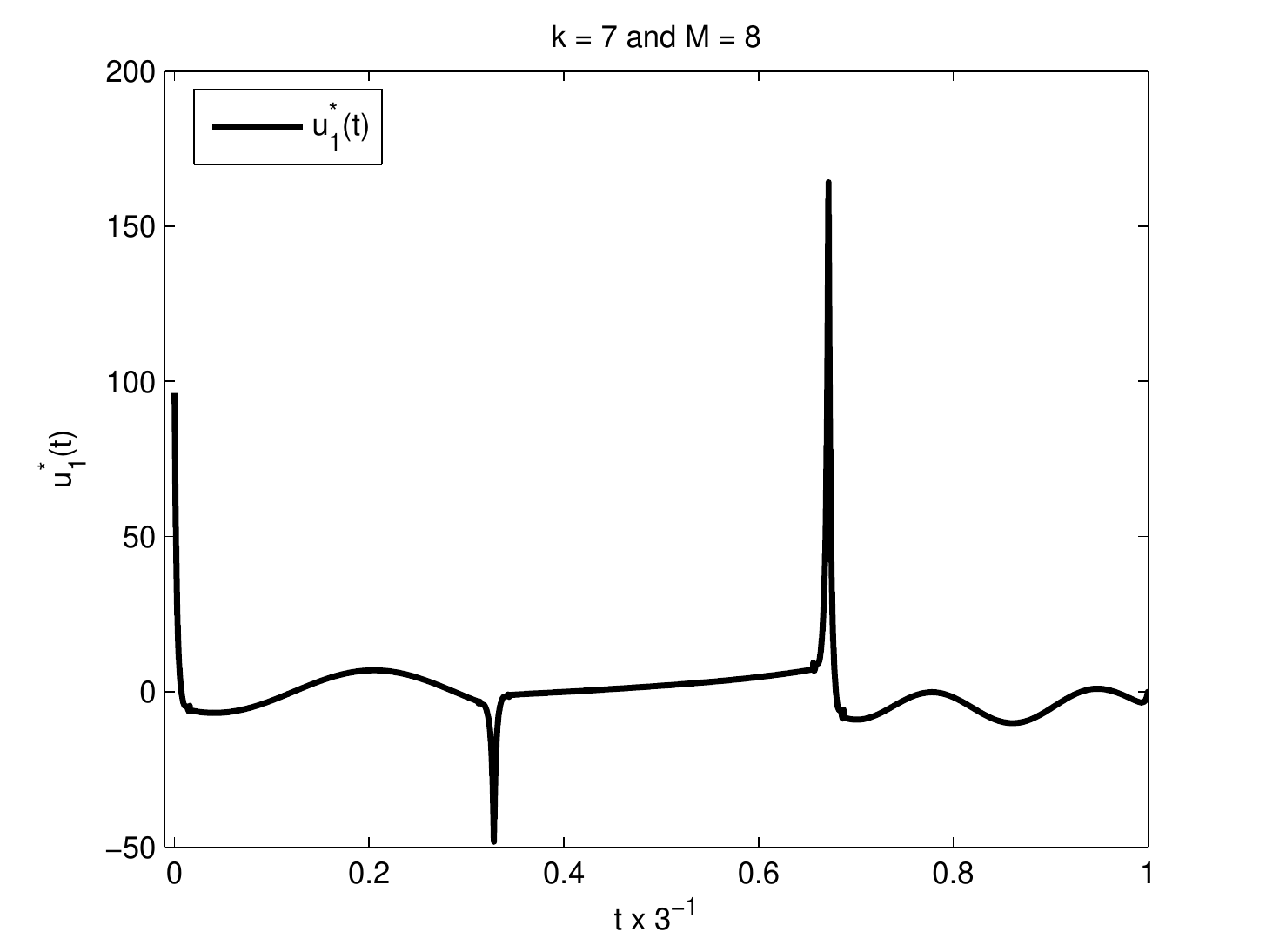}
    \label{Fig:8.2}
}
\subfigure[Optimal control $u^*_2(t)$]{
    \includegraphics[scale=.5]{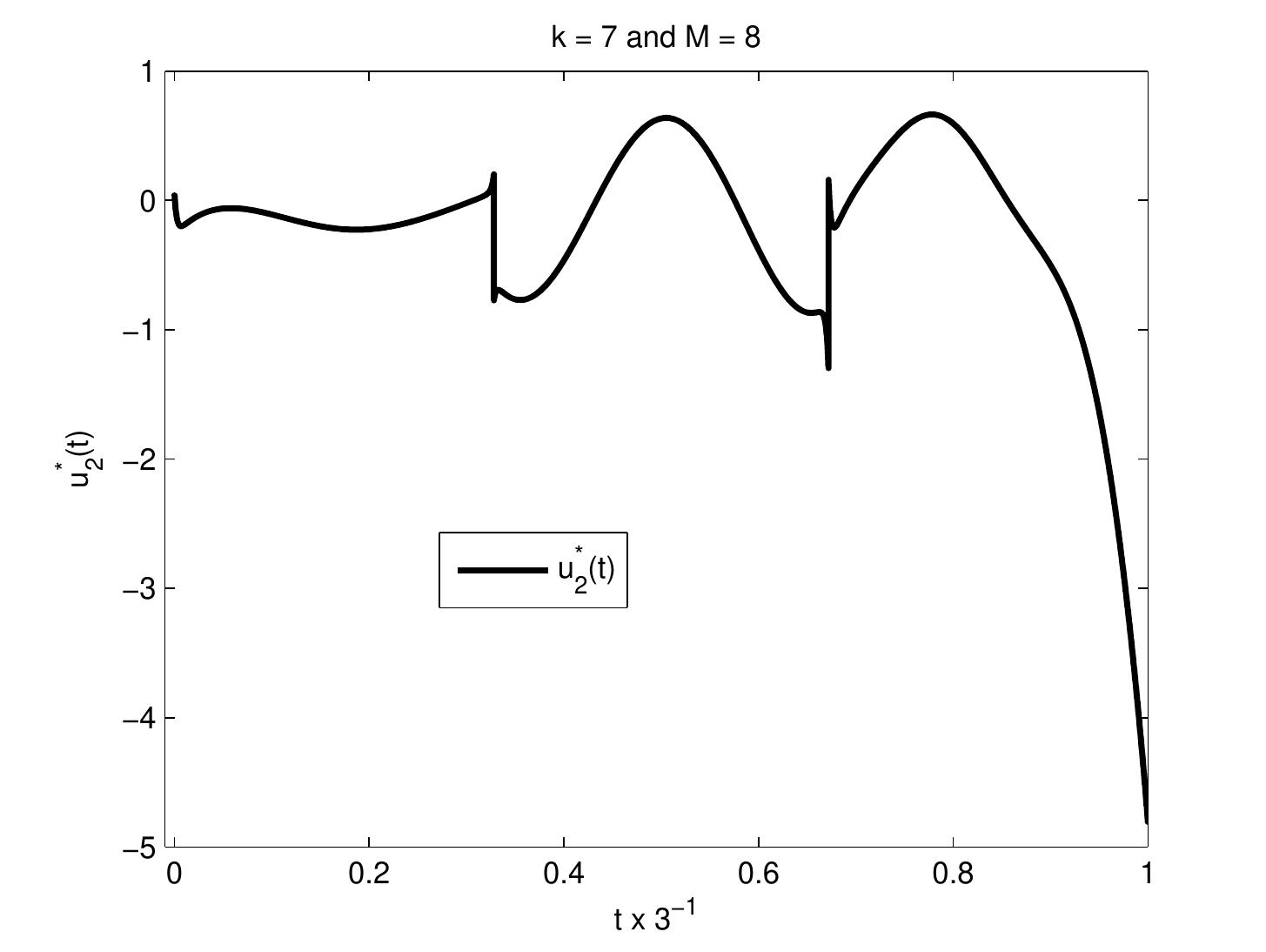}
    \label{Fig:8.3}
}
\subfigure[Optimal control $u^*_3(t)$]{
    \includegraphics[scale=.5]{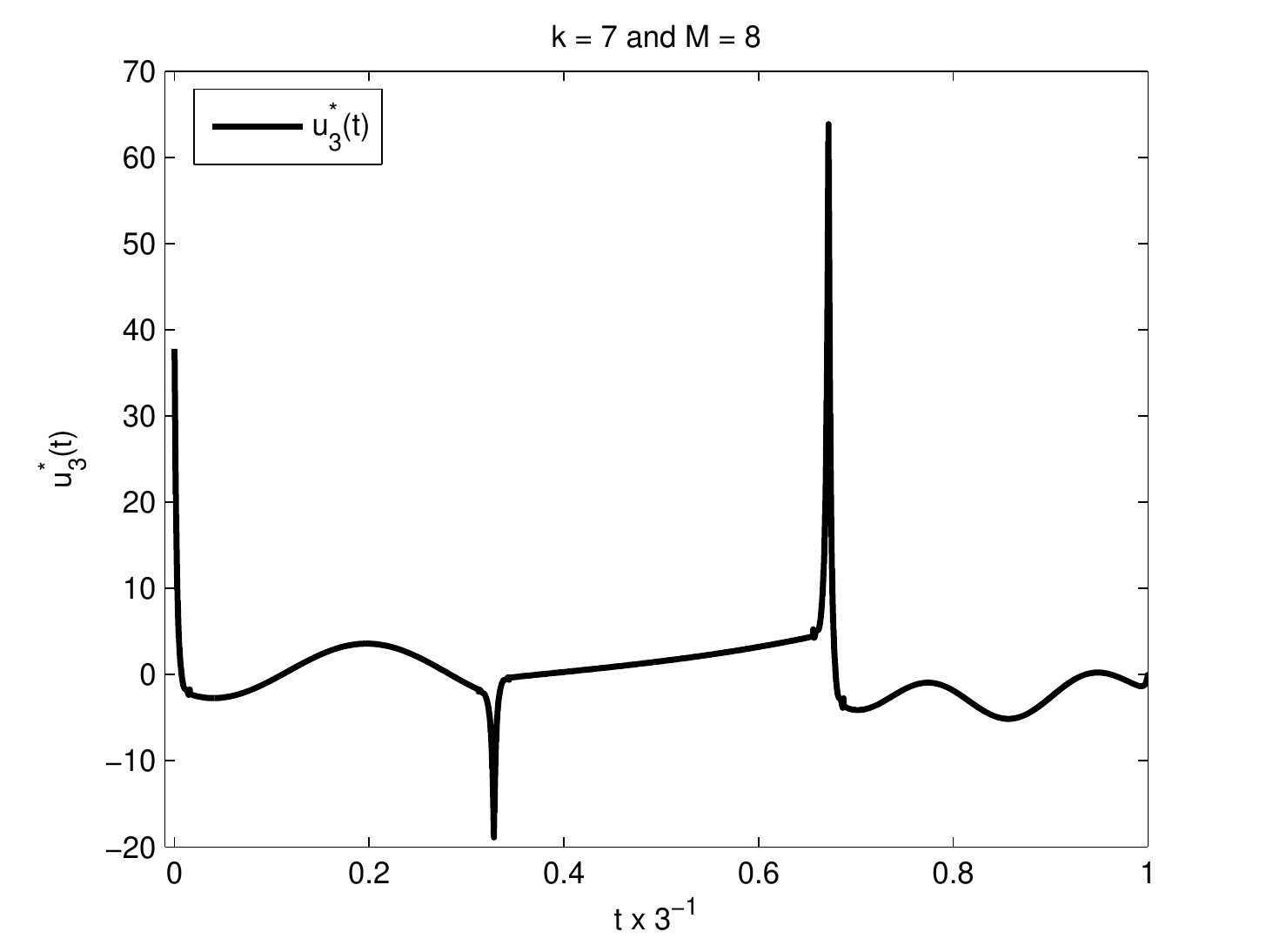}
    \label{Fig:8.4}
}
\caption[Optional caption for list of figures]{Optimal states and controls for Example 5.}
\label{Fig:8}
\end{figure}

\begin{figure}[!ht]
\centering
\subfigure[Output $y_1(t)$ and reference input $r_1(t)$]{
    \includegraphics[scale=.5]{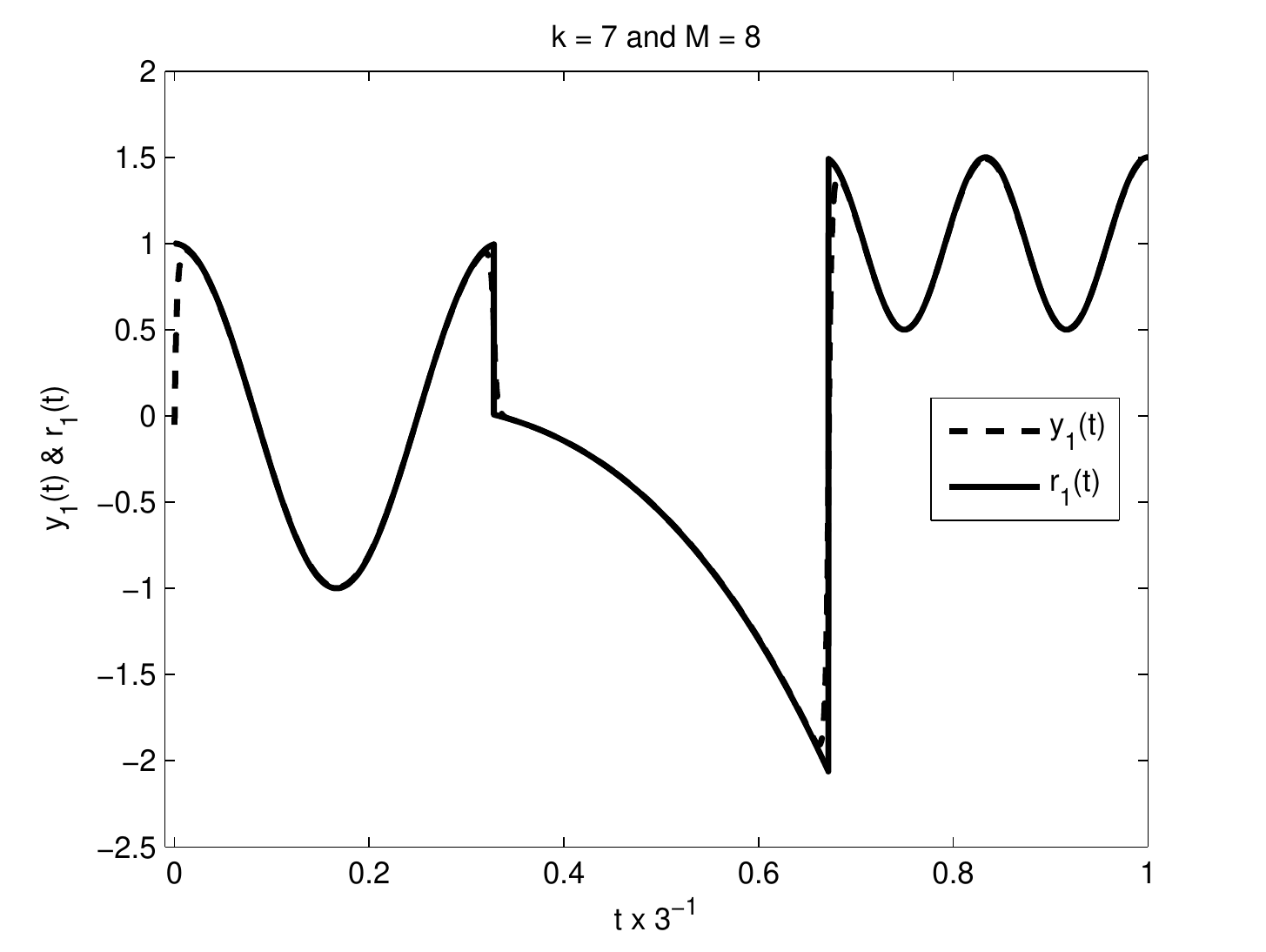}
    \label{Fig:9.1}
}
\subfigure[Output $y_2(t)$ and reference input $r_2(t)$]{
    \includegraphics[scale=.5]{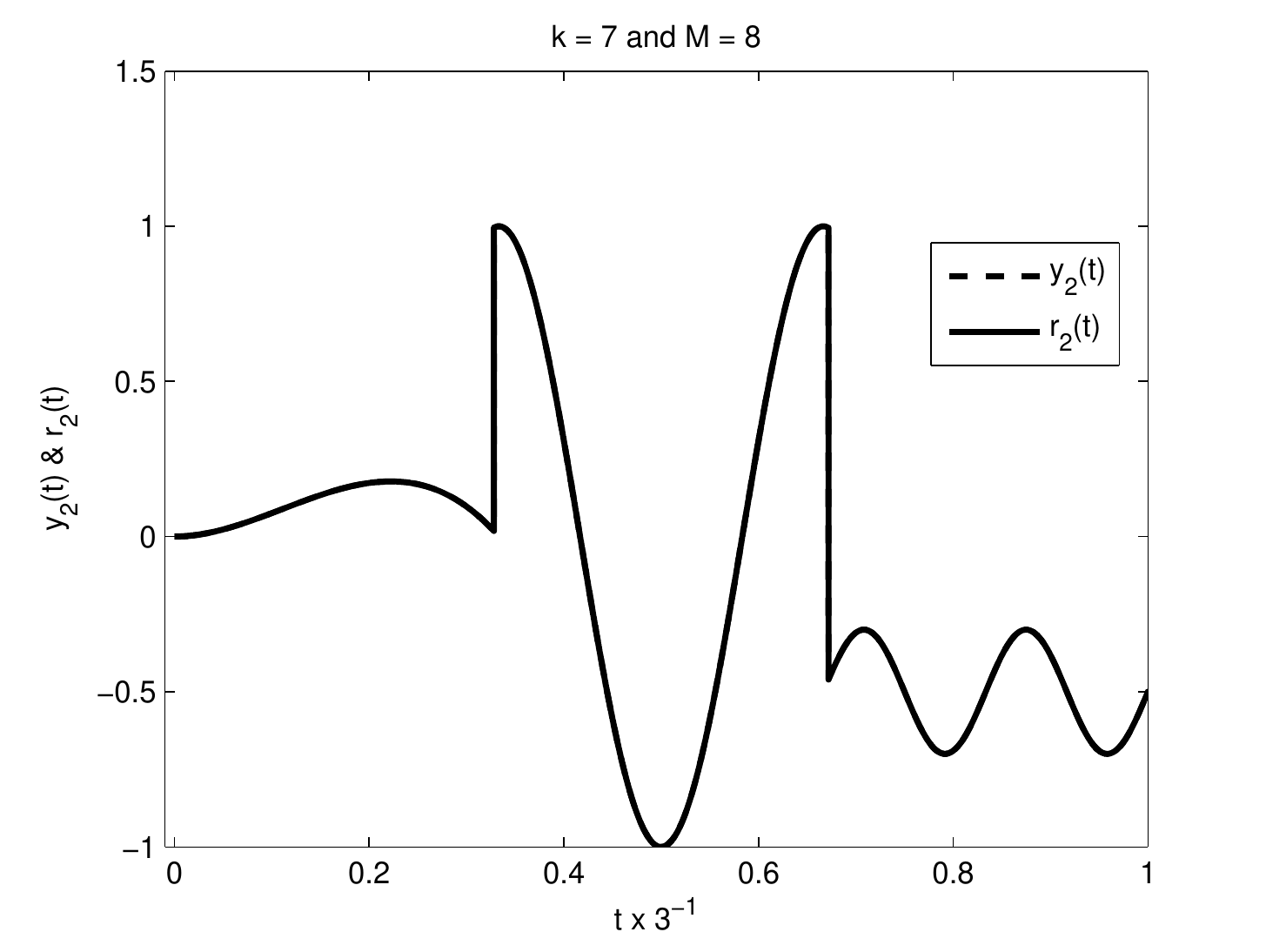}
    \label{Fig:9.2}
}
\caption[Optional caption for list of figures]{System output and reference input for Example 5.}
\label{Fig:9}
\end{figure}

Finally by selecting $k=7$ and $M=8$ we solve the transformed problem. The graphs of the optimal states and controls are given in Figs.\ref{Fig:8.1}--\ref{Fig:8.4}. Also Figs.\ref{Fig:9.1}--\ref{Fig:9.2} show the outputs of system and the reference inputs.
Since the desired input functions defined in \eqref{Ex5.4} are not continuous, we cannot use the compatibility constraint; also for this reason we observe two jumps in the optimal trajectories. Our algorithm implemented in MATLAB solves the problem within 6.436 seconds \footnote[5]{HP ENVY 15-j013cl Notebook PC} and gives 4608 computed parameters with $J^*=330.4858$.

\subsection{Example 6}
\noindent Consider a linear time-varying time-delay system described by
\begin{equation} \label{Ex6.1}
\dot{\mathbf{x}}(t)= \begin{bmatrix} {0} & {1} & {0} \\ {0} & {0} & {1} \\ {\cos t} & {0} & {0} \end{bmatrix} \mathbf{x}(t)+ \begin{bmatrix} {0} & \hspace{-2mm}{-1} & {0} \\ {-0.1t^2} & {0} & {0.5} \\ {{\text{e}}^{-t} } & {0} & {t} \end{bmatrix} \mathbf{x}(t-h_{x} )+ \begin{bmatrix} {0} \\ {0} \\ {2+\sin t} \end{bmatrix} u(t)\, ,\, \, 0\le t\le t_{f}
\end{equation}
\begin{equation} \label{Ex6.2} 
\mathbf{x}(t)= \begin{bmatrix} {1} & {0} & {\sin t} \end{bmatrix} ^{\top},\, \, -h_{x} \le t\le 0.
\end{equation} 
This system is to be controlled to minimize the performance index
\begin{equation}\label{Ex6.3}
J=[x_{1}(t_{f})-r(t_{f})]^{2}+\tfrac{1}{2}\int _{0}^{t_{f}}\left\{100[x_{1}(t)-r(t)]^{2}+u^{2} (t)\right\} dt
\end{equation}
in order that the state $x_{1}(t)$ tracks the desired trajectory $r(t)$, where $r(t) \in {\mathbb{R}}$ is
\[r(t)= \cos t.\]
The terminal time is $t_{f} =4$. In the following, we consider this optimal control problem with different time delays as case 1 and also with different constraints as case 2; we take:

Case 1: $\left\{\begin{array}{l}{{\text{a}}.\; h_{x} =0.5.}\\{{\text{b}}.\; h_x=1.0.}\\{{\text{c}}.\;\, h_x=2.0.}\end{array}\right.$

Case 2: $\left\{\begin{array}{l}{{\text{a}}.\;  h_{x} =0.5; x_{2}(h_{x})=-0.5\; \text{and} \; x_{3}(h_{x})=-1.5.}\vspace{1mm}\\{{\text{b}}.\; h_x=1.0; x_{2}(h_{x})=-1, x_{3}(h_{x})=-1\; \text{and} \; x_{3}(t_{f})=r(t_{f}).}\vspace{1mm}\\{{\text{c}}.\; h_x=2.0; \left\{\begin{array}{lll} x_{3}(t) \le r(t) & \text{when} & t \in [0\;\;\, h_{x}] \\ x_{2}(t) \le r(t) & \text{when} & t \in [h_{x}\,\,t_{f}] \end{array} \right.\; \text{and} \; x_{3}(t_{f})=0.}\vspace{1mm} \\{{\text{d}}. \, h_x=2.0;\left\{\begin{array}{lll} 0.0625t^{2}x_{2}(t)+(-0.05t+1)x_{3}(t)-u(t) \le 0.8 & \text{when} & t \in [0\;\;\, h_{x}] \\ x_{2}(t) \le r(t) & \text{when} & t \in [h_{x}\; t_{f}] \\ u(t) \le 0.5 & \text{when} & t \in [0\hspace{2.3mm} t_{f}].  \end{array} \right.} \end{array}\right.$\vspace{1mm}

\begin{figure}[!ht]
\centering
\subfigure[Optimal states and reference trajectory]{
    \includegraphics[scale=.5]{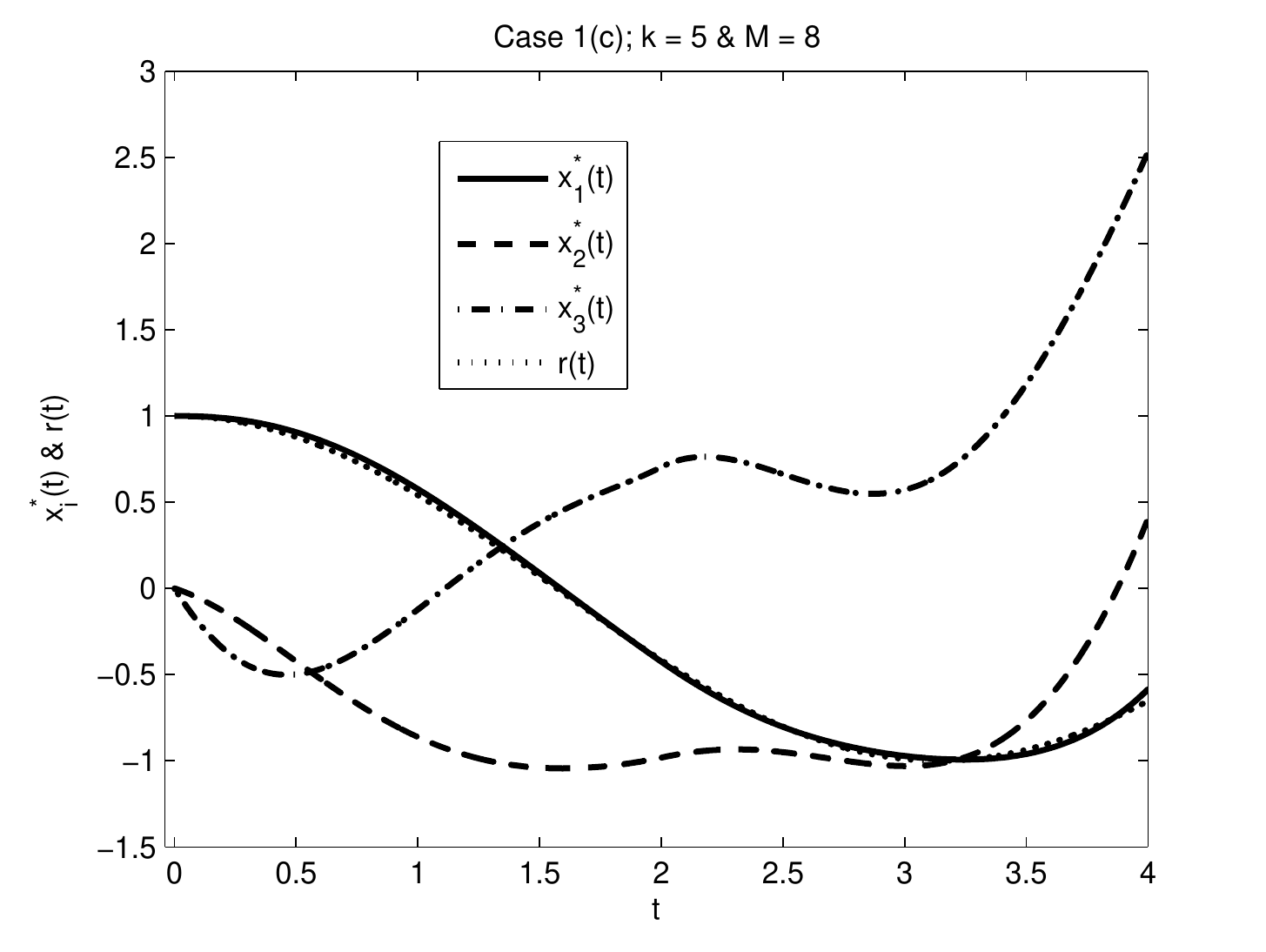}
    \label{Fig:e5.x1c}
}
\subfigure[Optimal control]{
    \includegraphics[scale=.5]{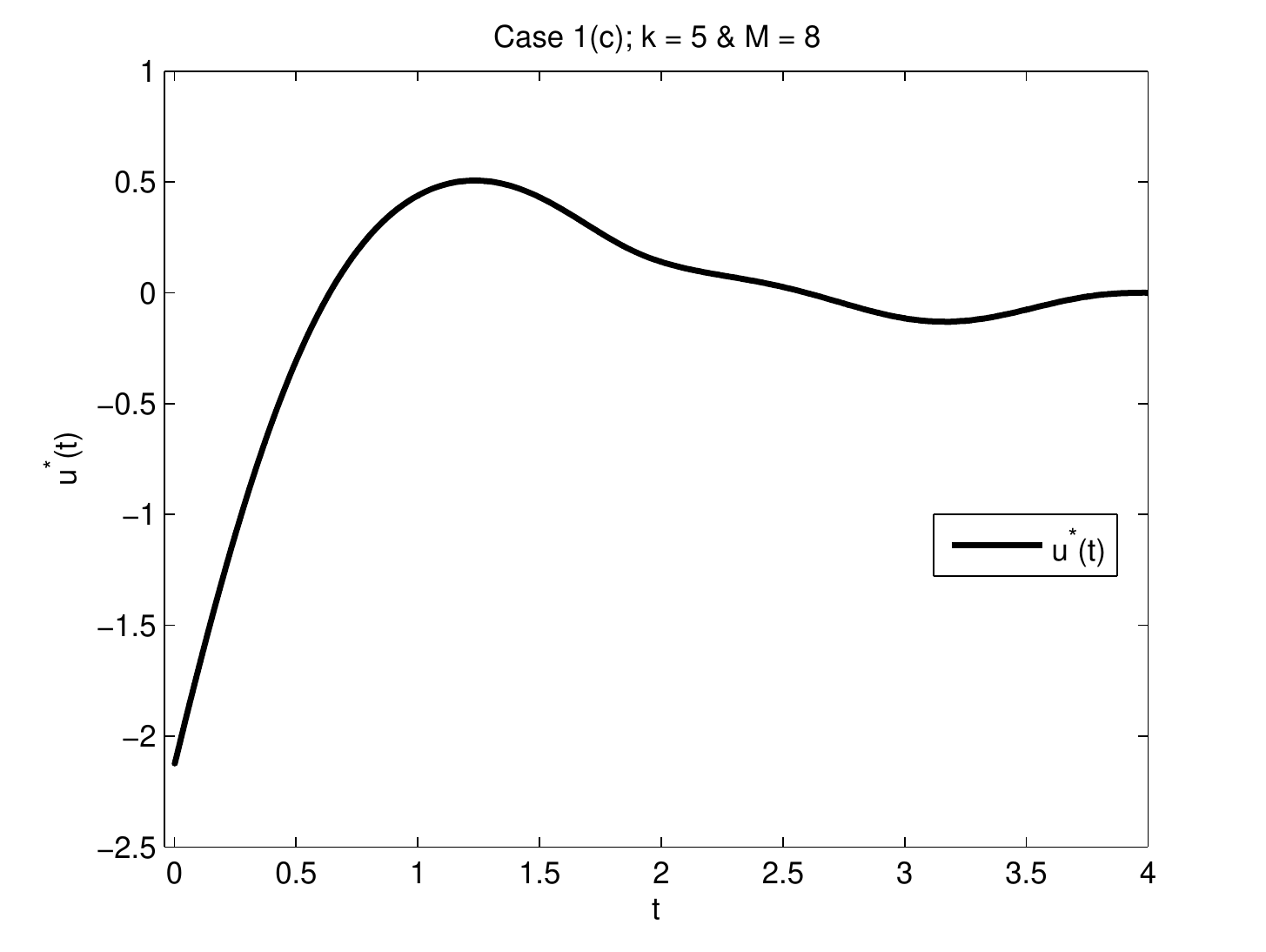}
    \label{Fig:e5.u1c}
}
\caption[Optional caption for list of figures]{Optimal states and control for Example 6, case 1(c).}
\label{Fig:e5.c}
\end{figure}

\begin{figure}[!ht]
\centering
\subfigure[$u^{*}(t)$, case 2]{
    \includegraphics[scale=.5]{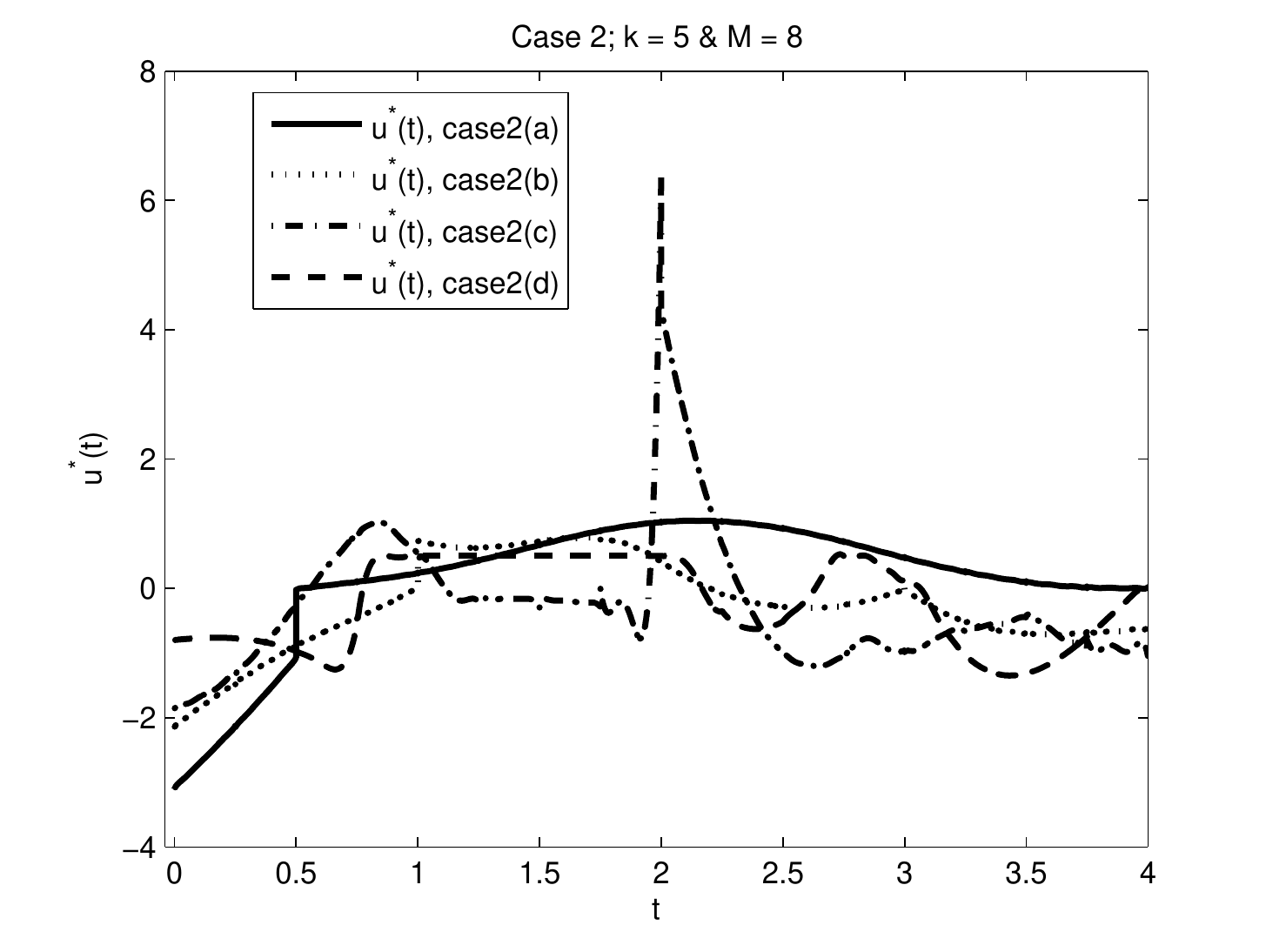}
    \label{Fig:e5.u2}
}
\subfigure[$\text{e}(t)$, case 2]{
    \includegraphics[scale=.5]{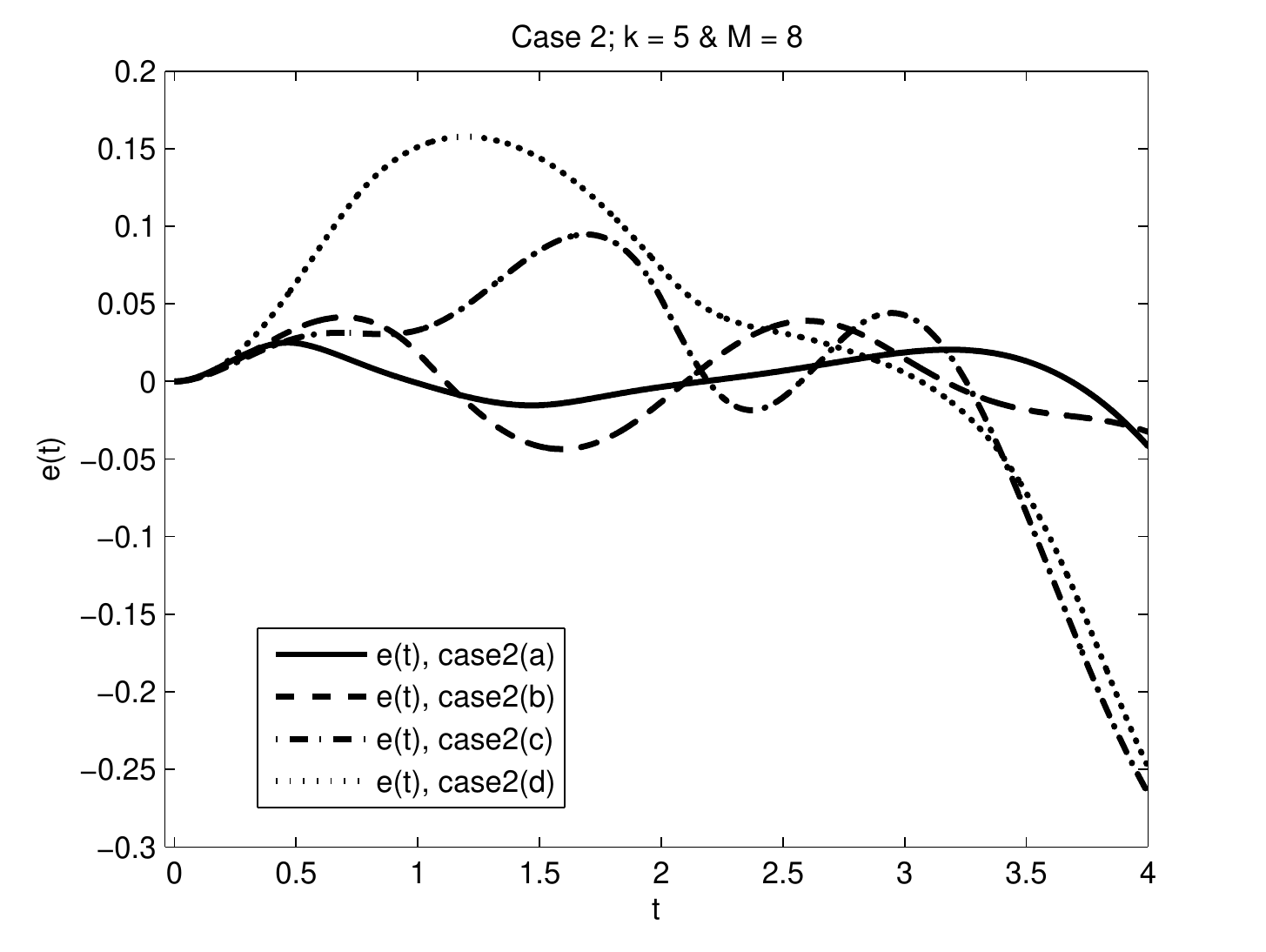}
    \label{Fig:e5.error}
}
\caption[Optional caption for list of figures]{Optimal controls and errors for Example 6, case 2.}
\label{Fig:e5.2abcd}
\end{figure}

By choosing $k=5$ and $M=8$, we give the simulation curves of cases 1(c) and 2 in Figs.\ref{Fig:e5.c} and \ref{Fig:e5.2abcd}. In all four cases (Case 2), the given constraints are satisfied and this shows the efficiency and applicability of the proposed method. The comparison which made with the performance indices in each case above are reported in Table \ref{tab:3}. As we see in Fig.\ref{Fig:e5.error}, imposing the given constraints can affect on the system error $\text{e}(t)$, $\text{e}(t)=x_{1}^*(t)-r(t)$. Setting $g(t)=0.0625t^{2}x_{2}^*(t)+(-0.05t+1)x_{3}^*(t)-u^*(t)$, then Table \ref{tab:4} explains how much the obtained results satisfy the first inequality constraint in case 2(d).

\begin{figure}[!ht]
\centering
\subfigure[$\mathbf{x}^{*}(t)$ and $r(t)$, case 2(a)]{
    \includegraphics[scale=.5]{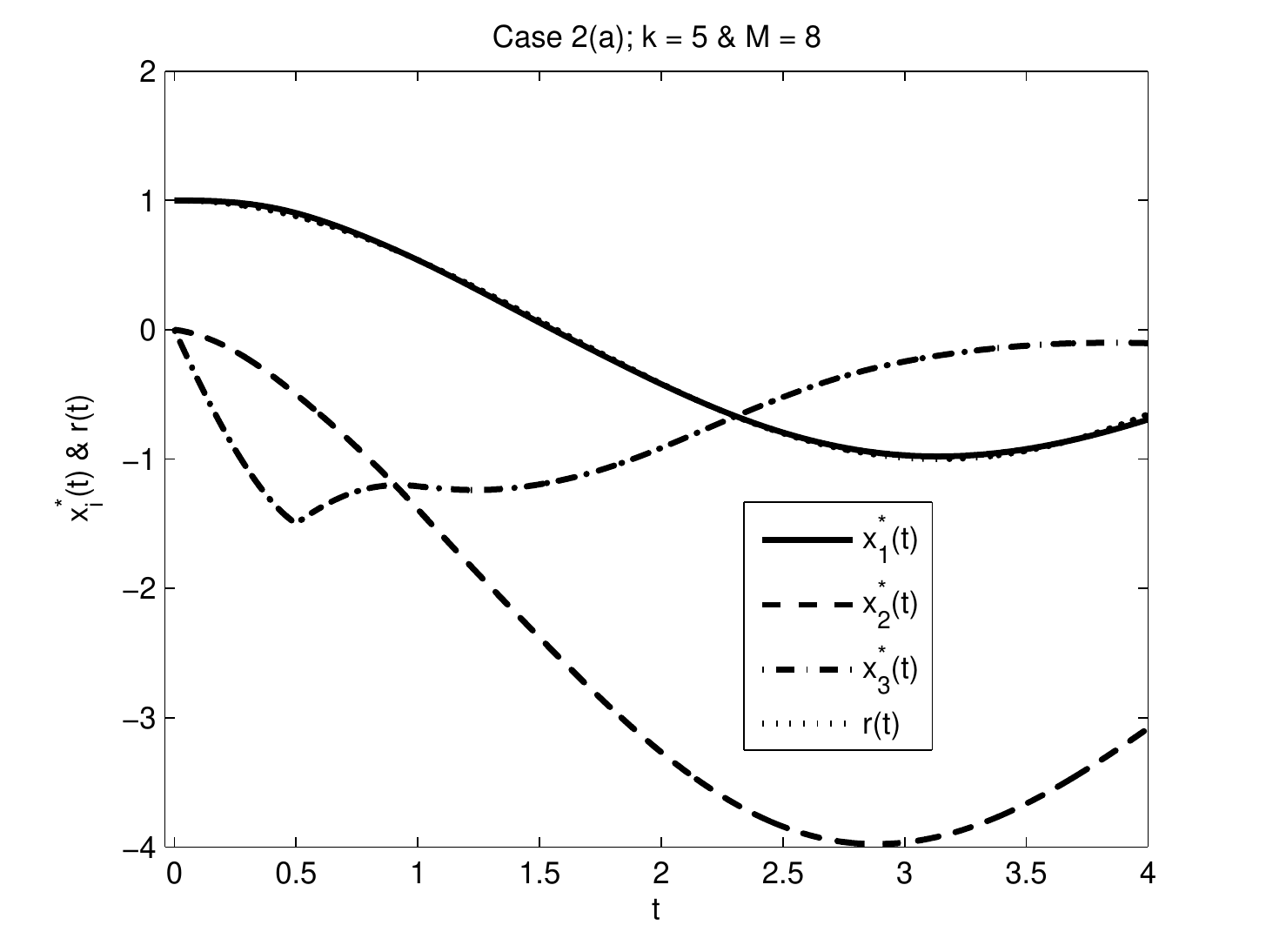}
    \label{Fig:e5.x1a}
}
\subfigure[$\mathbf{x}^{*}(t)$ and $r(t)$, case 2(b)]{
    \includegraphics[scale=.5]{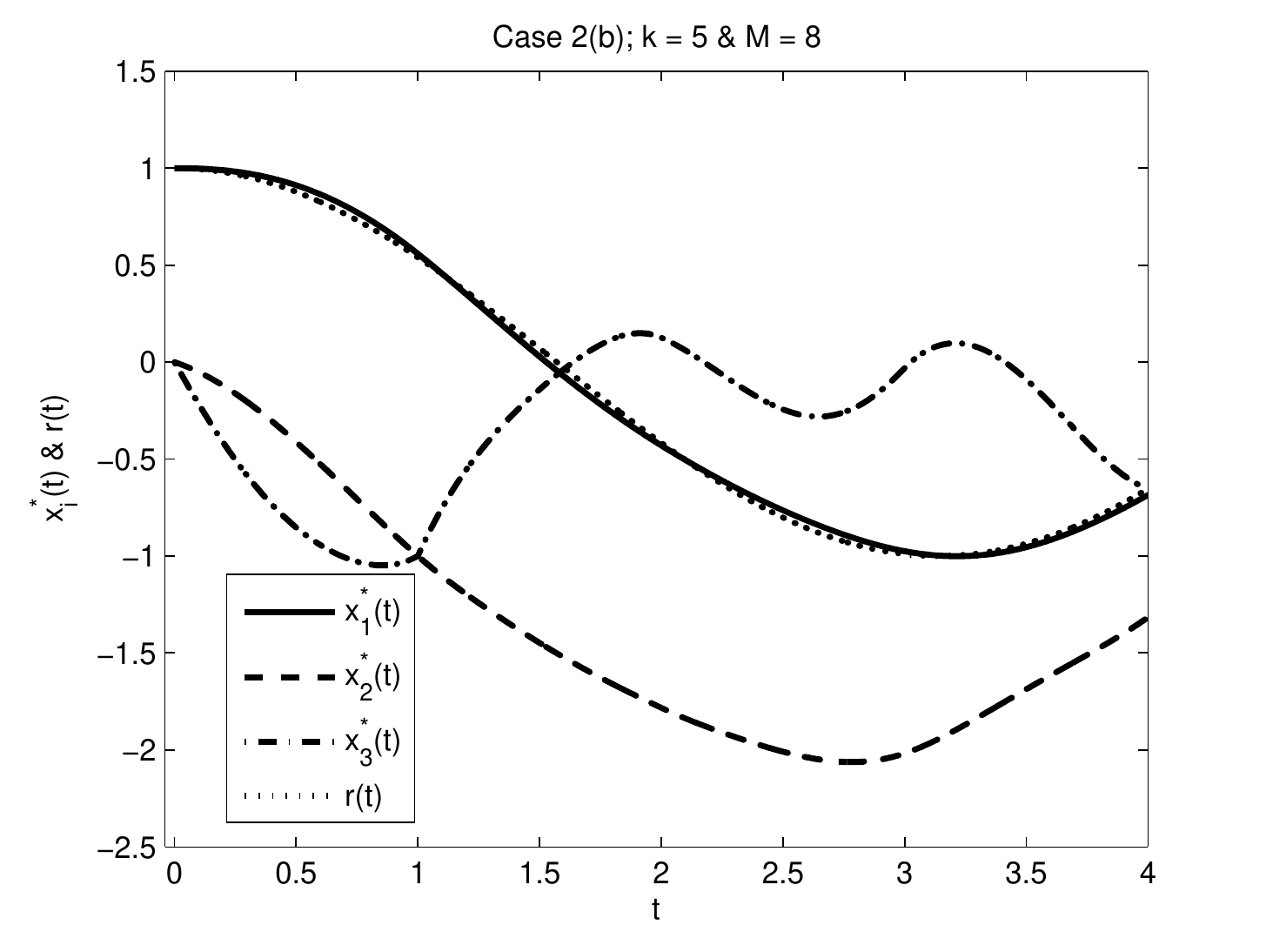}
    \label{Fig:e5.x2b}
}
\caption[Optional caption for list of figures]{Optimal states for Example 6, case 2(a) and 2(b).}
\label{Fig:e5.2ab}
\end{figure}
\begin{figure}[!ht]
\centering
\subfigure[$\mathbf{x}^{*}(t)$ and $r(t)$, case 2(c)]{
    \includegraphics[scale=.5]{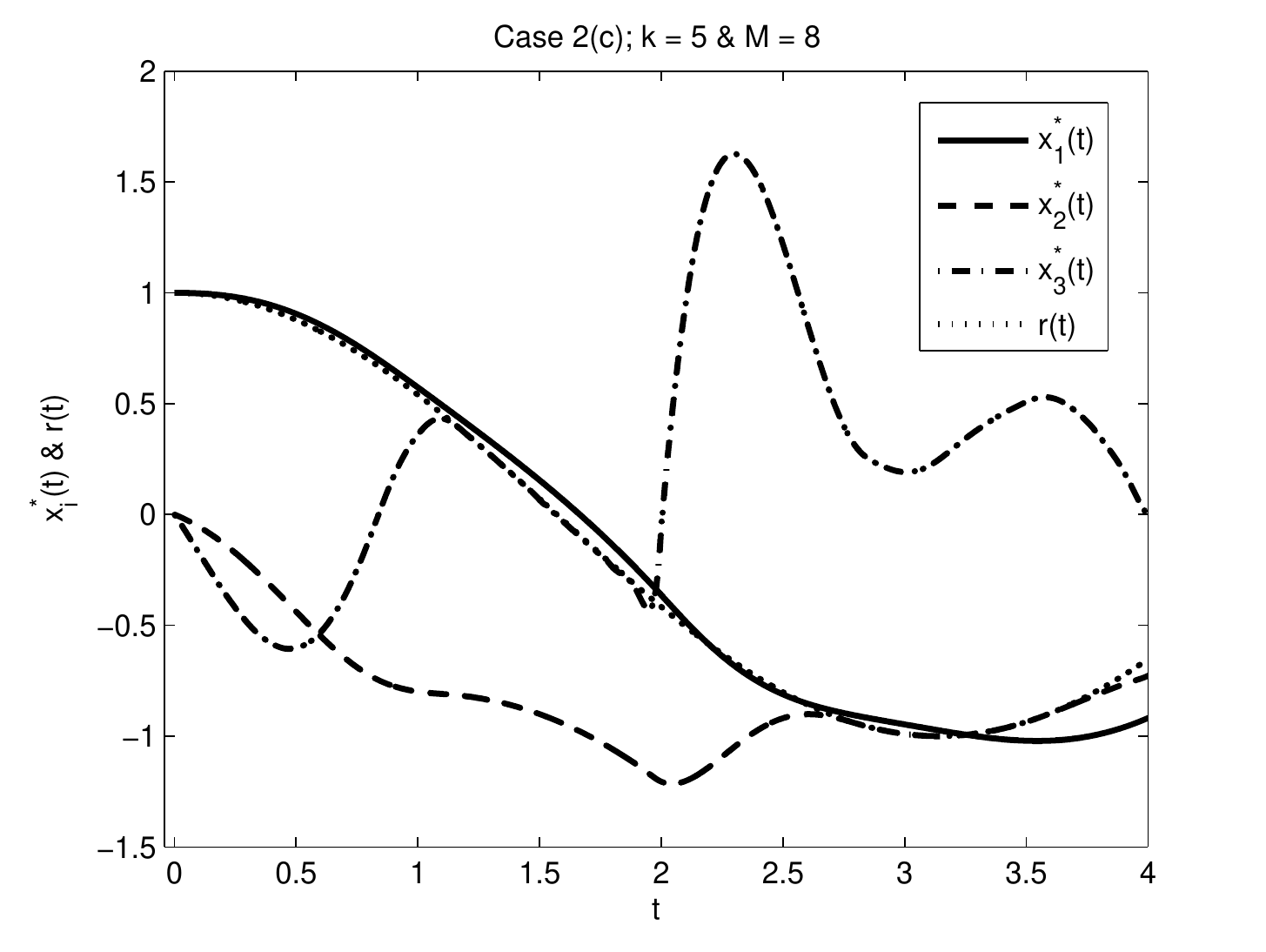}
    \label{Fig:e5.x2c}
}
\subfigure[$\mathbf{x}^{*}(t)$ and $r(t)$, case 2(d)]{
    \includegraphics[scale=.5]{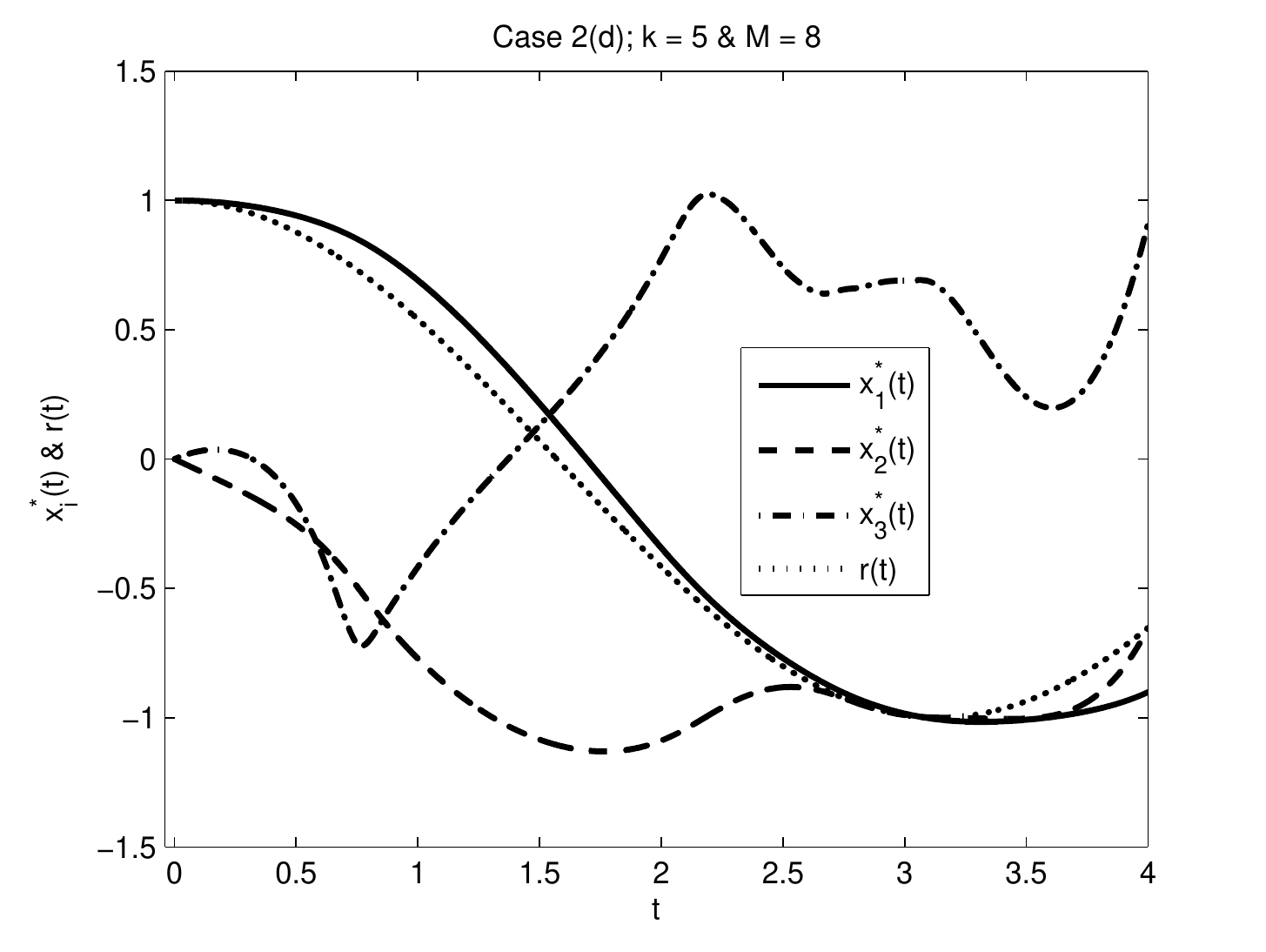}
    \label{Fig:e5.x2c}
}
\caption[Optional caption for list of figures]{Optimal states for Example 6, case 2(c) and 2(d).}
\label{Fig:e5.2cd}
\end{figure}
\begin{table}[h!]
\centering
\caption{ $J^*$ for Example 6}\label{tab:3}
\begin{tabular}{p{1.5in}  p{.57in}}
\toprule
Case & $J^{*}$ \\ 
\midrule
1(a) & 1.804925 \\
1(b) & 0.887031 \\
1(c) & 0.592368 \\
2(a) & 1.909284 \\
2(b) & 1.235810 \\ 
2(c) & 3.548268 \\
2(d) & 3.101320 \\
\bottomrule 
\end{tabular}
\end{table}
\begin{table}[h!]
\centering
\caption{Numerical results for Example 6; case 2(d)}\label{tab:4}
\begin{tabular}{p{.4in}  p{.8in} p{.8in} p{.8in} p{.55in}}
\toprule
t & $\hspace{1mm}x_{2}^{*}(t)$ & $\hspace{1mm}x_{3}^{*}(t)$& $\hspace{1mm}u^{*}(t)$ & \hspace{1mm}$g(t)$\\
\midrule
0 & \hspace{1.2mm}0.00000 & \hspace{1.2mm}0.00000 & -0.80000 & \hspace{1.2mm}0.80000 \\
0.25 & -0.11277 & \hspace{1.2mm}0.02668 & -0.77409 & \hspace{1.2mm}0.80000 \\
0.5 & -0.25283 & -0.17273 & -0.97235 & \hspace{1.2mm}0.79999 \\
0.75 & -0.49264 & -0.71443 & -0.55319 & -0.15177 \\
1 & -0.77177 & -0.41751 & \hspace{1.2mm}0.50000 & -0.94487 \\
1.25 & -0.96555 & -0.12408 & \hspace{1.2mm}0.50000 & -0.71062 \\
1.5 & -1.08493 & \hspace{1.2mm}0.13014 & \hspace{1.2mm}0.50000 & -0.53219 \\
1.75 & -1.13075 & \hspace{1.2mm}0.40720 & \hspace{1.2mm}0.50000 & -0.34486 \\
2 & -1.08918 & \hspace{1.2mm}0.77424 & \hspace{1.2mm}0.50000 & -0.07547 \\
2.25 & -0.96136 & \hspace{1.2mm}1.00549 & -0.35067 & \hspace{1.2mm}0.93886 \\
\bottomrule
\end{tabular}
\end{table}

\section{Conclusion}
An alternative method is introduced to find the optimal control, state and performance index of linear time-varying tracking systems with multiple state and input delays. In the proposed procedure, we can easily change the weighting matrices and impose the combined constraints. When we increase the value of the error weighted matrix, then the output is able to track the reference input better with lower output error, but we have to pay higher cost for larger control effort of the designed system. As can be seen, to better tracking we must try various values of the control weighted matrix. The significant disadvantage of this approach lies in the concept of Chebyshev wavelet, since its definition is slightly less sensitive to changes in time-delays. It should be noted that the method has the ability to implement by Legendre wavelets. The new optimal tracker presented by this paper can be successfully applied to the tracking system regardless of the system stability, minimum phase properties, the dimension of the system, equal number of input and output, the number of delays, and the types of desired states and initial functions.

\renewcommand\refname{REFERENCES}
 

\begin{thebibliography}{99}
{\small

\bibitem{Bellman.Cooke}
R. Bellman, K.L. Cooke
\emph{Differential-difference equations},
RAND Corporation, 1963.

\bibitem{Malek-Zavarei.Jamshidi}
M. Malek-Zavarei, M. Jamshidi,
\emph{Time-Delay Systems: Analysis, Optimization and Applications},
North-Holland, 1978.

\bibitem{Gorecki}
H. G{\'o}recki,
\textit{Analysis and synthesis of time delay systems},
John Wiley \& Sons Inc, 1989.

\bibitem{Naidu}
Desineni Subbaram Naidu,
\textit{OPTIMAL CONTROL SYSTEMS},
Idaho State University. Pocatello. Idaho. USA, CRC PRESS, 2003.

\bibitem{Daubechies}
I. Daubechies,
\textit{Ten Lectures on Wavelets},
SIAM, Philadelphia, 1992.

\bibitem{Mason.Handscomb}
J.C. Mason, D.C. Handscomb,
\textit{Chebyshev polynomials},
CRC Press, 2002.

\bibitem{Gollmann.Maurer}
L. G{\"o}llmann and H. Maurer, \textit{Theory and applications of optimal control problems with multiple time-delays},
Journal of Industrial and Management Optimization. 2014; 10: 413--441.

\bibitem{Brewer}
\newblock J. W. Brewer,
\newblock \emph{Kronecker product and matrix calculus in system theory},
\newblock IEEE Transactions on Circuits and Systems. 1978; 25(9): 772--781.

\bibitem{iman}
I. Malmir,
\emph{Optimal control of linear time-varying systems with state and input delays by Chebyshev wavelets},
Statistics, Optimization \& Information Computing. 2017; 5(4): 302--324.

\bibitem{Liu.Lin}
N. Liu, En-Bing Lin,
\textit{Legendre wavelet method for numerical solutions of partial differential equations},
Numerical Methods for Partial Differential Equations. 2010; 26(1): 81--94.

\bibitem{Gallina.Trevisani}
P. Gallina, A. Trevisani,
\textit{Delayed reference control of a two-mass elastic system},
Journal of Vibration and Control. 2004; 10(1): 135--159.

\bibitem{Liao.Tang.Wang}
F. Liao, Y.Y. Tang, H. Liu, Y. Wang,
\textit{Design of an optimal preview controller for continuous-time systems},
International Journal of Wavelets, Multiresolution and Information Processing. 2011; 9(4): 655--673.

\bibitem{Goodwin.Seron.Dona}
G. Goodwin, M.M. Seron, J.A. De Don{\'a},
\emph{Constrained control and estimation: an optimisation approach},
Springer Science \& Business Media. 2006.

\bibitem{Leondes.Shieh} 
C.T. Leondes, E. Shieh,
\textit{Suboptimal control of linear tracking systems with time delays},
International Journal of Control. 1984; 39(1): 173--180.

\bibitem{Tsai.Tsai.Guo.Chen} 
T.J. Tsai, J.S. H. Tsai, S. Guo, G. Chen,
\textit{OBSERVER-BASED OPTIMAL/SUB-OPTIMAL DIGITAL TRACKERS FOR ANALOG NEUTRAL SYSTEMS WITH MULTIPLE DISCRETE AND DISTRIBUTED TIME DELAYS},
Dynamics of Continuous, Discrete and Impulsive Systems, Series B: Applications \& Algorithms. 2006; 13: 743-789

\bibitem{Tang.Sun.Liu} 
G. Tang, H. Sun, Y. Liu,
\textit{OPTIMAL TRACKING CONTROL FOR DISCRETE TIME-DELAY SYSTEMS WITH PERSISTENT DISTURBANCES},
Asian Journal of Control. 2006; 8(2): 135-140

\bibitem{Tang.Sun} 
G. Tang, H. Sun,
\textit{Optimal tracking control for large-scale interconnected systems with time-delays},
Computers and Mathematics with Applications. 2007; 53: 80–88

\bibitem{Tang.Sun.Pang} 
G. Tang, H. Sun, H. Pang,
\textit{Approximately optimal tracking control for discrete time-delay systems with disturbances},
Progress in Natural Science. 2008; 18: 225–231

\bibitem{Zhang.Tang.Han} 
C. Zhang, G. Tang, S. Han,
\textit{Approximate design of optimal tracking controller for systems with delayed state and control},
IEEE International Conference on Control and Automation, Christchurch, New Zealand, December 9-11, 2009.

\bibitem{Chang.Shieh.Liu.Cofie} 
Y.P. Chang, L.S. Shieh, C.R. Liu, P. Cofie,
\textit{Digital Modeling and PID Controller Design for MIMO Analog Systems with Multiple Delays in States, Inputs and Outputs},
Circuits Syst Signal Process. 2009; 28: 111–145

\bibitem{Tang.Li.Zhao}
G. Tang, C. Li, Y. Zhao,
\textit{Approximate design of optimal tracking controller for time-delay systems},
Chinese Science Bulletin. 2006; 51(17): 2158--2163.

\bibitem{Huang.Tsai.Provence.Shieh}
C.M. Huang, J.S.H. Tsai, R.S. Provence, L.S. Shieh,
\textit{The observer-based linear quadratic sub-optimal digital tracker for analog systems with input and state delays},
Optim. Control Appl. Meth. 2003; 24: 197–236.

\bibitem{Tsai.Wu.Lee.Guo.Su}
J.S.H. Tsai, C.Y. Wu , C.H. Lee, S.M. Guo, T. J. Su,
\textit{A new optimal linear quadratic observer-based tracker under input constraint for the unknown system with a direct feed-through term},
Optimal Control Applications and Methods. 2016; 37: 34--71.

}
\end{thebibliography}
\end{document}